\documentclass[12pt]{article}
\usepackage[centertags]{amsmath}
\usepackage{amsfonts}
\usepackage{amssymb}
\usepackage{amsthm}
\usepackage{newlfont}
\usepackage{graphicx}
\usepackage{t1enc}
\usepackage[latin1]{inputenc}
\usepackage[english]{babel}
 \let\TEXTsymbol\ensuremath
\begin{document}

\begin{center}
\rule{0in}{0.3in}{\Large ARISTOTELIAN ASSERTORIC SYLLOGISTIC}\footnote{%
}

\rule{0in}{0.5in}\rule{0in}{0.5in}\rule{0in}{0.5in}\rule{0in}{0.5in}\rule%
{0in}{0.5in}\rule{0in}{0.5in}\rule{0in}{0.5in}\\[0pt]
{\large MOHAMED A. AMER}\\[0pt]

\rule{0in}{0.5in}{\large \rule{0in}{0.5in}\rule{0in}{0.5in}\rule{0in}{0.5in}%
\rule{0in}{0.5in}\rule{0in}{0.5in}\rule{0in}{0.5in}To Raouf Doss Who
introduced modern logic to Egypt}

\ 
\end{center}

\begin{quote}
$\mathbf{Abstract.}$ Aristotelian assertoric syllogistic, which is currently
of growing interest, has attracted the attention of the founders of modern
logic, who approached it in several (semantical and syntactical) ways.
Further approaches were introduced later on. These approaches (with few
exceptions) are here discussed, developed and interrelated.

\quad Among other things, different facets of soundness, completeness,
decidability and independence are investigated. Specifically
arith\allowbreak metization (Leibniz), algebraization (Leibniz and Boole),
and Venn models (Euler and Venn) are closely examined. All proofs are
simple. In particular there is no recourse to maximal nor minimal conditions
(with only one, dispensable, exception), which makes the long awaited
deciphering of the enigmatic Leibniz characteristic numbers possible. The
problem was how to look at matters from the right perspective.\newline
\end{quote}

\textbf{Introduction. }Aristotelian assertoric syllogistic (henceforth AAS),
which is currently of growing interest (Glashoff 2005), has attracted the
attention of the founders of modern logic. Leibniz, Boole, De Morgan, Venn,
Peirce, Frege, Hilbert, Russell and G\"{o}del all dealt with it. For some of
them it was the starting point (cf. Boole 1948).

Modern treatment of AAS started closer to what may be currently called the
semantical or model theoretic approach. This was threefold: arthimetical,
algebraic, and diagramatic (or set theoretic). The first trend was developed
by Leibniz (\L ukasiewicz 1998, pp. 126-9; Kneale and Kneale 1966, pp.
337-8; Glashoff 2002; and Sotirov 2015). The second was developed by Leibniz
(Kneale and Kneale 1966, pp. 338-45; Lenzen 2004) and after about two
centuries was again developed by Boole (1948), without mentioning the work
of Leibniz. The last trend was developed by Euler, then by Venn (Venn 1880).

With the rise of proof theory late nineteenth century, six syntactical
formalizations of AAS were developed:
\begin{description}
\item  (i) Monadic first order formalization which goes back to Frege (1967,
on p.28 the square of logical opposition may be found). This formalization
is adopted by Hilbert and Ackermann (1950, pp. 44-54).

\item  (ii) Sentential formalization which goes back to Peirce (Bellucci and
Pietarinen 2016, p. 226) and is adopted by G\"{o}del (Adzic and Dosen 2016,
p. 479). The most elaborate study of this formalization is that of \L %
ukasiewicz (1998).

\item  (iii) Dyadic first order formalization which goes back to Shepherdson
(1956). The novel idea of regarding categorical sentences (or propositions)
as binary relational sentences (or propositions) is due to De Morgan
(Valencia 2004, pp. 506-7).
\end{description}

It is worthwile to note here that, according to Boche\'{n}ski (1968, pp.
68-70), Aristotle dealt with the logic of relations among other topics which
Boche\'{n}ski (1968) puts (p. 63) collectively under the title
``Non-analytical laws and rules'', to be distinguished from syllogisms such
as those considered in this article, which Boche\'{n}ski (1968) terms (p.
42) ``analytic''.

\begin{description}
\item  (iv) Natural deduction formalization which goes back to Corcoran
(1972) and Smiley (1973).

\item  (v) First order many-sorted formalization which goes back to Smiley
(1962).

\item  (vi) A recent formalization based on Hilbert's epsilon and tau
quantifiers (Pasquali and Retor\'{e} 2016).
\end{description}

All of the above will be considered below with only two exceptions. The
first is the many-sorted formalization ((v) above), for it is a variant of
the monadic first order formalization mentioned in (i) above; moreover it
was, apparently, abandoned even by its own author (cf. Smiley 1973). The
formalization based on Hilbert's epsilon and tau quantifiers ((vi) above)
will take us far off the current mainstream of logic. So it will be the
second exception and will not be further considered here, though it may have
intrinsic merit especially for those who are interested in formalizing
natural languages.

With only one exception, the modern syntactical formalizations of AAS
degraded it to the rank of a secondary logic, a subordinate or a subsidiary
sublogic of a superior fundamental or principal primary logic. In contrast,
the natural deduction formalization ((iv) above, cf. Boche\'{n}ski (1968,
pp. 3, 31, 42, 49, 52)) rehabilitates it to a full-fledged primary logic, as
was probably designed by its founder: Aristotle, and as was taken for
granted for over two millennia. Accordingly, this formalization will be the
focus of this article. Through completeness we shall see that as far as the
basic sentences (to be defined in 1.6 below) are concerned other
formalizations add nothing new.

In the sequel we deal -from modern standpoints- with AAS, not with medieval
nor traditional syllogistic. In contrast to Boole (1948), Glashoff (2007),
Hilbert and Ackermann (1950), Russinoff (1999), Shepherdson (1956) and
Sotirov (1999), term negation (or complementation, to use a modern term; cf.
Boche\'{n}ski (1968, p. 50)) is not here permitted. Also, in contrast to
Hilbert and Ackermann (1950), \L ukasiewicz (1998), Shumann (2006),
Shepherdson (1956) and Sotirov (1999), Boolean combinations of categorical
sentences are not here permitted. So (the extensional aspect of) AAS will be
just the logic of, or the fragment of set theory which deals with, inclusion
(universal affirmative sentences) and exclusion (universal negative 
sentences) and their contradictories (particular sentences). However, some
exceptions may be appropriate as will be clear, or clarified, at the proper
places.

Among other things, different facets of soundness, completeness,
decidability, and independence are investigated. Particularly
arithmetization (Leibniz), algebraization (Leibniz and Boole) and Venn
models (Euler and Venn) are closely examined.

All proofs given here are simple. In contrast to Corcoran (1972), Glashoff
(2010), Martin (1997), Shepherdson (1956), Smiley (1973) and Smith (1983),
our proofs have no recourse to maximal nor minimal principles nor conditions
(with only one exception, which is indirect and may be dispensed with). This
makes the long awaited deciphering of the enigmatic Leibniz characteristic
numbers possible. The problem was how to look at matters from the right
perspective.

To specify, in section 3 below we provide a polynomial time algorithm to
decide for any finite set of categorical sentences whether it is consistent
and, if it is, to assign a Leibniz model (to be defined below) to it. This
settles positively problem 2 of Glashoff (2002) for finite sets. The general
case is discussed in section 2.

I hope that the simplicity of this exposition of AAS will help to
re-incorporate it into the mainstream of mathematical logic.

After this introduction, the structure of the rest of the article is as
follows:\newline

1. Formalizations of AAS

2. Semantics of AAS

3. Decidability

4. Basic equivalence of the four formalizations

5. Venn soundness and completeness

6. Direct way to Venn models

7. Variations on $NF(C)$

8. Direct completion of direct deduction

9. Models of $NF(C)$ revisited

10. Decidability revisited

11. Sorites

12. Independence

13. Algebraic semantics of AAS, a prelude

14. Algebraic interpretation of $NF(C)$

15. Annihilators: Embedding the partial into a total

16. Back to algebraic interpretation

17. Leibniz and Boole

18. Inadequacy: bounds of AAS

Acknowledgements

Appendix\newline

\textbf{1. Formalizations of AAS. }Formalizations of AAS differ with regard
to permitting the subject and the predicate of a formal symbolic categorical
sentence (henceforth categorical sentence) to be the same. Smith (1983) and
Glashoff (2010) follow Corcoran (1972) in not permitting sameness; as
accommodating sameness ``would entail rather more deviation from the
Aristotelian text'' says Corcoran (1972, p. 696).

On the other hand Smiley (1973) left the door open for permitting sameness,
noting (p. 144) that ``the variables he [Aristotle] uses for the major,
middle and minor terms are all distinct from one another [...]; though when
it comes to substituting actual terms in the resulting forms we are of
course at liberty to replace different variables by the same term (64a1).''.
Consequently, it seems that Aristotle excluded sameness, for technical -not
philosophical- reasons. This is, possibly, why \L ukasiewicz (1998) adopted
sameness (see pp. 77, 88); while Martin (1997) simultaneously considered two
systems, one of them is permitting sameness and the other is not.

In conformity with the current mainstream of mathematical logic, sameness is
here permitted. Excluding sameness, and other variations, will be considered
in section 7 below.\newline

\textbf{1.1. Monadic first order formalization of AAS.} The languge here is
a standard first order language, with or without equality, whose set \textbf{%
P} of non-logical constants has at least three elements, and all of its
elements are unary relational symbols. In the sequel ``$P$'', ``$Q$'' and ``$R$%
'' will be metalinguistic variables ranging over the elements of \textbf{P}.
With abuse of notation, ``\textbf{P}'' will denote this language too.\newline

\noindent ABBREVIATIONS 1.1.

``$APQ$'' \quad is an abbreviation for \qquad ``$\forall x(Px\rightarrow Qx)$%
''

``$EPQ$'' \quad is an abbreviation for \qquad ``$\forall x(Px\rightarrow
\urcorner Qx)$''

``$IPQ$'' \quad ~is an abbreviation for \qquad ``$\exists x(Px\wedge Qx)$''

``$OPQ$'' \quad is an abbreviation for \qquad ``$\exists x(Px\wedge
\urcorner Qx)$''\newline
\newline
\textit{\noindent \noindent }DEFINITION 1.2. $MF$(\textbf{P}) is the theory
based on \textbf{P} with only one non -logical axiom schema, namely, $%
APQ\rightarrow IQP$ (which is equivalent to the schema $\exists x$ $Px$).%
\newline

\noindent PROPOSITION 1.3. The following are theorem schemata of $MF($%
\textbf{P}$)$:

\begin{description}
\item  1. $EPQ\leftrightarrow$ $\urcorner$ $IPQ$ \qquad \qquad 2. $%
OPQ\leftrightarrow \urcorner APQ$

\item  3. $APP$ \qquad \qquad \qquad \qquad 4. $APQ\rightarrow IQP$

\item  5. $EPQ\rightarrow EQP$ \qquad \qquad 6. $APQ\wedge AQR\rightarrow APR
$

\item  7. $APQ\wedge EQR\rightarrow EPR$ \newline
\end{description}

\noindent \textit{Proof.} Routine. \hspace{10cm} $\square$ 
\newline

\textbf{1.2. Sentential formalization of AAS.} The symbols ``$A$'', ``$E$'',
``$I$'' and ``$O$'' were made use of in section 1.1, in this section they
will be made use of differently. This abuse of notation is benign as long as
the intended denotation is clear from the context, so it will be here
permitted. Such abuses of notation may be permitted later on without further
notice.

Let $J$ be a set (whose elements are to correspond to categorical constants)
having at least three elements, let $A,E,I$ and $O$ be four injective
functions of pairwise disjoint ranges, each of domain $J\times J$, and let $%
AS(J)$ be the union of their ranges. In the sequel ``$i$'', ``$j$'' and ``$k$%
'' will be metalinguistic variables ranging over the elements of $J$.

The language here is a standard sentential language whose set of sentential
symbols is $AS(J)$. With abuse of notation ``$J$'' will denote this language
too.\newline

\noindent DEFINITION 1.4. $SF(J)$ is the theory based on $J$ with the
following non - logical axiom schemata:

1. \textit{Ei}$j\leftrightarrow \urcorner Iij\qquad \qquad $\qquad 2. $%
Oij\leftrightarrow \urcorner Aij$

3. $Aii$ $\qquad \qquad \qquad \quad $\qquad 4. $Aij\rightarrow Iji$

5. \textit{Ei}$j\rightarrow Eji\qquad \qquad $\qquad 6. $Aij\wedge
Ajk\rightarrow Aik$

7. $Aij\wedge Ejk\rightarrow $ \textit{Ei}$k$\newline

The proof machinery is modus ponens together with any standard set of
sentential logical axiom schemata.\newline

\noindent REMARK 1.5. There are two kinds of substitution: sentences for
sentences and indices for indices. Each may be permitted, under some
conditions, as a derived rule of inference (cf. \L ukasiewiez 1998, p. 88;
see section 1.5 below).\newline

\textbf{1.3. Dyadic first order formalization of AAS. }The language here is
a standard first order language, with or without, equality whose non-logical
constants are four binary relation symbols ``$A$'',``$E$'',``$I$'' and ``$O$'', together with a
set $C$ of individual (or categorical) constants having at least three
elements. In the sequel ``$a$'', ``$b$'', and ``$c$'' will be metalinguistic
variables ranging over the elements of $C$. With abuse of notation ``$C$''
will denote this language too.\newline

\noindent DEFINITION 1.6. $DF(C)$ is the theory based on $C$ whose
non-logical axioms are the universal closures of:\newline

1. $Exy\leftrightarrow \urcorner Ixy\qquad \qquad $\qquad 2. $%
Oxy\leftrightarrow \urcorner Axy$

3. $Axx\qquad \qquad \qquad \quad ~$\qquad 4. $Axy\rightarrow Iyx$

5. $Exy\rightarrow Eyx\qquad \qquad $\qquad 6. $Axy\wedge Ayz\rightarrow Axz$

7. $Axy\wedge Eyz\rightarrow Exz$ \newline

\textbf{1.4. Natural deduction formalization of AAS.} The language here is a
sublanguage of the language defined in 1.3. The alphabet is the four binary
relation symbols $A,E,I$ and $O$, together with a set $C$ of individual (or
categorical) constants having at least three elements. The sentences are the
equality free atomic sentences of 1.3., viz. a sentence is a string $Yab$
where $Y\in \{A,E,I,O\}$ and $a,b\in C$. By abuse of notation ``$C$'' will
denote this language too, and the set of all sentences will be denoted by ``$%
S(C)$''. In the sequel ``$\alpha $'', ``$\beta $'', ``$\gamma $'', ``$\delta 
$'', ``$\sigma $'' and ``$\rho $'' will be metalinguistic variables ranging
over the elements of $S(C)$.

Sentences starting with $A$ or $E$ are called universal, those starting with $%
I$ or $O$ are called particular. Also, sentences starting with $A$ or $I$
are called affirmative, those starting with $E$ or $O$ are called negative.
For $W\in \{A,E,I,O\}$, sentences starting with $W$ are called $W$-sentences.%
\newline

\noindent DEFINITION 1.7. $NF(C)$ is the logical system based on the
language $C$ with the following deduction rules (or enrichments thereof, see
sections 8 and 11 below):

\begin{description}
\item  0. $\frac {}{Aaa}$ $\quad \quad $~$(A$-$Id)$ \qquad \qquad \qquad 1. $%
\frac{Aab}{Iba}$ $\quad \quad \ (Apc)$

\item  2. $\frac{Eab}{Eba}$ $\quad \quad $\ $(Ec)$ \qquad \qquad \quad
\qquad 3. $\frac{Aab,Abc}{Aac}$ $\quad $(Barbara)

\item  4. $\frac{Aab,Ebc}{Eac}$ $\quad $(Celarent).
\end{description}

``Barbara'' and ``Celarent'' are, respectively, the medieval names of the
rules 3 and 4; ``$A$-$Id"$, ``$Apc$'' and ``$Ec$'' are, respectively,
abbreviations for ``A-identity'', ``A-partial conversion'' and
``E-conversion''. For simplicity, we may write ``rules'' instead of
``deduction rules''.\newline

\noindent DEFINITION 1.8. A direct deduction (or d-deduction) of $\sigma
(\in S(C))$ from $\Gamma (\subseteq S(C))$ is a sequence $<\rho _{i}>_{i\in
k}(k\in N^{+})$ such that $\rho _{k-1}=\sigma $ and for each $i\in k,\sigma
_{i}\in \Gamma $ or is the consequent of some rule of $NF(C)$ whose
antecedents are previous terms of the sequence. In this case we write $%
\Gamma \stackrel{d}{\vdash }\sigma $, and $\sigma $ is said to be a direct
consequence (or theorem) of $\Gamma $. Also 
$<\rho _{i}>_{i\in
k}$ is said to be a direct (or $d$-) deduction from $\Gamma $. From now on,
the rules 0-4 given above will be called also ``$d$-rules''.

Regarding the current mainstream of mathematical logic, this definition is a
typical definition. In contrast, corresponding definitions given in Corcoran
(1972), Glashoff (2010), Martin (1997), Smiley (1973) and Smith (1983) are
atypical, each has its own peculiarity.

To get closer to the Aristotelian tradition, a more restricted definition of
direct deduction is presented in section 11 below, and its relationship to
the above one is investigated there.

As usual, the contradictory $\widehat{\sigma }$ of $\sigma (\in S(C))$ is
defined as follows:

$\widehat{Aab}=Oab$ \qquad $\widehat{Eab}=Iab$ \qquad $\widehat{Iab}=Eab$
\qquad $\widehat{Oab}=Aab$

so $\widehat{\widehat{\sigma }}=\sigma $. 

A set $\Gamma (\subseteq S(C))$ is
said to be $d$-inconsistent (or $d$-contradictory) if $\Gamma \stackrel{d}{%
\vdash }\sigma $ and $\Gamma \stackrel{d}{\vdash }\widehat{\sigma }$, for
some $\sigma \in S(C)$; otherwise $\Gamma $ is said to be $d$-consistent.%
\newline

\noindent DEFINITION 1.9. The general (or $g$-) deduction relation $%
\stackrel{g}{\vdash }(\subseteq \wp (S(C))\times S(C))$ is defined as
follows:

$\Gamma \stackrel{g}{\vdash }\sigma $ \qquad iff \qquad $\Gamma \cup \{%
\widehat{\sigma }\}$ is $d$-inconsistent.\newline

For $\Gamma \subseteq S(C)$, ``$\Gamma $ is $g$-inconsistent (or $g$%
-contradictory)'' and ``$\Gamma $ is $g$-consistent'' may be defined along
the above lines, replacing ``$d$'' by ``$g$''. Obviously $\stackrel{d}{%
\vdash }$ $\subseteq $ $\stackrel{g}{\vdash }$, so if $\Gamma (\subseteq
S(C))$ is $d$-inconsistent, it is $g$-inconsistent.\newline

For $e\in \{d,g \}$ and $\Gamma \subseteq S(C)$, ``$\Gamma ^{e}$'' will denote
the closure of $\Gamma $ under $\stackrel{e}{\vdash }$, i.e.

$\Gamma ^{e}$= the smallest $\Delta \subseteq S(C)$ such that $\Gamma
\subseteq \Delta $ and for every $\sigma \in S(C)$, $\sigma \in \Delta $
whenever $\Delta \stackrel{e}{\vdash }\sigma $.\newline

\noindent LEMMA 1.10. Let $\Sigma ,\Sigma ^{\prime }$ be subsets of $S(C)$
such that for every $\sigma \in \Sigma ,$ $\Sigma ^{\prime }\stackrel{d}{%
\vdash }\sigma $. For every $d$-deduction $<\rho _{i}>_{i\in k}$ from $\Sigma 
$, there are a $k^{\prime }(\geq k)$, a $d$-deduction $<\rho _{i}^{\prime
}>_{i\in k^{\prime }}$ from $\Sigma ^{\prime }$ and a strictly increasing
function $f:k\rightarrow k^{\prime }$ such that $f(k-1)=k^{\prime }-1$ and
for every $i\in k$, $\rho _{i}=\rho _{f(i)}^{\prime }$. Hence for every $%
\alpha \in S(C)$, $\Sigma ^{\prime }\stackrel{d}{\vdash }\alpha $ whenever $%
\Sigma \stackrel{d}{\vdash }\alpha $.\newline

\noindent \textit{Proof.} By induction on $\ell $, the number of times of
making use in $<\rho _{i}>_{i\in k}$ of assumptions from $\Sigma $.

Basis: $\ell =0$; take $k^{\prime }=k,<\rho _{i}^{\prime }>_{i\in k^{\prime
}}$ $=$ $<\rho _{i}>_{i\in k}$ and $f$ the identity function on $k$.

Induction step: assume the required for $\ell =m$. Let $\ell =m+1$ and let $%
j(\in k)$ be the last line in which an assumption from $\Sigma $ is made use
of. The case $j=0$ is easier than the case $j>0$, so we shall deal only with
the latter.

By the induction hypothesis, there are a $j^{\prime }(\geq j)$, a deduction $%
<\alpha _{i}>_{i\in j^{\prime }}$, from $\Sigma ^{\prime }$ and a strictly
increasing function $g:j\rightarrow j^{\prime }$ such that $g(j-1)=j^{\prime
}-1$ and for every $i\in j,$ $\rho _{i}=\alpha _{g(i)}$.

Let $<\beta _{i}>_{i\in m}$ be a deduction of $\rho _{j}$ from $\Sigma
^{\prime }$. Put:

$k^{\prime }=j^{\prime }+m+k-j-1,$

$\gamma _{i}=\rho _{i+j+1}$ \quad for $i\in k-j-1$,

$<\rho _{i}^{\prime }>_{i\in k^{\prime }}=<\alpha _{i}>_{i\in j^{\prime
}}\frown $ $<\beta _{i}>_{i\in m}\frown $ $<\gamma _{i}>_{i\in k-j-1},$

where ``$\frown $'' is the concatenation operation symbol.

Evidently $k^{\prime }\geq k$. The completion of the proof is now easy. \hspace{2cm}$%
\square $ \newline

Parts 3 and 4 of the next proposition are, respectively, reformulations of
lemmata M$_{1}$ and M$_{2}$ of Corcoran (1972).\newline

\noindent PROPOSITION 1.11. Let $\Gamma \cup \{\sigma \}\subseteq S(C)$, $%
\Gamma ^{\prime }=\{\rho \in S(C):\Gamma \stackrel{g}{\vdash }\rho \}$, $%
U\in \{A,E\},W\in \{I,O\},e\in \{d,g\}$ and $a,b,c\in C$, then:

1. $\Gamma ^{d}=\{\rho \in S(C):\Gamma \stackrel{d}{\vdash }\rho \}$.

2. If $\Gamma \cup \{Wab\}\stackrel{d}{\vdash }\sigma $ and $\sigma \neq Wab$%
, then $\Gamma \stackrel{d}{\vdash }\sigma $.

3. If $\Gamma $ is $d$-consistent, then $\Gamma \stackrel{g}{\vdash }Uab$
iff $\Gamma \stackrel{d}{\vdash }Uab$.

4. $\Gamma $ is $d$-consistent iff it is $g$-consistent.

5. If $\rho \in \Gamma ^{\prime }$ and $\Gamma ,\rho \stackrel{g}{\vdash }%
\sigma $, then $\Gamma \stackrel{g}{\vdash }\sigma $.

6. $\Gamma ^{\prime g}=\Gamma ^{\prime }$, hence $\Gamma ^{g}=\Gamma
^{\prime }$.

7. $\Gamma ^{ee}=\Gamma ^{e}$; hence, for every $d$-rule $r$, if each
antecedent of $r$ belongs to $\Gamma ^{e}$, then so also does its consequent.%
\newline

\newpage

\noindent \textit{Proof.}\newline

1. The required is a corollary of the above lemma.

2. Generalizing upon the metalinguistic variable ``$\sigma $'', the resulting
sentence may be proved by course of values induction on the length of the $d$%
-deduction from $\Gamma \cup \{Wab\}$, noticing that $Wab$ is not a premise
of any rule of $NF(C)$.

3. Let $\Gamma \stackrel{g}{\vdash }Uab$, then for some $\alpha \in S(C)$, $%
\Gamma $, $\widehat{U}ab\stackrel{d}{\vdash }\alpha ,\widehat{\alpha }$. So,
by part 2, if $\Gamma $ is $d$-consistent then $\widehat{U}ab\in \{\alpha ,%
\widehat{\alpha }\}$, then $Uab\in \{\alpha ,\widehat{\alpha }\}$. So, by
part 2 again, $\Gamma \stackrel{d}{\vdash }Uab$. The other direction is
obvious.

4. Let $\Gamma $ be $d$-consistent and $g$-inconsistent, then there is a
universal $\alpha \in S(C)$ such that $\Gamma \stackrel{g}{\vdash }\alpha ,%
\widehat{\alpha }$, then by part 3, $\Gamma \stackrel{d}{\vdash }\alpha $.
Also there is $\beta \in S(C)$ such that $\Gamma ,\alpha \stackrel{d}{\vdash 
}\beta ,\widehat{\beta }$, so, by lemma 1.10., $\Gamma \stackrel{d}{\vdash }%
\beta ,\widehat{\beta }$, hence $\Gamma $ is $d$-inconsistent. Consequently,
if $\Gamma $ is $g$-inconsistent it is $d$-inconsistent. The other direction
is obvious.

5. Obvious if $\Gamma $ is $d$-inconsistent, so let $\Gamma $ be $d$%
-consistent and let $\rho \in \Gamma ^{\prime }$ and $\Gamma ,\rho \stackrel{%
g}{\vdash }\sigma $. There is $\alpha \in S(C)$ such that $\Gamma ,\rho ,%
\widehat{\sigma }\stackrel{d}{\vdash }\alpha ,\widehat{\alpha }$. If $\rho $
is universal, then by part 3 and lemma 1.10, $\Gamma ,\widehat{\sigma }%
\stackrel{d}{\vdash }\alpha ,\widehat{\alpha }$. Also, if $\rho $ is
particular and $\rho \notin \{\alpha ,\widehat{\alpha }\}$, then by part 2, $%
\Gamma ,\widehat{\sigma }\stackrel{d}{\vdash }\alpha ,\widehat{\alpha }$. In
both cases $\Gamma \stackrel{g}{\vdash }\sigma $. In the remaining case $%
\rho $ must be particular and $\Gamma ,\rho ,\widehat{\sigma }\stackrel{d}{%
\vdash }\rho ,\widehat{\rho }$, then by part 2, $\Gamma ,\widehat{\sigma }%
\stackrel{d}{\vdash }\widehat{\rho }$. But there is $\beta \in S(C)$ such
that $\Gamma ,\widehat{\rho }\stackrel{d}{\vdash }\beta ,\widehat{\beta }$,
then $\Gamma ,\widehat{\sigma }\stackrel{d}{\vdash }\beta ,\widehat{\beta }$%
, hence $\Gamma \stackrel{g}{\vdash }\sigma $, which completes the proof.

6. By induction, part 5 may be generalized to: for every finite $\Delta
\subseteq \Gamma ^{\prime }$, $\Gamma \stackrel{g}{\vdash }\sigma $ whenever 
$\Gamma \cup \Delta \stackrel{g}{\vdash }\sigma $. From this it readily
follows that $\Gamma \stackrel{g}{\vdash }\sigma $ whenever $\Gamma ^{\prime
}\stackrel{g}{\vdash }\sigma $, hence the result.

7. By part 1 and lemma 1.10, $\Gamma ^{dd}=\Gamma ^{d}$, and by part 6, $%
\Gamma ^{gg}=\Gamma ^{\prime g}=\Gamma ^{\prime }=\Gamma ^{g}$. To prove the
last clause, let $r$ be a $d$-rule. If each antecedent of $r$ belongs to $%
\Gamma ^{e}$, then its consequent belongs to $\Gamma ^{ee}=\Gamma ^{e}$. \hspace{4cm}$%
\square $\newline

In view of part 4 of the above proposition, for $e\in \{d,g\}$, the prefix ``%
$e$-'' may be deleted from ``$e$-consistent'', ``$e$-inconsistent'' and ``$e$%
-contradictory''.\newline

\newpage

\textbf{1.5. Equality / equivalence.} For $\Gamma \subseteq S(C)$, $e\in
\{d,g\}$ and $a,b\in C$, $\Gamma \stackrel{e}{\vdash }Aab,Aba$ is equivalent
to each of:

1. $\Gamma \stackrel{e}{\vdash }Aac$ \qquad iff \qquad $\Gamma \stackrel{e}{%
\vdash }Abc$ \qquad all $c\in C$,

2. $\Gamma \stackrel{e}{\vdash }Aca$ \qquad iff \qquad $\Gamma \stackrel{e}{%
\vdash }Acb$ \qquad all $c\in C$.

Thus, $\Gamma \stackrel{e}{\vdash }Aab$, $Aba$ imply the substitutability of 
$a,b$ for each other in universal positive sentences. This will be
generalized below to all universal sentences, respectively all sentences,
for $e=d$, respectively $e=g$. The last generalization is the essence of
equality (congruence or equivalence, depending on the situation). It holds,
in the respective appropriate forms, for the other formalizations as is
shown in the following:\newline

\noindent THEOREM 1.12.\newline

1. Let $P,Q\in $\textbf{P}, then:

$MF(\mathbf{P})\vdash (APQ\wedge AQP)\rightarrow (\varphi \leftrightarrow
\varphi ^{\prime })$ for every form $\varphi $ of \textbf{P}, where $\varphi
^{\prime }$ is a form obtained from $\varphi $ by substituting some
occurrences of ``$P(x)$'' in $\varphi $ by ``$Q(x)$'', or vice versa.\newline

2. Let $i,j\in J$, then:

$SF(J)\vdash (Aij\wedge Aji)\rightarrow (\alpha \leftrightarrow \alpha
^{\prime })$ for every sentence $\alpha $ of $J$, where $\alpha ^{\prime }$
is a sentence obtained from $\alpha $ by substituting some occurrences of ``$%
i$'' in $\alpha $ by ``$j$'', or vice versa.\newline

3. Let $a,b\in C$, then:

$DF(C)\vdash (Aab\wedge Aba)\rightarrow (\varphi \leftrightarrow \varphi
^{\prime })$ for every form $\varphi $ of $C$, where $\varphi ^{\prime }$ is
a form obtained from $\varphi $ by substituting some occurrences of ``$a$''
in $\varphi $ by ``$b$'', or vice versa; provided -for languages with
equality- no substitution takes place in a form or a subform of the form $%
t=t^{\prime }$, where $t$ and $t^{\prime }$ are terms.\newline

4. Let $\Gamma \subseteq S(C)$, $e\in \{d,g\}$ and $a,b\in C$, and let $%
\Gamma \stackrel{e}{\vdash }Aab$, $Aba$, then for all $c\in C$:

$\Gamma \stackrel{e}{\vdash }Yac$ iff $\Gamma \stackrel{e}{\vdash }Ybc$ and $%
\Gamma \stackrel{e}{\vdash }Yca$ iff $\Gamma \stackrel{e}{\vdash }Ycb$,
where $Y\in \{A,E\},\{AE,I,O\}$ for $e=d,g$ respectively.\newline

\noindent \textit{\noindent Proof.} The first three parts may be proved by
the standard methods developed in the respective formal systems.

For the last part, part 7 of proposition 1.11 secures the required for $e\in
\{d,g\}$ and $Y\in \{A,E\}$.

It remains to consider the cases where $e=g$ and $Y\in \{I,O\}$. If $\Gamma $
is inconsistent the required follows by the definition of $\stackrel{g}{%
\vdash }$, so let $\Gamma $ be consistent. Assume $\Gamma \stackrel{g}{%
\vdash }Aab$, $Aba$ and $\Gamma \stackrel{g}{\vdash }Ica$, then by part 3 of
proposition 1.11, $\Gamma \stackrel{d}{\vdash }Aab$, $Aba$, hence $\Gamma $, 
$Ecb\stackrel{d}{\vdash }Eca$. But there is $\alpha \in S(C)$ such that $%
\Gamma $, $Eca\stackrel{d}{\vdash }\alpha ,\widehat{\alpha }$, consequently $%
\Gamma $, $Ecb\stackrel{d}{\vdash }\alpha ,\widehat{\alpha }$ and, by
definition, $\Gamma \stackrel{g}{\vdash }Icb$. The other cases are similar
or easier. \hspace{7.4cm}$\square $\newline

\textbf{1.6. Basic sentences.} In each of the four formalizations $MF($%
\textbf{P}$)$, $SF(J)$, $DF(C)$ and $NF(C)$ the sentences to be made use of
in the Aristotelian syllogistic will be called basic (or categorical)
sentences. The sets of basic sentences will be denoted, respectively, by ``$%
BM($\textbf{P}$)$'', ``$BS(J)$'', ``$BD(C)$'' and ``$BN(C)$''. That is:

$BM(\mathbf{P})=\{YPQ:Y\in \{A,E,I,O\}$ and $P,Q\in $\textbf{P}$\mathbf{\}}$.

$BS(J)=AS(J)$ (= the set of all atomic sentences of $J$).

$BD(C)=$ the set of all (equality free) atomic sentences of $C$.

$BN(C)=S(C)$ ($=BD(C)$).\newline

\textbf{1.7. Interpretation.} Let $h:C\rightarrow J$, with abuse of notation
(no confusion will ensue) we define another function $h:BN(C)\rightarrow
BS(J)$ by $h(Yab)=Yhahb$, for $Y\in \{A,E,I,O\}$ and $a,b\in C$. As usual,
for $\Gamma \subseteq BN(C)$, the image of $\Gamma $ under $h$ is denoted by
``$h(\Gamma)$''; also we may write ``$Y^{h}ab$'', ``$\Gamma ^{h}$'' for ``$%
h(Yab)$'', ``$h(\Gamma )$'' respectively. The function $h$ is said to be an
interpretation of $BN(C)$ in $BS(J)$. Similarly $BM($\textbf{P}$)$, $BS(J)$
and $BN(C)$ ($=BD(C)$) may be interpreted in each other.\newline

\noindent PROPOSITION 1.13. Let $\Gamma \cup \{\sigma \}\subseteq BN(C)$,
let $h$ and $H$ be interpretations of $BN(C)$ in $BS(J)$ and $BM($\textbf{P}$%
)$ respectively, and let $<\Gamma ^{\prime },\sigma ^{\prime },\widehat{%
\sigma }^{\prime },T>\in \{<\Gamma ,\sigma ,\widehat{\sigma },DF(C)>$, $%
<\Gamma ^{h},\sigma ^{h},\widehat{\sigma }^{h},SF(J)>$, $<\Gamma ^{H},\sigma
^{H},\widehat{\sigma }^{H},MF($\textbf{P}$)>\}$. Then:\newline

1. $T\vdash \widehat{\sigma }^{\prime }\leftrightarrow \urcorner \sigma
^{\prime }$,

2. $T\cup \Gamma ^{\prime }\vdash \sigma ^{\prime }$ \qquad whenever $\qquad
\Gamma \stackrel{g}{\vdash }\sigma $.\newline

\noindent \textit{\noindent Proof.}

1. Easy.

2. By proposition 1.3 for $T=MF($\textbf{P}$)$, and by the definitions for
the other cases. \hspace{12.2cm}$\square $\newline

\noindent PROPOSITION 1.14. Let $T\in \{MF($\textbf{P}$),DF(C)\}$, let $h$
be an interpretation of $BS(J)$ in the set of basic sentences of $T$, and
let $\Gamma \cup \{\sigma \}\subseteq BS(J)$. Then:
\newline
$T\cup \Gamma ^{h}\vdash \sigma ^{h}$ \qquad whenever $\qquad SF(J)\cup
\Gamma \vdash \sigma $.\newline

\noindent \textit{Proof.} The interpretations of the axioms of $SF(J)$ are
theorems of $T$, and the proof machinery of $T$ is not weaker than that of $%
SF(J)$. \hspace{2.3cm}$\square $\newline

To investigate the converses of proposition 1.14 and part 2 of proposition
1.13, we first go to:\newline

\textbf{2. Semantics of AAS. }The theories $MF(\mathbf{P})$ and $DF(C)$ are
first order, and the theory $SF(J)$ is sentential; so each has its usual
class of models with respect to which it is sound and complete.\newline

\textbf{2.1. Models of }$MF($\textbf{P}$)$. A model $\frak{B}$ of $MF($%
\textbf{P}$)$ is an ordered pair $<B,\mu >$ where $B$ is a non-empty set and 
$\mu $ maps \textbf{P} into $\wp (B)-\{\phi \}$. $B$ is called the universe,
or the base, of $\frak{B}$ and may be denoted also by 
``\TEXTsymbol{\vert}$\mathfrak{B}$\TEXTsymbol{\vert}''.
\newline
\textbf{2.2. Models of }$SF(J)$. A model $\frak{B}$ of $SF(J)$ is a mapping
form $AS(J)$ into 2 ($=\{0,1\}$), which satisfies all the axioms of $SF(J)$. 
For $\sigma \in S(J)$, $\frak{B}\vDash \sigma $ means that $\sigma $ takes
the value 1 under the usual extension of $\frak{B}$.\ \newline
\textbf{2.3. Models of }$DF(C)$. A structure $\frak{B}$ of the dyadic
language $C$ (or a $DF(C)$-structure $\frak{B}$) is a 6-tuple $%
<B,A^{*},E^{*},I^{*},O^{*},\mu >$ where $B$ is a non-empty set and $%
A^{*},E^{*},I^{*}$ and $O^{*}$ are binary relations on $B$ corresponding to
the relation symbols ``$A$'',``$E$'',``$I$'' and ``$O$'' respectively, and $\mu $ is a mapping
of $C$ into $B$. $B$ is called the universe, or the base, of $\frak{B}$ and
may be denoted also by 
``\TEXTsymbol{\vert}$\frak{B}$\TEXTsymbol{\vert}''. 
$\mathfrak{B}$ is a model of $DF(C)$ if it satisfies its axioms.

Since, by axioms 1 and 2, $E^{*}=I^{*c}$ and $O^{*}=A^{*c}$ (where ``$c$''
denotes the complement with respect to $B\times B$) we may -by abuse of
notation- say that $<B,A^{*},I^{*},\mu >$ is a model of $DF(C)$ whenever the
expansion $<B,A^{*},I^{*c},I^{*},A^{*c},\mu >$ is a model of $DF(C)$.\newline

\noindent PROPOSITION 2.1. Let $B\neq \phi $, $\mu :C\rightarrow B$ and $%
R_{1},R_{2}\subseteq B\times B$. Then $<B,R_{1},R_{2},\mu >$ is a model of $%
DF(C)$ iff:\newline

1. $R_{1}$ is reflexive and transitive (i.e. $R_{1}$ is a pre-ordering on $B$%
),

2. $R_{2}$ is symmetric,

3. $R_{1}\subseteq R_{2}$,

4. $R_{2}|R_{1}\subseteq R_{2}$, where $R_{2}|R_{1}$ is the relative product
of $R_{2}$ and $R_{1}$.\newline

\noindent \textit{Proof.}\textbf{\ }$<B,R_{1},R_{2}^{c},R_{2},R_{1}^{c},\mu
> $ satisfies axioms 1-6 iff conditions 1-3 above are satisfied.

Axiom 7 is equivalent to $\forall x\forall z[\exists y(Axy\wedge
Eyz)\rightarrow Exz]$. So axiom 7 is satisfied iff $R_{1}|R_{2}^{c}\subseteq
R_{2}^{c}$ which, in the presence of condition 2, is equivalent to condition
4.$ \hspace{12cm}\square 
$\newline

\textbf{2.4. Models of}$NF(C)$\textbf{.} The structures in which $NF(C)$
may be interpreted (henceforth $NF(C)$-structures) are exactly the $DF(C)$%
-structures. For $\Gamma \cup \{\sigma \}\subseteq BN(C)$ we write $\Gamma
\vDash _{\frak{B}}\sigma $ to mean that $\frak{B}\vDash \sigma $ whenever $%
\frak{B}\vDash \Gamma $, where $\frak{B}$ is an $NF(C)$-structure and $\frak{%
B}\vDash \sigma $, $\frak{B}\vDash \Gamma $ are defined as usual.

Two $NF(C)$-structures $\frak{B}$, $\frak{B}^{\prime }$ are said to be $%
BN(C) $-equivalent, basically equivalent, or (for short) $B$-eq if for every 
$\sigma \in BN(C)$, $\frak{B}\vDash \sigma $ iff $\frak{B}^{\prime }\vDash
\sigma $; in this case we may say also that $\frak{B}$ is $B$-eq to $\frak{B}%
^{\prime }$. This notion may be extended in an obvious way to the other
formalization of $AAS$.\newline

\noindent DEFINITION 2.2. An $NF(C)$-structure $\frak{B}$ is said to be a
direct model (or, for short, a $d$-model) if for every $\Gamma \cup \{\sigma
\}\subseteq BN(C)$, $\Gamma \vDash _{\frak{B}}\sigma $ whenever $\Gamma 
\stackrel{d}{\vdash }\sigma $.\newline

The proof of the following is straightforward.\newline

\noindent PROPOSITION 2.3. An $NF(C)$-structure $%
<B,A^{*},E^{*},I^{*},O^{*},\mu >$ is a $d$-model iff all of the rules of
inference of $NF(C)$ are valid in it (in the sense that if the antecedents
are true in it, then so also is the consequent), iff:\newline

1. $A^{*\mu }$ is reflexive on $\mu (C)$ and transitive (equivalently, $%
A^{*\mu }$ is a pre-ordering on $\mu (C)$),

2. $A^{*\mu }\subseteq \stackrel{\smallsmile }{I^{*\mu }}$,

3. $E^{*\mu }$ is symmetric,

4. $A^{*\mu }|E^{*\mu }\subseteq E^{*\mu }$,

\noindent where, for a set $X$, $X^{\mu }=X\cap (\mu (C)\times \mu (C)$); and for a
binary relation $R$, $\stackrel{\smallsmile }{R}$ is its converse.$ \hspace{11.2cm}\square $%
\newline

\noindent DEFINITION 2.4. The canonical structure $\mathfrak{B}_{\Gamma }$
corresponding to $\Gamma (\subseteq BN(C)$) is the $NF(C)$-structure $%
<B,A^{*},E^{*},I^{*},O^{*},\mu >$ satisfying:\newline

1. $B=C$,

2. $\mu =\mathfrak{l}_{C}$, where for a set $X$, $\mathfrak{l}_{X}$ is the identity
function on $X$,

3. for every $Y\in \{A,E,I,O\}$, $Y^{*}=\{<a,b>\in C\times C:Yab\in \Gamma
\} $.

An $NF(C)$-structure is said to be canonical if it is equal to $\mathfrak{B}%
_{\Gamma }$, for some $\Gamma \subseteq BN(C)$.

The basic property of $\mathfrak{B}_{\Gamma }$ is:
\newline
$\mathfrak{B}_{\Gamma }$ $\vDash \sigma $ \quad iff $\quad \sigma \in \Gamma $
\quad all $\sigma \in BN(C)$.
\newline
Every $NF(C)$-structure $\mathfrak{B}=<B,A^{*},E^{*},I^{*},O^{*},\mu >$ in which 
$\mu =\mathfrak{l}_{C}$ is an extension of a canonical structure; namely the
canonical structure corresponding to $\{\sigma \in BN(C):\mathfrak{B}\vDash
\sigma \}$.\newline

\noindent LEMMA 2.5. For $\Gamma \subseteq BN(C)$, $\mathfrak{B}_{\Gamma ^{d}}$
is a $d$-model (of $\Gamma ^{d}$ hence of $\Gamma $).\newline
\noindent \textit{Proof.}\textbf{\ }Let $\Delta \cup \{\sigma \}\subseteq
BN(C)$ and let $\Delta \stackrel{d}{\vdash }\sigma $. If $\mathfrak{B}_{\Gamma
^{d}}\vDash \Delta $ then $\Delta \subseteq \Gamma ^{d}$, hence $\sigma \in
\Gamma ^{d}$, consequently $\mathfrak{B}_{\Gamma ^{d}}\vDash \sigma $. \hspace{8cm}$\square $%
\newline

\noindent THEOREM 2.6. (Direct soundness and completeness). Direct deduction
is sound and complete with respect to the class of all direct models. That
is, for every 
$\Gamma \cup \{ \sigma \}
\subseteq BN(C)
$,
 
$$
\Gamma \stackrel{d}{\vdash}\sigma 
\text{\qquad iff \qquad} 
\Gamma 
\stackrel{d}{ \vDash} \sigma$$

\noindent where ``$\Gamma \stackrel{d}{\vDash }\sigma $'' means that $\Gamma 
\underset{\mathfrak{B}}{\vDash }\sigma $ for every $d$-model $\mathfrak{B}$.%
\newline

\noindent \textit{Proof.}\textbf{\ }Soundness is immediate by the
definition. To prove completeness assume $\Gamma \stackrel{d}{\vDash }\sigma 
$ then, in particular, $\Gamma \vDash _{\mathfrak{B}_{\Gamma ^{d}}}\sigma $. But 
$\mathfrak{B}_{\Gamma ^{d}}\vDash \Gamma $, then $\mathfrak{B}_{\Gamma ^{d}}\vDash
\sigma $, consequently $\sigma \in \Gamma ^{d}$, hence $\Gamma \stackrel{d}{%
\vdash }\sigma $.$ \hspace{8.4cm}\square $
\newline

\noindent DEFINITION 2.7. (General models). An $NF(C)$-structure $\mathfrak{B}$
is said to be a general model (or, for short, a $g$-model) if for every $%
\Gamma \cup \{\sigma \}\subseteq BN(C)$, $\Gamma \underset{\mathfrak{B}}{%
\vDash }\sigma $ whenever $\Gamma \stackrel{g}{\vdash }\sigma $.\newline

\noindent LEMMA 2.8. For $\Gamma \subseteq BN(C)$, $\mathfrak{B}_{\Gamma ^{g}}$
is a $g$-model (of $\Gamma ^{g}$ hence of $\Gamma $).\newline

\noindent \textit{Proof.}\textbf{\ }Replace ``$d$'' by ``$g$'' in the proof
of lemma 2.5. \hspace{4.3cm}$\square $
\newline

\noindent THEOREM 2.9. (General soundness and completeness). General
deduction is sound and complete with respect to the class of all general
models. That is, for every $\Gamma \cup \{\sigma \}\subseteq BN(C)$, 
\[
\Gamma \stackrel{g}{\vdash }\sigma \text{ \qquad iff \qquad }\Gamma 
\stackrel{g}{\vDash }\sigma 
\]
where ``$\Gamma \stackrel{g}{\vDash }\sigma $'' means that $\Gamma 
\underset{\mathfrak{B}}{\vDash }\sigma $ for every $g$-model $\mathfrak{B}$.%
\newline
\noindent \textit{Proof.}\textbf{\ }Replace ``$d$'' by ``$g$'' in the proof
of lemma 2.6. \hspace{4.3cm}$\square $
\newline

\noindent THEOREM 2.10. ($NF(C)$-compactness). For every $e\in \{d,g\}$, for
every $\Gamma \cup \{\sigma \}\subseteq BN(C)$, 
\[
\Gamma \stackrel{e}{\vDash }\sigma \text{ \qquad iff \qquad for some finite }%
\Gamma _{1}\subseteq \Gamma ,\text{ }\Gamma _{1}\stackrel{e}{\vDash }\sigma
. 
\]

\noindent \textit{Proof.}\textbf{\ }By $e$-soundness and $e$-completeness. \hspace{5.5cm}$%
\square $
\newline

In theorem 9.4 below the $g$-models will be fully characterized. 
Now we confine ourselves to the following:\newline

\noindent REMARKS and definitions 2.11.\newline

1. Every $g$-model is a $d$-model (obvious) but not vice versa. For, let $%
\Gamma =\{Eaa\}$ for some $a\in C$, then $\mathfrak{B}_{\Gamma ^{d}}$ is a $d$%
-model but not a $g$-model.\newline

2. Every model of $DF(C)$ is obviously a $g$-model (hence a $d$-model) but
not vice versa. For let $\Gamma =\{Oaa\}$ for some $a\in C$, then $\mathfrak{B}
_{\Gamma ^{g}}$ is a $g$-model but not a model of $DF(C)$.\newline

3. For every $g$-model $\mathfrak{B}=<B,A^{*},E^{*},I^{*},O^{*},\mu >$ in which $%
\mu $ is surjective, $A^{*}=E^{*}=I^{*}=O^{*}=B\times B$ iff ($A^{*}\cap
O^{*}\neq \phi $ or $E^{*}\cap I^{*}\neq \phi $) iff for some $\sigma \in
BN(C)$, $\mathfrak{B}\vDash \sigma ,\widehat{\sigma }$. Such models are called
full models.

If $A^{*}\cup O^{*}=B\times B=E^{*}\cup I^{*}($respectively $A^{*}\cap
O^{*}=\phi =E^{*}\cap I^{*})$, $\mathfrak{B}$ is said to be complete
(respectively consistent). Thus $\mathfrak{B}$ is not full iff it is consistent. 
$\mathfrak{B}_{\phi ^{g}}$ is an example of a $g$-model in which $\mu $ is
bijective, while it is not complete.

These notions may be generalized to all $NF(C)$-structures, $\mu $ does not
have to be surjective.\newline

4. Let $\mathfrak{B}=<B,A^{*},E^{*},I^{*},O^{*},\mu >$ be an $NF(C)$-structure
in which $\mu $ is surjective. Then $\mathfrak{B}$ is not complete iff for some $%
\sigma \in BN(C)$, ($\mathfrak{B}\nvDash \sigma $ and $\mathfrak{B}\nvDash \widehat{%
\sigma }$) iff for some $\sigma \in BN(C)$, $\mathfrak{B}\nvDash \sigma $ and
for all $\rho \in BN(C)$, $\{\widehat{\sigma }\}\vDash _{\mathfrak{B}}\rho ,%
\widehat{\rho }$.\newline

5. For $\Gamma \subseteq BN(C)$, $\Gamma $ is consistent iff it has a
consistent $d$-model.\newline

6. Direct deduction is sound with respect to any class of models with
respect to which general deduction is sound.\newline

\textbf{2.5. Order models and Venn models.} To the best of my knowledge
Shepherdson (1956) was the first to make use of a version of order models;
the ordering was pre-ordering (reflexive and transitive, but not necessarily
antisymmetric) and the context was the semantics of a version of $DF(C)$.
In the context of the semantics of $NF(C)$, or versions thereof, versions of
order models were made use of in Martin (1997) and in Glashoff (2002). The
former required a model to be some variation on a lower semi-lattice with a
smallest element, the latter relaxed these conditions; none of them
mentioned that Shepherdson (1956) made use of order models.

Following Shepherdson (1956), let $R_{1}$ be a pre-ordering on a non-empty
set $B$; and following Glashoff (2002), put: 
\[
R_{2}=\{<x,y>\in B\times B:\{x,y\}\text{ has an }R_{1}\text{-lower bound}\}. 
\]
Let $\mu $ be a function from $C$ to $B$, then $<B,R_{1},R_{2},\mu >$ is a
model of $DF(C)$, hence a $g$-model and a $d$-model. Such models are said to
be order models. If $R_{1}$ is a partial ordering (equivalently,
antisymmetric) the order model will also be called partial (or
antisymmetric). $<B,R_{1}>$ is said to be the order structure underlying
the order model $<B,R_{1},R_{2},\mu >$. Notice that if $<B,R_{1},R_{2}^{%
\prime },\mu >$ is a model of $DF(C)$, then $R_{2}\subseteq R_{2}^{\prime }$.

A concrete order model (henceforth c.o.m, and c.o.ms for the plural) is an
order model in which $B$ is a collection of non-empty sets and $R_{1}$ is $%
\subseteq $, so the c.o.ms are partial. If $R^{\prime }$ is defined on such
a $B$ by $xR^{\prime }y$ iff , $x\cap y\neq \phi $ then for every $\mu
:C\rightarrow B,<B,\subseteq ,R^{\prime },\mu >$ is a model of $DF(C)$. Such
models are said to be Venn models.

In an order model $\mathfrak{B}=<B,R_{1},R_{2},\mu >$, $R_{2}$ is determined by $%
B$ and $R_{1}$, so we may write -for short- ``$\mathfrak{B}=<B,R_{1},\mu >$''.
For similar reasons we may write ``$<B,\mu >$'' to denote the Venn model $%
<B,\subseteq ,R^{\prime },\mu >$; the c.o.m with the same $B$ and $\mu $ is
denoted by ``$<B,\subseteq ,\mu >$''.

Let $B=\{\{1,2\},\{2,3\}\}$, then for every $\mu :C\rightarrow
B,<B,\subseteq ,\mu >$ is a c.o.m but not a Venn model and $<B,\mu >$ is a
Venn model but not a c.o.m. This is not always the case, for if the universe
is the set of all non-empty subsets of a non-empty set, then the model is
both a Venn model and a c.o.m. Every Venn model is embeddable in such a
model which is $B$-eq to it. A Venn model with universe $B$ is a c.o.m iff
for every $b,b^{\prime }\in B$ there is $c\in B$ such that $c\subseteq b\cap
b^{\prime }$ whenever $b\cap b^{\prime }\neq \phi $.

Let $\mathfrak{C},\mathfrak{C}^{\prime }\in $\{the class of all Venn models, the class of all c.o.ms, the class of all partial
order models, the
class of all order models\}, then for every $\mathfrak{B}\in \mathfrak{C}$ there is $\mathfrak{B}%
^{\prime }\in \mathfrak{C}^{\prime }$ which is $B$-eq to it, hence $\vDash _{%
\mathfrak{C}}$ $=$ $\vDash _{\mathfrak{C}^{\prime }}$ (with the usual meaning). This
is a corollary of the above discussion and the following observation.

Let $\mathfrak{B}=<B,R_{1},R_{2},\mu >$ be an order model. Define the function $%
^{\prime }$ from $B$ to $\wp (B)$ by $b^{\prime }=$ $_{R_{1}}[b]$, where for
a binary relation $\rho $, $_{\rho }[y]=\{x\in $ Domain $\rho :x\rho y\}$.
Let $B^{\prime }$ be the range of $^{\prime }$ and define the function $\mu
^{\prime }$ from $C$ to $B^{\prime }$ by $\mu ^{\prime }(c)=\mu (c)^{\prime
} $. Then $\mathfrak{B}^{\prime }=<B^{\prime },\subseteq ,\mu ^{\prime }>$ is a
c.o.m which is a Venn model and $^{\prime }$ is a homomorphism from $\mathfrak{B}$ onto $\mathfrak{B}^{\prime }$. It is an isomorphism iff $R_{1}$ is
antisymmetric. In all cases $\mathfrak{B}$ and $\mathfrak{B}^{\prime }$ are $B$-eq. \newline

\textbf{2.6. Models and interpretations.}
\newline

\textit{2.6.1.} $MF($\textbf{P}$)$ and $SF(J)$\textbf{.} Let $f$
be an interpretation of $BM(\mathbf{P})$ in $BS(J)$ and let $\mathfrak{B}$ be a
model of $SF(J)$. Put: 
\begin{eqnarray*}
\mu &:&\mathbf{P}\rightarrow \wp (J)-\{\phi \}. \\
\mu (Q) &=&\{j\in J:\mathfrak{B}\vDash Ajf(Q)\}
\end{eqnarray*}
then $\mathfrak{B}^{\prime }=<J,\mu >$ is a model of $MF(\mathbf{P})$. It is
easy to see that:\newline

1. For every positive universal $\alpha \in BM(\mathbf{P})$, $\mathfrak{B}%
^{\prime }\vDash \alpha $ iff $\mathfrak{B}\vDash \alpha ^{f}$.\newline

2. For every positive particular $\alpha \in BM(\mathbf{P})$, $\mathfrak{B}%
^{\prime }\vDash \alpha $ only if $\mathfrak{B}\vDash \alpha ^{f}$. The other
direction holds iff for every $i,j\in $ Range $f$ there is $k\in J$ such
that both $\mathfrak{B}\vDash Aki$ and $\mathfrak{B}\vDash Akj$ whenever $\mathfrak{B}%
\vDash Iij$. In this case: 
\[
\mathfrak{B}^{\prime }\vDash \alpha \text{ \quad iff \quad }\mathfrak{B}\vDash
\alpha ^{f}\text{ \quad for every }\alpha \in BM(\mathbf{P}). 
\]

On the other hand, let $h$ be an interpretation of $BS(J)$ in $BM(\mathbf{P})$
and let $\mathfrak{B}$ be a model of $MF(\mathbf{P})$. Put: 
\begin{eqnarray*}
\mathfrak{B}^{\prime } &:&BS(J)\rightarrow 2 \\
\mathfrak{B}^{\prime }(\alpha ) &=&1\text{ \qquad iff \qquad }\mathfrak{B}\vDash
\alpha ^{h}\text{,}
\end{eqnarray*}
then $\mathfrak{B}^{\prime }$ is a model of $SF(J)$.\newline

\textit{2.6.2.} $SF(J)$ and $DF(C)$\textbf{.} Let $f$ be
an interpretation of $BS(J)$ in $BD(C)$ and let $\mathfrak{B}$ be a model of $%
DF(C)$. Put: 
\begin{eqnarray*}
\mathfrak{B}^{\prime } &:&BS(J)\rightarrow 2 \\
\mathfrak{B}^{\prime }(\alpha ) &=&1\text{ \qquad iff \qquad }\mathfrak{B}\vDash
\alpha ^{f}\text{,}
\end{eqnarray*}
then $\mathfrak{B}^{\prime }$ is a model of $SF(J)$.\newline

On the other hand, let $h$ be an interpretation of $BD(C)$ in $BS(J)$ and let $\mathfrak{B}$ be a model of $SF(J)$. Define: 
\begin{eqnarray*}
R_{1} &=&\{<i,j>\in J\times J:\mathfrak{B}(Aij)=1\}, \\
R_{2} &=&\{<i,j>\in J\times J:\mathfrak{B}(Iij)=1\},
\end{eqnarray*}
then $\mathfrak{B}^{\prime }=<J,R_{1},R_{2},h>$ is a model of $DF(C)$ and for
every $\alpha \in BD(C)$ 
\[
\mathfrak{B}^{\prime }\vDash \alpha \text{ \qquad iff \qquad }\mathfrak{B}\vDash
\alpha ^{h}\text{.} 
\]

\textit{2.6.3. }$DF(C)$ and $MF(\mathbf{P})$\textbf{.}
Let $f$ be an interpretation of $BD(C)$ in $BM(\mathbf{P})$ and let $\mathfrak{B}%
=<B,\mu >$ be a model of $MF(\mathbf{P})$, then $\mathfrak{B}^{\prime }=<\wp
(B)-\{\phi \},\mu \circ f>$ is a Venn model of $DF(C)$ and for every $\alpha
\in BD(C)$, $\mathfrak{B}^{\prime }\vDash \alpha $ iff $\mathfrak{B}\vDash \alpha
^{f}$.\newline

On the other hand, let $h$ be an interpretation of $BM(\mathbf{P})$ in $BD(C)$
and let $\mathfrak{B}=<B,R_{1},R_{2},\mu >$ be a model of $DF(C)$. Put: 
\begin{eqnarray*}
\mu ^{\prime } &:&\mathbf{P}\rightarrow \wp (B)-\{\phi \} \\
\mu ^{\prime }(Q) &=&_{R_{1}}[\mu h(Q)]\text{,}
\end{eqnarray*}
then $\mathfrak{B}^{\prime }=<B,\mu ^{\prime }>$ is a model of $MF(\mathbf{P})$.
It is easy to see that:\newline

1. For every positive universal $\alpha \in BM(\mathbf{P})$, $\mathfrak{B}%
^{\prime }\vDash \alpha $ iff $\mathfrak{B}\vDash \alpha ^{h}$.\newline

2. For every positive particular $\alpha \in BM(\mathbf{P})$, $\mathfrak{B}%
^{\prime }\vDash \alpha $ only if $\mathfrak{B}\vDash \alpha ^{h}$. The other
direction holds if $\mathfrak{B}$ is an order model, in this case:

\[
\mathfrak{B}^{\prime }\vDash \alpha \text{ \qquad iff \qquad }\mathfrak{B}\vDash
\alpha ^{h}\text{ \qquad for every }\alpha \in BM(\mathbf{P}). 
\]
\newline

\textbf{2.7. Leibniz models. }Let $\eta $ be the partial ordering defined on
the set $\Bbb{N}$ of natural numbers by $m\eta n$ iff $m$ is a multiple of $%
n $, and define the partial ordering $R$ on $\Bbb{N}\times \Bbb{N}$ by $%
<m_{1},n_{1}>R<m_{2},n_{2}>$ iff $m_{1}\eta m_{2}$ and $n_{1}\eta n_{2}$.
Denote the binary operations of the greatest common divisor and the least
common multiple on $\Bbb{N}$ by $\stackrel{\circ }{\wedge }$ and $\stackrel{%
\circ }{\vee }$ respectively, and put: 
\[
B=\{<m,n>\in \Bbb{N\times N}:m\stackrel{\circ }{\wedge }n=1\}. 
\]

The restriction of $R$ on $B$, to be also denoted by ``$R$'', partially
orders $B$. So, for every $\mu :C\rightarrow B$, $<B,R,\mu >$ is an order
model of $DF(C)$, hence a $g$-model and a $d$-model. Such models are called
Leibniz models, for they were first introduced -in a different setting- by
him in 1679, as may be learned from \L ukasiewicz (1998, pp. 126-9), Kneale
and Kneale (1966, pp. 337-8) and Glashoff (2002).
Leibniz practically defines $A^{*}$ to be $R$, but he sets $%
I^{*}<m_{1},n_{1}>$ 
$<m_{2},n_{2}>$ iff $m_{1}\stackrel{\circ }{\wedge }n_{2}=1=n_{1}\stackrel{%
\circ }{\wedge }m_{2}$. To show that this gives rise to an order model as
defined in 2.5 above, notice that $\{<m_{1},n_{1}>,<m_{2},n_{2}>\}$ has an $%
R $-lower bound iff there is $<m_{3},n_{3}>\in B$ such that $<m_{3},n_{3}>R$%
$<m_{1},n_{1}>$ and $<m_{3},n_{3}>R<m_{2},n_{2}>$, which is equivalent to $%
m_{3}\eta (m_{1}\stackrel{\circ }{\vee }m_{2})$ and $n_{3}\eta (n_{1}%
\stackrel{\circ }{\vee }n_{2})$. But $m_{3}\stackrel{\circ }{\wedge }n_{3}=1$%
, so the condition is equivalent to $(m_{1}\stackrel{\circ }{\vee }m_{2})%
\stackrel{\circ }{\wedge }(n_{1}\stackrel{\circ }{\vee }n_{2})=1$. The
l.h.s. $=(m_{1}\stackrel{\circ }{\wedge }n_{1})\stackrel{\circ }{\vee }(m_{1}%
\stackrel{\circ }{\wedge }n_{2})\stackrel{\circ }{\vee }(m_{2}\stackrel{%
\circ }{\wedge }n_{1})\stackrel{\circ }{\vee }(m_{2}\stackrel{\circ }{\wedge 
}n_{2})$. But $m_{1}\stackrel{\circ }{\wedge }n_{1}=1=m_{2}\stackrel{\circ }{%
\wedge }n_{2}$, so the condition is equivalent to $(m_{1}\stackrel{\circ }{%
\wedge }n_{2})\stackrel{\circ }{\vee }(m_{2}\stackrel{\circ }{\wedge }%
n_{1})=1$ which is equivalent to Leibniz condition. Via reductio ad absurdum
Glashoff (2002) gave a different proof of the same result.

Every Leibniz model is isomorphic to a Venn model. The converse is not true,
for $B$ is denumerable while there are non-denumerable Venn models.\newline

\textit{2.7.1. Assigning Leibniz models. }For $\Gamma \subseteq
BN(C)$ put: \newline

$C_{\Gamma }=\{c\in C:c$ occurs in some element of $\Gamma $ having two
distinct categorical constants$\}$.

$\Gamma $ will be called essentially finite if $C_{\Gamma }$ is finite. This
notion may be generalized to subsets of $BM(\mathbf{P})$ and $BS(J)$.\newline

\noindent LEMMA 2.12. $C_{\Gamma ^{d}}=C_{\Gamma }$.\newline

\noindent \textit{Proof}. By proposition 1.11 and induction on the length of
the deduction. \hspace{0.3cm} $\square $
\newline

\noindent THEOREM 2.13. To each consistent essentially finite $\Gamma
\subseteq BN(C)$ a Leibniz model of $\Gamma $ may be assigned (cf. Glashoff
2010, Lemma 3.4).\newline

\noindent \textit{Proof}. Let $<B,R>$ be the order structure underlying the
Leibniz models. Put $\ell =|C_{\Gamma }|$ and let $<c_{i}>_{i\in \ell }$ be
an injective enumeration of $C_{\Gamma }$, $<p_{i}>_{i\in \ell }$ be an
injective $\ell $-sequence of primes and $b\in B$. Define $\mu :C\rightarrow
B$ as follows: 
\[
\mu (c)=b\text{ \qquad if \qquad }c\in C-C_{\Gamma }\text{,} 
\]
and for $i\in \ell $, $\mu (c_{i})=<m_{i},n_{i}>$ where 
\[
m_{i}=\underset{
\begin{array}{c}
j\in \ell \\ 
Ac_{i}c_{j}\in \Gamma ^{d}
\end{array}
}{\prod }p_{j}\text{ \qquad , \qquad }n_{i}=\underset{
\begin{array}{c}
j\in \ell \\ 
Ec_{i}c_{j}\in \Gamma ^{d}
\end{array}
}{\prod }p_{j}\text{;} 
\]
$m_{i}$ and $n_{i}$ are square free finite products (the empty product is
equal to 1). By the consistency of $\Gamma $, $m_{i}\stackrel{\circ }{\wedge 
}n_{i}=1$ for all $i\in \ell $. Therefore $\mathfrak{B}=<B,R,\mu >$ is a Leibniz
model.\newline

To show that $\mathfrak{B}\vDash \Gamma $:\newline

1. Let $i,k\in \ell $, then $Ac_{i}c_{k}\in \Gamma ^{d}$ only if $(\forall
j\in \ell )[(Ac_{k}c_{j}\in \Gamma ^{d}\rightarrow Ac_{i}c_{j}\in \Gamma
^{d})\wedge (Ec_{k}c_{j}\in \Gamma ^{d}\rightarrow Ec_{i}c_{j}\in \Gamma
^{d})]$ only if $<m_{i},n_{i}>R<m_{k},n_{k}>$ only if $m_{i}\eta p_{k}$ only
if $Ac_{i}c_{k}\in \Gamma ^{d}$. So $Ac_{i}c_{k}\in \Gamma ^{d}$ iff $\mathfrak{B}\vDash Ac_{i}c_{k}$. Consequently, for every $c,c^{\prime }\in C$, $\mathfrak{B}%
\vDash Acc^{\prime }$ if $Acc^{\prime }\in \Gamma $.

Moreover if for some $c,c^{\prime }\in C$, $Occ^{\prime }\in \Gamma $, then
by the consistency of $\Gamma $ there are $i,k\in \ell $ such that $i\neq k$
and $c=c_{i}$, $c^{\prime }=c_{k}$.

Again by the consistency of $\Gamma $, $Ac_{i}c_{k}\notin \Gamma ^{d}$,
hence $\mathfrak{B}\nvDash Ac_{i}c_{k}$, consequently $\mathfrak{B}\vDash
Oc_{i}c_{k} $.\newline

2. Let $i,k\in \ell $ be such that $i\neq k$, then $Ec_{i}c_{k}\in \Gamma
^{d}$ only if $n_{i}\eta p_{k}$ only if $n_{i}\stackrel{\circ }{\wedge }%
m_{k}\neq 1$ only if $\mathfrak{B}\nvDash Ic_{i}c_{k}$ only if $\exists j\in
\ell [Ac_{i}c_{j}$, $Ec_{k}c_{j}\in \Gamma ^{d}\vee Ac_{k}c_{j}$, $%
Ec_{i}c_{j}\in \Gamma ^{d}]$ only if $Ec_{i}c_{k}\in \Gamma ^{d}$. From this
and the consistency of $\Gamma $ it follows that for every $c,c^{\prime }\in
C$, $\mathfrak{B}\vDash Ecc^{\prime }$, $\mathfrak{B}\vDash Icc^{\prime }$ whenever $%
Ecc^{\prime }\in \Gamma $, $Icc^{\prime }\in \Gamma $ respectively. \hspace{11.34cm}$\square $%
\newline

\textit{2.7.2. Leibniz soundness and completeness.} For $e \in \{d,g\}$, $e$%
-deduction is sound with respect to the set of all Leibniz models (to be
denoted, henceforth, by ``$L$'') as they are order models. Regarding
completeness, for $\Gamma \cup \{\sigma \}\subseteq BN(C)$ put: 
\[
\Gamma \vDash _{L}\sigma \text{ \quad iff \quad }\Gamma \vDash _{\mathfrak{B}%
}\sigma \text{ for every }\mathfrak{B}\in L\text{.} 
\]

\noindent THEOREM 2.14. If $\Gamma $ is essentially finite, then $\Gamma 
\stackrel{g}{\vdash }\sigma $ whenever $\Gamma \underset{L}{\vDash }\sigma 
$.\newline

\noindent \textit{Proof}. Obvious if $\Gamma $ is inconsistent. Let $\Gamma $
be consistent and $\Gamma \underset{L}{\vDash }\sigma $, then by theorem
2.13 $\Gamma \cup \{\widehat{\sigma }\}$ is inconsistent, from which the
result follows.$ \hspace{2cm} \square $
\newline

\noindent REMARKS 2.15.\newline

1. From \L ukasiewicz (1998, pp. 126-9) it follows that: 
\[
SF(J)\vdash \alpha \text{ \qquad iff \qquad }\vDash _{L^{\prime }}\alpha 
\]
where $\alpha $ is any sentence (not necessarily basic) of the language $J$,
and $L^{\prime }$ is the obvious adaptation of $L$ to $J$. Consequently, for
every $\Gamma \cup \{\alpha \}\subseteq BS(J)$: 
\[
SF(J)\cup \Gamma \vdash \alpha \text{ \qquad only if \qquad }\Gamma \vDash
_{L^{\prime }}\alpha \text{,} 
\]
the other direction holds if $\Gamma $ is essentially finite.\newline

2. In the above remark, as well as in theorem 2.14, only square free Leibniz
models (with the obvious definition) may be taken into consideration.\newline

\textit{2.7.3. Generalization.} Theorem 2.13 cannot be unconditionally
generalized to infinite $C_{\Gamma}$. For, let $<c_{i}>_{i\in \Bbb{N}}$ be an
injective enumeration of some denumerable subset of $C$. Put: 
\[
\Gamma =\{Ac_{i}c_{i+1}:i\in \Bbb{N}\}\cup \{Oc_{i+1}c_{i}:i\in \Bbb{N}\} 
\]
then $\Gamma $ is consistent but has no Leibniz model, though it has a Venn
model. The following theorem gives a sufficient condition for $\Gamma $ to
have a Leibniz model if $C_{\Gamma }$ is denumerable.\newline

\noindent THEOREM 2.16. Let $C_{\Gamma }$ be denumerable and let $%
<c_{i}>_{i\in \Bbb{N}}$ be an injective enumeration of it. Then $\Gamma $
has a Leibniz model if it is consistent and for every $i\in \Bbb{N}$, $%
\{q\in \Bbb{N}:Ac_{i}c_{q}\in \Gamma ^{d}\}$ is finite.\newline

\noindent \textit{Proof}. Along the lines of the proof of theorem 2.13 with
the following modifications. Let $<p_{i}>_{i\in \Bbb{N}}$ be an injective
enumeration of the primes, put:

\[
m_{i}=\underset{Ac_{i}c_{j}\in \Gamma ^{d}}{\prod }p_{j}\text{ \qquad ,
\qquad }n_{i}=\underset{
\begin{array}{c}
Ec_{i}c_{j}\in \Gamma ^{d} \\ 
j<\max \{q\in \Bbb{N}:Ac_{i}c_{q}\in \Gamma ^{d}\}
\end{array}
}{\prod }p_{j}\text{.} \hspace{1.8cm} \square
\]  
\newline

To see that the condition of the above theorem is essentially necessary, 
define the equivalence relation $\sim _{\Gamma }$ on $C_{\Gamma }$ by $a\sim
_{\Gamma }b$ iff $Aab,Aba\in \Gamma ^{d}$ (cf. section 1.5 above).\newline

\noindent THEOREM 2.17. If $\Gamma $ has a Leibniz model then there is a
consistent extension $\Gamma ^{\prime }$ of $\Gamma $ such that $C_{\Gamma
^{\prime }}=C_{\Gamma },$ $C_{\Gamma ^{\prime }}/\sim _{\Gamma ^{\prime }}$
is countable and for every $a\in C_{\Gamma ^{\prime }}$, with at most two
exceptions, $Q_{a}(=\{c\in C_{\Gamma ^{\prime }}:Aac\in \Gamma ^{\prime
d}\}/\sim _{\Gamma ^{\prime }})$ is finite.\newline

\noindent \textit{Proof}. Let $\mathfrak{B}(=<B,R,\mu >)$ be a Leibniz model of $%
\Gamma $. Put: 
\[
\Gamma ^{\prime }=\Gamma \cup \{\sigma \in BN(C_{\Gamma }):\mathfrak{B}\vDash
\sigma \}\text{,} 
\]
then $C_{\Gamma ^{\prime }}=C_{\Gamma }$ and $C_{\Gamma ^{\prime }}/\sim
_{\Gamma ^{\prime }}$ is countable.

If for some $a\in C_{\Gamma ^{\prime }}$, $Q_{a}$ is infinite, then $\mu
(a)\in \{<0,1>,<1,0>\}$ from which the last part of the theorem follows. \hspace{5.8cm} $%
\square $
\newline

\noindent REMARKS 2.18.\newline

1. In the underlying order structure of a Leibniz model $\mathfrak{B}=<B,R,\mu >$%
, $<1,1>$ is the greatest element and $<0,1>,<1,0>$ are the only minimal
elements. Let $a,c\in C$. If $\mu (a)=<1,1>$ then $\mathfrak{B}\vDash Aca$. Also
assuming that $\mu (a)\in \{<0,1>,<1,0>\}$, then $\mathfrak{B}\vDash Aca$
implies that $\mu (c)=\mu (a)$ hence $\mathfrak{B}\vDash Aac$, and $\mathfrak{B}%
\vDash Ica$ implies that $\mathfrak{B}\vDash Aac$.\newline

2. There would be no exceptions in the above theorem had $\Bbb{N}$ been
replaced by $\Bbb{N}^{+}$ in the definition of Leibniz models, which is
equivalent to excluding $<0,1>$ and $<1,0>$ from the universe of Leibniz
models.\newline

3. Noticing that $c\sim _{\Gamma }c^{\prime }$ forces $c,c^{\prime }$ to be
assigned the same value in any Leibniz model of $\Gamma $, with a slight
modification of its proof, theorem 2.16 may be strengthened as follows:

$\Gamma $ has a Leibniz model if there is a consistent extension $\Gamma
^{\prime }$ of $\Gamma $ such that:

1. $C_{\Gamma ^{\prime }}=C_{\Gamma }$.

2. $C_{\Gamma ^{\prime }}/\sim _{\Gamma ^{\prime }}$ is countable.

3. For every $a\in C_{\Gamma },$ $\{c/\sim _{\Gamma }:Aac\in \Gamma ^{\prime
d}\}$ is finite.\newline

4. The above strengthening is very close to be the converse of theorem 2.17.
As a matter of fact, it is its converse had $\Bbb{N}$ been replaced by $\Bbb{N%
}^{+}$ in the definition of Leibniz models.\newline

5. The completeness theorem 2.14 may be generalized in line with the above
generalizations.\newline

\textit{2.7.4. Logico-philosophical discussion of Leibniz models.}
``It is strange that his [Leibniz's] philosophic intuitions, which guided
him in his research, yielded such a sound result.'' says \L ukasiewicz
(1998, p. 126). Hopefully the above reasoning would make matters less
strange.

Following is a further discussion taking into consideration the Liebnizian
correlation between prime and composite numbers on one hand and atomic and
composite sentences, propositions, concepts or attributes on the other hand
(cf. Glashoff 2002, 2010).

If the primes $p_{1},p_{2}$ correspond, respectively, to the atomic
sentences $p_{1}^{\prime },p_{2}^{\prime }$, it is natural to let the
composite number $p_{1}p_{2}$ correspond to the composite sentence $%
p_{1}^{\prime }\wedge p_{2}^{\prime }$ . The difficulty here is that $%
p_{1}^{2}$, which is not equal to a prime, would correspond to the sentence $%
p_{1}^{\prime }\wedge p_{1}^{\prime }$, which is equivalent to an atomic
sentence; as conjunction of sentences is idempotent, while multiplication of
numbers is not. Obviously this difficulty will not arise for square free
numbers.

Notice that in the definitions of $\mu :C\rightarrow B$ given above, the
values assigned by $\mu $ to the elements of $C_{\Gamma }$ are always
ordered pairs of square free numbers. Extending this property to all
elements of $C$, after relaxing it to permit $<0,1>,<1,0>$ also to be taken
as values, gives rise to what will be called essentially square free Leibniz
models.

To investigate the relationship between the Leibniz models and the
essentially square free Leibniz models, let $<q_{ij}>_{i,j\in \Bbb{N}}$ be
an injective double sequence of primes. Map the k$^{th}$ power of the i$%
^{th} $ prime $p_{i}$ on $\underset{j\in k}{\prod }q_{ij}$. This mapping
may be extended in the obvious way to an injection $\nu $ from $\Bbb{N}$ to $%
\Bbb{N} $ such that $\nu (0)=0$ and for $n\geq 1$, $\nu (n)$ is square free
(being the empty product of primes, $\nu (1)=1$). The mapping $\nu $ may be
further extended, in the obvious way, to $\Bbb{N\times N}$, the extension
also will be denoted by ``$\nu $''. It may be easily seen that $\nu
(B)\subseteq B$ and that $<m_{1},n_{1}>R<m_{2},n_{2}>$ iff $\nu
(<m_{1},n_{1}>)R\nu (<m_{2},n_{2}>)$ for every $<m_{1},n_{1}>,<m_{2},n_{2}>%
\in B$. So for every Leibniz model $\mathfrak{B}(=<B,R,\mu >)$, $\nu $ is a
monomorphism from $\mathfrak{B}$ into $\mathfrak{B}^{\nu }(=<B,R,\nu \mu >)$ which
is essentially square free and is basically equivalent to $\mathfrak{B}$.
Moreover if in $\mathfrak{B}^{\nu }$, $B$ is replaced by $\nu (B)$ and $R$ by $%
R\cap (\nu (B)\times \nu (B))$, then $\nu $ will be an isomorphism. Such
models will be called proper Leibniz models. Since every Leibniz model is
isomorphic to a proper Leibniz model, attention may be confined to the latter.

Let $<q_{ij}^{\prime }>_{i,j\in \Bbb{N}}$ be an injective double sequence of
atomic sentences in some sentential language. For $<m,n>\in \nu (B)$ put: 
\[
\lambda (<m,n>)=\underset{m\eta q_{ij}}{\bigwedge }q_{ij}^{\prime }\text{ }%
\wedge \underset{n\eta q_{ij}}{\bigwedge }\urcorner q_{ij}^{\prime }. 
\]

As $1\eta p$ for no prime $p$, $\underset{1\eta q_{ij}}{\bigwedge }%
q_{ij}^{\prime }=\underset{1\eta q_{ij}}{\bigwedge }\urcorner
q_{ij}^{\prime }=$ the empty conjunction, which is always true. So $\lambda
(<1,1>)$ is always true.

On the other hand, $0 \eta p$ for every prime $p$, so 
\[
\lambda (<0,1>)=\underset{i,j\in \Bbb{N}}{\bigwedge }q_{ij}^{\prime }\text{
\qquad and \qquad }\lambda (<1,0>)=\underset{i,j\in \Bbb{N}}{\bigwedge }%
\urcorner q_{ij}^{\prime }. 
\]

These are the only infinitary sentences to be considered.\newline

It may be easily seen that for every proper Leibniz model $\mathfrak{B}$, $\mathfrak{B}\vDash Ac_{1}c_{2}$ iff $\lambda \mu c_{1}\rightarrow \lambda \mu c_{2}$
is a tautology, and $\mathfrak{B}\vDash Ec_{1}c_{2}$ iff $\lambda \mu
c_{1}\wedge \lambda \mu c_{2}$ is a contradiction.

As a matter of fact $<0,1>,<1,0>$ and $<1,1>$ are not indispensable as
elements of the universe of proper Leibniz models. To keep them or not is a
philosophical choice. Rejecting them is probably more compatible with the
Aristotelian legacy.

Following Boole (1948, p. 49), to each $<m,n>\in \nu (B)$ the set $\theta
(<m,n>)$ of all truth assignments which satisfy $\lambda (<m,n>)$ may be
appropriated. For every proper Leibniz model (hence for every Leibniz model) 
$\mathfrak{B}$, $\theta $ induces an isomorphism of $\mathfrak{B}$ onto a Venn model
which is a concrete order model.\newline

\textbf{\noindent \noindent 3. Decidability.}\newline

\noindent REMARKS 3.1. Let $a,b,c\in C$ and $\Gamma \subseteq BN(C)$.\newline

1. If $\{Ecc,Occ\}\cap \Gamma \neq \phi $, $\Gamma $ may be easily seen to
be contradictory. In such a case $\Gamma $ is said to be plainly
contradictory.\newline

2. $\Gamma \stackrel{d}{\vdash }Oab$ iff $Oab\in \Gamma $.\newline

3. If $\Gamma \stackrel{d}{\vdash }Ecc$ then $Ecc\in \Gamma $ or $c\in
C_{\Gamma }$.\newline

4. If $\Gamma $ is not plainly contradictory, then $\Gamma $ is
contradictory iff there are $\sigma ,\widehat{\sigma }\in (\Gamma ^{d}\cap
BN(C_{\Gamma }))$.\newline

5. $\Gamma ^{d}\cap BN(C_{\Gamma })=(\Gamma \cap BN(C_{\Gamma }))^{d}\cap
BN(C_{\Gamma })$.\newline

6. In a different context, Glashoff (2005) presents an algorithm which may
be regarded as a prelude to the one given below. Roughly speaking, it
amounts -in our terminology- to: For a finite $\Gamma (\subseteq BN(C))$, $%
\Gamma ^{d}\cap BN(C_{\Gamma })$ may be obtained from $\Gamma $ in finitely
many steps.\newline

\noindent THEOREM 3.2. There is a polynomial (of degree 8) time algorithm to
decide for any essentially finite $\Gamma (\subseteq BN(C))$ which is not
plainly contradictory whether it is contradictory, and to assign a Leibniz
model to it if it is not.\newline

\noindent \textit{Proof}. Let $\Gamma $ satisfy the conditions of the
theorem, then $BN(C_{\Gamma })$ is finite. Put $\Gamma ^{\prime }=\Gamma
\cap BN(C_{\Gamma })$ and $\Delta =\Gamma ^{\prime d}\cap BN(C_{\Gamma })$.

The input of the algorithm is $\Gamma ^{\prime }$ structured as a list $%
<\gamma _{i}>_{i\in n}$ where $n=|\Gamma ^{\prime }|$, and for every $i\in n$%
, $\gamma _{i}=<\gamma _{ij}>_{j\in 3}$ where $\gamma _{io}\in \{A,E,I,O\}$
and $\gamma _{i1},\gamma _{i2}\in C_{\Gamma }$. $C_{\Gamma }$ may be
obtained from $\Gamma ^{\prime }$ or supplied as a secondary input. $%
|C_{\Gamma }|\leq 2n$ and $|BN(C_{\Gamma })|\leq 16n^{2}$.

The next step is to extract for each $Y\in \{A,E,I,O\}$, $\Gamma
_{Y}^{\prime }$ (the set of all elements of $\Gamma ^{\prime }$ starting
with $Y$) which may be done through a simple scanning procedure in a linear
time. Then construct $\Delta _{Y}$ (with the obvious meaning) for each $Y\in
\{A,E,I,O\}$.

Notice that $\Delta _{O}=\Gamma _{O}^{\prime }$ and $\Delta _{A}$ is needed
to construct each of $\Delta _{E}$ and $\Delta _{I}$. To construct $\Delta
_{A}$ start with the list $\Gamma _{A}^{\prime }$. At most $16n^{4}$
comparisons are needed to determine all the possible applicabilities of
Barbara. And for each possible applicability at most $4n^{2}$ comparisons are
needed to check whether the consequent is already there. If not, append it.

It is needed to repeat this process at most $4n^{2}$ times to cover all the
required applications of Barbara. In addition, for each $c\in C_{\Gamma }$
at most $4n^{2}$ comparisons are needed to check whether $Acc$ is listed; if
not, append it. It is easy to see that this completes the construction of $%
\Delta _{A}$.

By simple variations on the above procedure $\Delta _{E}$ may be
constructed. Constructing $\Delta _{I}$ is much simpler.

$\Gamma $ is contradictory iff $\sigma ,\widehat{\sigma }\in \Delta $ for
some $\sigma $, which needs at most $32n^{4}$ comparisons to check.

If $\Gamma $ is consistent assign to it a Leibniz model along the lines of
the proof of theorem 2.13 (in the appendix a polynomial (in $n$ of degree 6)
time algorithm will be presented to generate the first n primes).

The total running time is bounded above by a polynomial (in $n$) of degree 8.%
\hspace{13cm} $\square $\noindent 
\newline

\noindent

\textbf{4. Basic equivalence of the four formalizations. }Let $\Gamma \cup
\{\sigma \}\subseteq BN(C)$, let $h$ and $H$ be bijective interpretations of 
$BN(C)$ in $BS(J)$ and $BM(\mathbf{P})$ respectively, and let $<\Gamma
^{\prime },\sigma ^{\prime },T>\in \{<\Gamma ,\sigma ,DF(C)>$, $<\Gamma
^{h},\sigma ^{h},$ $SF(J)>$, $<\Gamma ^{H},\sigma ^{H},MF(\mathbf{P})>\}$.%
\newline

\noindent THEOREM 4.1.\newline

1. $\Gamma $ is consistent \quad iff $\quad \Gamma ^{\prime }\cup T$ is.%
\newline

2. $\Gamma \vdash ^{g}\sigma $ \qquad iff $\qquad \Gamma ^{\prime }\cup
T\vdash \sigma ^{\prime }$.\newline

\noindent \textit{Proof}. From proposition 1.13, if of $(1)$ and only if of $%
(2)$ follow.

The other two directions for $<\Gamma ^{\prime },\sigma ^{\prime
},T>=<\Gamma ^{h},\sigma ^{h},SF(J)>$ follow from the corresponding
directions for $<\Gamma ^{\prime },\sigma ^{\prime },T>=<\Gamma ,\sigma
,DF(C)>$, this is a consequence of proposition 1.14. So it remains to prove
these two other directions for $<\Gamma ^{\prime },\sigma ^{\prime },T>\in
\{<\Gamma ,\sigma ,DF(C)>$, $<\Gamma ^{H},\sigma ^{H},MF(\mathbf{P})>\}$.

Only if of (1): Assume $\Gamma $ is consistent. Let $\Delta $ be a finite
subset of $\Gamma $, then by theorem 2.13 it has a Leibniz model, $\mathfrak{B}$
say. $\mathfrak{B}$ is a model of $\Delta \cup DF(C)$; from this the consistency
of $\Gamma \cup DF(C)$ follows. Since $\Delta $ may be assumed to be the
inverse image of some finite $\Delta ^{\prime }\subseteq \Gamma ^{H}$, then
by subsection 2.6.3, $\mathfrak{B}$ induces a model of $\Delta ^{\prime }\cup MF(%
\mathbf{P})$. From this the consistency of $\Gamma ^{H}\cup MF(\mathbf{P})$
follows.

If of (2): Let $\Gamma ^{\prime }\cup T\vdash \sigma ^{\prime }$, then $%
\Gamma ^{\prime }\cup \{\widehat{\sigma }^{\prime }\}\cup T$ is
inconsistent. By part (1), $\Gamma \cup \{\widehat{\sigma }\}$ is
inconsistent, hence $\Gamma \vdash ^{g}\sigma $. \hspace{6.5cm} $\square $
\newline

\noindent REMARKS 4.2.\newline

1. As far as the basic sentences are concerned, the four formalizations are
equivalent in the sense expressed by part 2 of the above theorem; so it may
be said, for brevity, that they are basically equivalent.\newline

2. In the above theorem, the only if direction of (1) and the if direction
of (2) may be directly proved for $<\Gamma ^{\prime },\sigma ^{\prime
},T>=<\Gamma ^{h},\sigma ^{h},SF(J)>$.\newline

\textbf{\noindent 5. Venn soundness and completeness. }Let $\Gamma \cup
\{\sigma \}\subseteq BN(C)$.\newline

\noindent DEFINITION 5.1.\textbf{\ }$\Gamma \underset{V}{\vDash }\sigma $
iff $\Gamma \underset{\mathfrak{B}}{\vDash }\sigma $ for every Venn model $%
\mathfrak{B}$.\newline

\noindent THEOREM 5.2.\textbf{\ }(Venn soundness and completeness). General
deduction is sound and complete with respect to the class of Venn models.
That is $\Gamma \stackrel{g}{\vdash }\sigma $ iff $\Gamma \underset{V}{%
\vDash }\sigma $.

\noindent \textit{Proof}. Every Venn model is a $DF(C)$ model (subsection
2.5), then a $g$-model (2 of remarks 2.11). This guarantees soundness.

To prove completeness, let $\Gamma \stackrel{g}{\nvdash }\sigma $, then $%
\Gamma \cup \{\widehat{\sigma }\}$ is consistent then, by theorem 4.1, $%
(\Gamma \cup \{\widehat{\sigma }\})^{H}\cup MF(\mathbf{P})$ is consistent
where $H$ is a bijective interpretation of $BN(C)$ in $BM(\mathbf{P})$, for
some appropriate $\mathbf{P}$. By well known results in first order logic, $%
(\Gamma \cup \{\widehat{\sigma }\})^{H}\cup MF(\mathbf{P})$ has a model. By
subsection 2.6.3 and 2 of remarks 2.11, $\Gamma \cup \{\widehat{\sigma }\}$
has a Venn model, hence $\Gamma \underset{V}{\nvDash }\sigma $. \hspace{2cm} $\square $%
\newline

Alternatively theorem 6.3 below may be made use of to directly show that $%
\Gamma \cup \{\widehat{\sigma }\}$ has a Venn model.\newline

\noindent REMARKS 5.3.\textbf{\ }\newline

1. In view of subsection 2.5, the above theorem entails that general
deduction is sound and complete with respect to each of the classes of
order, partial order, and concrete order models.\newline

2. Direct ways to Venn models on one hand, and to order and partial order
models on the other hand, will be presented in sections 6 and 9 respectively.%
\newline

3. For the Venn soundness and completeness of \L ukasiewicz's system,
Shepherdson (1956) may be consulted.\newline

\textbf{\noindent 6. Direct way to Venn models. }Let $\Gamma \subseteq BN(C)$%
, put $D=\{<a,b>\in C\times C:\{Iab,Iba\}\cap \Gamma ^{d}\neq \phi \}$, $%
B=\wp (D)-\{\phi \}$. Define the function $\mu $ from $C$ to $B$ by: 
\[
\mu (c)=\{<a,b>\in D:\{Aac,Abc\}\cap \Gamma ^{d}\neq \phi \}. 
\]
Then $<B,\mu >$ is a Venn model (which is a concrete order model), denote it
by ``\noindent $\mathfrak{B}^{\Gamma }$''.\newline

\noindent LEMMA 6.1.\textbf{\ }For every $c,c^{\prime }\in C$, the following
are equivalent:\newline

1. $Acc^{\prime }\in \Gamma ^{d}$,

2. $\mu (c)\subseteq \mu (c^{\prime })$ (which is equivalent to $\mathfrak{B}%
^{\Gamma }\vDash Acc^{\prime }$),

3. $<c,c>\in \mu (c^{\prime })$.\newline

\noindent \textit{Proof. }Straightforward.$ \hspace{8.8cm} \square $
\newline

\noindent LEMMA 6.2.\textbf{\ }Let $c,c^{\prime }\in C$, consider:\newline

1. $<c,c^{\prime }>\in D$,

2. $<c,c^{\prime }>\in \mu (c)\cap \mu (c^{\prime })$,

3. $\mu (c)\cap \mu (c^{\prime })\neq \phi $ (which is equivalent to $\mathfrak{B}^{\Gamma }\vDash Icc^{\prime },Ic^{\prime }c$),

4. $\{Ecc^{\prime },Ec^{\prime }c\}\cap \Gamma ^{d}=\phi $,

\noindent then 1 is equivalent to 2 which implies 3 which, for consistent $%
\Gamma $, implies 4.\newline

\noindent \textit{Proof. }The first two parts are easy to see. For the last
part assume $\mu (c)\cap \mu (c^{\prime })\neq \phi \neq \{Ecc^{\prime
},Ec^{\prime }c\}\cap \Gamma ^{d}$, then $Ecc^{\prime }\in \Gamma ^{d}$ and
there are $a,b\in C$ such that $(Iab\in \Gamma ^{d}$ or $Iba\in \Gamma ^{d})$%
, $(Aac\in \Gamma ^{d}$ or $Abc\in \Gamma ^{d})$ and $(Aac^{\prime }\in
\Gamma ^{d}$ or $Abc^{\prime }\in \Gamma ^{d})$. So there are eight cases to
consider. We deal only with the case $Iab,Aac,Abc^{\prime }\in \Gamma ^{d}$;
the other cases are similar or easier. In this case $Aac,Ecc^{\prime }\in
\Gamma ^{d}$, then $Eac^{\prime }\in \Gamma ^{d}$, then $Ec^{\prime }a\in
\Gamma ^{d}$, but $Abc^{\prime }\in \Gamma ^{d}$, then $Eba\in \Gamma ^{d}$,
then $Eab\in \Gamma ^{d}$ which contradicts that $Iab\in \Gamma ^{d}$; from
this and the consistency of $\Gamma $ the result follows. \hspace{7cm} $\square $
\newline

\noindent THEOREM 6.3.\textbf{\ }(Existence of Venn models)\textbf{. }Let $%
\Gamma \subseteq BN(C)$ be consistent, then $\mathfrak{B}^{\Gamma }$ is a Venn
model (which is a concrete order model) of $\Gamma $.\newline

\noindent \textit{Proof. }Let $c,c^{\prime }\in C$. By lemma 6.1, $%
Acc^{\prime }\in \Gamma ^{d}$ iff $\mathfrak{B}^{\Gamma }\vDash Acc^{\prime }$.
From this it follows that for $Y\in \{A,O\}$, $\mathfrak{B}^{\Gamma }\vDash
Ycc^{\prime }$ if $Ycc^{\prime }\in \Gamma ^{d}$.

Moreover, if $Icc^{\prime }\in \Gamma ^{d}$ then $<c,c^{\prime }>\in D$
then, by lemma 6.2, $\mathfrak{B}^{\Gamma }\vDash Icc^{\prime }$. Finally, if $%
Ecc^{\prime }\in \Gamma ^{d}$ then, by lemma 6.2, $\mu (c)\cap \mu
(c^{\prime })=\phi $ then $\mathfrak{B}^{\Gamma }\vDash Ecc^{\prime }$.$\square $%
\newline

Lemma 6.1 syntactically characterizes $\{Acc^{\prime }\in BN(C):\mathfrak{B}%
^{\Gamma }\vDash Acc^{\prime }\}$, hence it syntactically characterizes $%
\{Occ^{\prime }\in BN(C):\mathfrak{B}^{\Gamma }\vDash Occ^{\prime }\}$. The
following syntactical characterization of $\{Icc^{\prime }\in BN(C):\mathfrak{B}%
^{\Gamma }\vDash Icc^{\prime }\}$, hence of $\{Ecc^{\prime }\in BN(C):\mathfrak{B}^{\Gamma }\vDash Ecc^{\prime }\}$, 
\[
\mathfrak{B}^{\Gamma }\vDash Icc^{\prime }\text{ \quad iff \quad }\Gamma 
\stackrel{d^{\prime }}{\vdash }Icc^{\prime } 
\]
is an immediate consequence of lemma 8.3 below; the definition of ``$%
\stackrel{d^{\prime }}{\vdash }$'' may be found at the beginning of section
8 below.

Slightly modifying the above construction, light may be shed on the role
played by the Venn models among the models of $DF(C)$.\newline

\noindent THEOREM 6.4.\textbf{\ }For every $DF(C)$ model $\mathfrak{B}%
=<B,R_{1},R_{2},\mu >$ there is a Venn model $\mathfrak{B}^{\prime }=<B^{\prime
},\mu ^{\prime }>$ and surjection $h:B\rightarrow B^{\prime }$ such that:%
\newline

1. $\mu ^{\prime }=h\mu $ and for every $b_{1},b_{2}\in B$:

$b_{1}R_{1}b_{2}$ iff $h(b_{1})\subseteq h(b_{2})$ \quad , $\quad
b_{1}R_{2}b_{2}$ iff $h(b_{1})\cap h(b_{2})\neq \phi $,\newline

2. $\mathfrak{B}$ and $\mathfrak{B}^{\prime }$ are basically equivalent,\newline

3. $h$ is an isomorphism iff $R_{1}$ is antisymmetric.\newline

\noindent \textit{Proof. }Put: 
\begin{eqnarray*}
h &:&B\rightarrow \wp (\wp (B)) \\
h(b) &=&\{\{b_{1},b_{2}\}\in \wp (B):(b_{1}R_{2}b_{2}\text{ or }%
b_{2}R_{2}b_{1})\text{ } \\
&&\text{and }(b_{1}R_{1}b\text{ or }b_{2}R_{1}b)\}\text{,}
\end{eqnarray*}
\[
B^{\prime }=h(B)\text{ \qquad , \qquad }\mu ^{\prime }=h\mu \text{.} 
\]

The rest of the proof is easy. \hspace{7.2cm} $\square $
\newline

\textbf{7. Variations on }$NF(C)$\textbf{.} As was promised in section 1, we
follow in subsection 7.1 the long standing tradition of not permitting the
subject and the predicate of a categorical sentence to be the same. The
resulting formalization, $WF(C)$, and its relationship to $NF(C)$ are
discussed.

In subsection 7.2 the standpoint that $Acc^{\prime }$ requires that all $c$
are $c^{\prime }$ but not vice versa, will be considered.\newline

\textbf{7.1. Weak natural deduction formalization of AAS.} The alphabet of
the logical system $WF(C)$, the weak natural deduction formalization of AAS,
is the same as the alphabet of $NF(C)$. The set $W(C)$ of sentences of $%
WF(C) $ is defined as follows: 
\[
W(C)=S(C)-\{Ycc:Y\in \{A,E,I,O\}\text{ and }c\in C\}\text{.} 
\]

In accordance with subsection 1.6, the set $BW(C)$ of basic sentences of $%
WF(C)$ is $W(C)$ itself.

The rules of inference of $WF(C)$ are those of $NF(C)$ after dropping the
first one ($\frac {}{Aaa}$). The weak direct and general deduction relations
are respectively denoted by ``$\stackrel{wd}{\vdash }$'' and ``$\stackrel{wg%
}{\vdash }$'' and are defined along the lines of definitions 1.8 and 1.9
respectively. The definition of the other notions introduced in the theory
of $NF(C)$ may be modified in the obvious way to render the corresponding
definitions for the theory of $WF(C)$.

The theory of $WF(C)$ may be obtained from that of $NF(C)$ by making the
obvious modifications. The key observations are the following, where $\Gamma
\cup \{\sigma \}\subseteq BW(C)$.\newline

\noindent PROPOSITION 7.1. 
\[
\Gamma \stackrel{wd}{\vdash }\sigma \text{ \qquad iff \qquad }\Gamma \stackrel{d}{%
\vdash }\sigma . 
\]
\newline
\noindent \textit{Proof.} The only if direction is obvious. To prove the
other direction let $<\rho _{i}>_{i\in n}$ be a $d$-deduction of $\sigma $
from $\Gamma $. We show by induction that for every $i\in n$, $\Gamma 
\stackrel{wd}{\vdash }\rho _{i}$ if $\rho _{i}\in BW(C)$.

Distinguish between four cases:\newline

1. $\rho _{i}=Occ^{\prime }$, then $\rho _{i}\in \Gamma $, then $\Gamma 
\stackrel{wd}{\vdash }\rho _{i}$.\newline

2. $\rho _{i}=Acc^{\prime }$, then the result follows by the induction
hypothesis.\newline

3. $\rho _{i}=Icc^{\prime }$, then $Icc^{\prime }\in \Gamma $ or $\rho
_{j}=Ac^{\prime }c$ for some $j<i$. From this and part 2 the result follows
by the induction hypothesis.\newline

4. $\sigma =Ecc^{\prime }$, then the result follows by the induction
hypothesis noting that if $Ecc^{\prime }$ is obtained via applying $\frac{%
Acc^{\prime },Ec^{\prime }c^{\prime }}{Ecc^{\prime }}$ then the first
occurrence of $Ec^{\prime }c^{\prime }$ in the deduction must be obtained
via $\frac{Ac^{\prime }c^{\prime \prime },Ec^{\prime \prime }c^{\prime }}{%
Ec^{\prime }c^{\prime }}$ for some $c^{\prime \prime }\neq c^{\prime }$. By
part 2 and the induction hypothesis $\Gamma \stackrel{wd}{\vdash }%
Acc^{\prime },Ac^{\prime }c^{\prime \prime },Ec^{\prime \prime }c^{\prime }$%
, hence $\Gamma \stackrel{wd}{\vdash }Ecc^{\prime }$.$\square $\newline

\noindent PROPOSITION 7.2.

$\Gamma $ is $wd$-consistent iff $\quad \Gamma $ is $d$-consistent (hence $%
\Gamma \stackrel{wg}{\vdash }\sigma $ iff \quad $\Gamma \stackrel{g}{\vdash }%
\sigma $).\newline

\noindent \textit{Proof.} The if direction easily follows from proposition
7.1. To prove the other direction assume that $\Gamma $ is $d$-inconsistent,
then $\Gamma \stackrel{d}{\vdash }Ycc^{\prime }$, $\widehat{Y}cc^{\prime }$
for some $Y\in \{A,E,I,O\}$ and some $c,c^{\prime }\in C$. If $c\neq
c^{\prime }$ the result follows by the previous proposition. Else,
distinguish between two cases:\newline

1. $Occ\in \{Ycc,\widehat{Y}cc\}$, then $Occ\in \Gamma $ which is not
permitted.\newline

2. $Ecc\in \{Ycc,\widehat{Y}cc\}$. In this case there is a $d$-deduction of $%
Ecc$ from $\Gamma $. The rule made use of to justify the first occurrence of 
$Ecc$ in this deduction must be $\frac{Acc^{\prime \prime },Ec^{\prime
\prime }c}{Ecc}$ for some $c^{\prime \prime }\neq c$. By the previous
proposition $\Gamma \stackrel{wd}{\vdash }Ic^{\prime \prime }c,Ec^{\prime
\prime }c$, hence $\Gamma $ is $wd$-inconsistent. \hspace{5.2cm} $\square $
\newline

\noindent COROLLARY 7.3.\textbf{\ }$\Gamma $ is $wd$-consistent \quad iff $%
\quad \Gamma $ is $d$-consistent \quad iff $\quad \Gamma $ is $g$-consistent
\quad iff $\quad \Gamma $ is $wg$-consistent.\newline

\noindent \textit{Proof. }$\Gamma $ is $wd$-consistent only if $\Gamma $ is $%
d$-consistent only if $\Gamma $ is $g$-consistent only if $\Gamma $ is $wg$%
-consistent only if $\Gamma $ is $wd$-consistent. \hspace{4cm} $\square $
\newline

\noindent REMARK 7.4.\textbf{\ }From propositions 7.1 and 7.2 it follows
that the results concerning Leibniz soundness and completeness (subsection
2.7.2) and Venn and order soundness and completeness (section 5) apply to $%
WF(C)$ after replacing $d,g$ and $BN(C)$ by $wd,wg$ and $BW(C)$ respectively.%
\newline

\textbf{7.2. Proper natural deduction formalization of AAS.} AAS may be
interpreted to require $Acc^{\prime }$ to hold iff all $c$ are $c^{\prime }$ 
\textit{but not} vice versa. That is, extensionally, the denotation of ``$c$%
'' is required to be a proper subclass of the denotation of ``$c^{\prime }$%
''.

To satisfy this requirement introduce the logical system $PF(C)$, the proper
natural deduction formalization of AAS, based on the same language as the
system $NF(C)$. So the set $P(C)$ of sentences of $PF(C)$ is the same as $%
S(C)$. In accordance with subsection 1.6 the set $BP(C)$ of basic sentences
of $PF(C)$ is $P(C)$ itself. For ``$Occ^{\prime }$'' to remain to be the
contradictory of ``$Acc^{\prime }$'', it must be interpreted as \textit{some 
}$c$ \textit{are not }$c^{\prime }$\textit{\ or (all }$c$\textit{\ are }$%
c^{\prime }$\textit{\ and vice versa)}. The rules of inference of $PF(C)$
are to be obtained from those of $NF(C)$ by dropping the first one and
augmenting the remaining ones by $\frac {}{Icc}$($I$-$Id$) $\quad $and $%
\quad \frac {}{Occ}$($O$-$Id$).

The proper direct and general deduction relations are respectively denoted
by ``$\stackrel{pd}{\vdash }$'' and ``$\stackrel{pg}{\vdash }$'' and are
defined along the lines of definitions 1.8 and 1.9 respectively. The
definitions of the other notions introduced in the theory of $NF(C)$ may be
modified in the obvious way to render the corresponding definitions for the
theory of $PF(C)$.\newline

\noindent PROPOSITION 7.5. For $\Gamma \cup \{\sigma \}\subseteq BW(C)$:%
\newline

1. $\Gamma \stackrel{pd}{\vdash }\sigma $ $\quad $iff $\quad \Gamma 
\stackrel{wd}{\vdash }\sigma $ $\quad $(iff $\Gamma \stackrel{d}{\vdash }%
\sigma $).\newline

2. If $\Gamma $ is $pd$-consistent then it is $wd$-consistent (equivalently $%
d$-consistent), but not always vice versa.\newline

3. If $\Gamma \stackrel{wg}{\vdash }\sigma $ (equivalently $\Gamma \stackrel{%
g}{\vdash }\sigma $) then $\Gamma \stackrel{pg}{\vdash }\sigma $, but not
always vice versa.\newline

\noindent \textit{Proof.}\newline

1. The proof of part 1 is similar to that of proposition 7.1.\newline

2. That $\Gamma $ is $wd$-consistent if it is $pd$-consistent easily follows
from part 1. To see that the other direction does not always hold consider $%
\{Acc^{\prime },Ac^{\prime }c\}$ for some $c,c^{\prime }\in C$ such that $%
c\neq c^{\prime }$; this proves part 2.\newline

3. Part 3 is a direct consequence of part 2.$ \hspace{5.1cm} \square $ 
\newline

\noindent PROPOSITION 7.6.\textbf{\ }$\Delta (\subseteq BP(C))$ is $pd$%
-consistent iff it is $pg$-consistent.\newline

\noindent \textit{Proof.} Along the lines of the proof of part 4 of
proposition 1.11. \hspace{2cm} $\square $
\newline

Order models, Leibniz models, and Venn models are \textit{not} $e(\in
\{pd,pg\})$-models, so it does not make sense to ask whether $e$ is sound or
complete with respect to any of these classes. However, with some
modifications, to be shown below, everything goes as expected.

Let $\mathfrak{B}=<B,A^{*},E^{*},I^{*},O^{*},\mu >$ be a $d$-model of $\Gamma
(\subseteq BW(C))$ such that $\mu $ is injective on $C_{\Gamma }$ and $%
A^{*\mu }$ is antisymmetric, and let $A^{p},O^{p}$ be subsets of $B\times B$
such that $A^{*\mu }-\mathfrak{l}_{\mu (C)}\subseteq A^{p\mu }\subseteq A^{*\mu
} $ and $O^{*\mu }\cup \mathfrak{l}_{\mu (C)}\subseteq O^{p\mu }$. Put: 
\[
\mathfrak{B}^{p}=<B,A^{p},E^{*},I^{*},O^{p},\mu >. 
\]

\noindent PROPOSITION 7.7.\textbf{\ }$\mathfrak{B}^{p}$ is a $pd$-model of $%
\Gamma $.\newline

\noindent \textit{Proof.} Assume $\mathfrak{B}\vDash \Gamma $. To show that $%
\mathfrak{B}^{p}\vDash \Gamma $, let $\gamma \in \Gamma $ then $\gamma
=Ycc^{\prime }$ for some $Y\in \{A,E,I,O\}$ and some $c,c^{\prime }\in C$
such that $c\neq c^{\prime }$, hence $\mu (c)\neq \mu (c^{\prime })$. If $%
Y=A $ then $<\mu (c),\mu (c^{\prime })>\in A^{*\mu }-\mathfrak{l}_{\mu
(C)}\subseteq A^{p\mu }$. The other cases are obvious.

To show that $\mathfrak{B}^{p}$ is a $pd$-model, assume $\mathfrak{B}^{p}\vDash
Acc^{\prime },Ac^{\prime }c^{\prime \prime }$. If not $\mathfrak{B}^{p}\vDash
Acc^{\prime \prime }$ then $\mu (c)=\mu (c^{\prime \prime })$, then $<\mu
(c),\mu (c^{\prime })>,<\mu (c^{\prime }),\mu (c)>\in A^{*\mu }$, then $\mu
(c^{\prime })=\mu (c)=\mu (c^{\prime \prime })$ which is absurd. So Barbara
is valid. The other rules are easier to deal with. \hspace{9.1cm}$\square $
\newline

Accordingly, it is legitimate to adopt in the sequel the following
modifications: 
\[
A^{p}=A^{*}-\mathfrak{l}_{\mu (C)}\text{ \quad , \quad }O^{p}=O^{*}\cup \mathfrak{l}%
_{\mu (C)}. 
\]
\newline

\noindent THEOREM 7.8.\textbf{\ }Let $\Delta (\subseteq BP(C))$ be $pd$%
-consistent, then it has a modified Venn model which is a modified c.o.m
and which is also a $pg$-model. If, in addition, $\Delta $ is essentially
finite then it has also a modified Leibniz model which is a $pg$-model.%
\newline

\noindent \textit{Proof.} Put $\Gamma =\Delta \cap BW(C)$, then $\Gamma $ is 
$pd$-consistent, hence it is $d$-consistent, hence $\mathfrak{B}^{\Gamma }$ is a 
$d$-model of $\Gamma $ in which $A^{*}$ is antisymmetric.

To show that $\mu $ is injective let $c,c^{\prime }\in C$ be such that $%
c\neq c^{\prime }$ and $\mu (c)=\mu (c^{\prime })$, then $\Gamma \stackrel{d%
}{\vdash }Acc^{\prime },Ac^{\prime }c$, then $\Gamma \stackrel{wd}{\vdash }%
Acc^{\prime },Ac^{\prime }c$, then $\Gamma \stackrel{pd}{\vdash }Acc^{\prime
},Ac^{\prime }c$, then $\Gamma \stackrel{pd}{\vdash }Acc$ which contradicts
that $\Gamma $ is $pd$-consistent.

Therefore $\mathfrak{B}^{\Gamma p}$ is a modified Venn model (which is also a
modified c.o.m) of $\Gamma $. By proposition 7.7 it is a $pd$-model of $%
\Gamma $, from which it may be easily seen that it is a $pg$-model of $%
\Delta $.
The proof of the additional result in case $\Delta $ is essentially finite
is almost the same. The only major difference is that $\mu $ may not be
injective. But its restriction to $C_{\Delta }$ is injective, which is
sufficient for our purpose. \hspace{10.8cm} $\square $
\newline

\noindent REMARK 7.9.\textbf{\ }The last theorem shows that remark 7.4
applies to $PF(C)$ after making the obvious modifications.\newline

\textbf{8. Direct completion of direct deduction. }In this section the five
rules of inference given in definition 1.7 are augmented by five more rules,
in order that $\Gamma ^{g}$ may be directly obtained from $\Gamma $ in case $%
\Gamma $ is consistent (cf. Glashoff (2005) where related problems are dealt
with by brute force via a computer program). The additional five rules are:%
\newline

5. $\frac{Iab}{Iba}$ $\quad \qquad (Ic)$\quad \qquad \qquad $~$6. $\frac{%
Iab,Abc}{Iac}$ $\quad ~(Darii)$\newline

7. $\frac{Iab,Ebc}{Oac}$ $\quad ~~(Ferio)$\qquad $\quad ~~$8. $\frac{Oab,Acb%
}{Oac}$ $\quad (Baroco)$\newline

9. $\frac{Aba,Obc}{Oac}$ $\quad ~(Bocardo).$\newline

Taking the ten rules of inference into consideration, the $d^{\prime }$%
-deduction relation ``$\stackrel{d^{\prime }}{\vdash }$'' may be defined
along the lines of the definition of ``$\stackrel{d}{\vdash }$''. Likewise,
all other definitions involving ``$d$'' may be modified in an obvious way to
give corresponding definitions involving ``$d^{\prime }$''.\newline

\noindent PROPOSITION 8.1.

1. $\Gamma ^{d^{\prime }}=\{\sigma \in S(C):\Gamma \stackrel{d^{\prime }}{%
\vdash }\sigma \}$.

2. $C_{\Gamma ^{d^{\prime }}}=C_{\Gamma }$.\newline

\noindent \textit{Proof. }Along the lines of the proofs of the corresponding
results for $d$: Part 1 of proposition 1.11 and lemma 2.12, respectively. \hspace{4.1 cm} $%
\square $
\newline

The next definition and parts 1,2 of the next lemma are essentially due to
Smith (1983).\newline

\noindent DEFINITION 8.2.\textbf{\ }Let $a,b,a^{\prime },b^{\prime }\in C$.
An $a$-$b$ chain is a sequence $<c_{i}>_{i\in n}\in $ $^{n}C$ for some $n\in 
\Bbb{N}^{+}$, such that $c_{o}=a$ and $c_{n-1}=b$. This chain is said to be
a $\Gamma $-chain, or a chain in $\Gamma $, if $\{Ac_{i}c_{i+1}:i\in
n-1\}\subseteq \Gamma $; it is said to be an $<a^{\prime },b^{\prime }>$
chain if there is $i\in n-1$ such that $c_{i}=a^{\prime }$ and $%
c_{i+1}=b^{\prime }$.\newline

\noindent LEMMA 8.3.\textbf{\ }For $a,b\in C$ and $\Gamma \subseteq BN(C)$:%
\newline

1. $\Gamma \stackrel{d}{\vdash }Aab$ \quad iff \quad $\Gamma \stackrel{%
d^{\prime }}{\vdash }Aab$ \quad iff \quad there is an $a$-$b$ chain in $%
\Gamma $.\newline

2. $\Gamma \stackrel{d}{\vdash }Eab$ \quad iff \quad $\Gamma \stackrel{%
d^{\prime }}{\vdash }Eab$ \quad iff there is \quad $Ea^{\prime }b^{\prime
}\in BN(C)$ such that $\{Ea^{\prime }b^{\prime },Eb^{\prime }a^{\prime
}\}\cap \Gamma \neq \phi $ and $\Gamma \stackrel{d}{\vdash }Aaa^{\prime
},Abb^{\prime }$.\newline

3. $\Gamma \stackrel{d}{\vdash }Iab$ \quad iff \quad $Iab\in \Gamma $ or $%
\Gamma \stackrel{d}{\vdash }Aba$.\newline

3$^{\prime }$. $\Gamma \stackrel{d^{\prime }}{\vdash }Iab$ \quad iff \quad
for some $a^{\prime },b^{\prime }\in C$, \quad $\Gamma \stackrel{d^{\prime }%
}{\vdash }Aa^{\prime }a,Aa^{\prime }b$ \quad or \newline
$\{Ia^{\prime }b^{\prime },Ib^{\prime }a^{\prime }\}\cap \Gamma \neq \phi $
and $\Gamma \stackrel{d^{\prime }}{\vdash }Aa^{\prime }a,Ab^{\prime }b$.%
\newline

4. $\Gamma \stackrel{d}{\vdash }Oab$ \quad iff \quad $Oab\in \Gamma $.%
\newline

4$^{\prime }$. $\Gamma \stackrel{d^{\prime }}{\vdash }Oab$ \quad iff \quad
for some $a^{\prime },b^{\prime }\in C$, \quad $\Gamma \stackrel{d^{\prime }%
}{\vdash }Ia^{\prime }a,Ea^{\prime }b$ \quad or\newline
$Oa^{\prime }b^{\prime }\in \Gamma $ and $\Gamma \stackrel{d^{\prime }}{%
\vdash }Aa^{\prime }a,Abb^{\prime }$.\newline

\noindent \textit{Proof.}\newline

1. If is easy to show that the first statement implies the second. By
induction it may be shown that the second statement implies the third. Again
by induction it may be shown that the third statement implies the first.%
\newline

2. It is easy to show that the first statement implies the second and that
the third implies the first. By induction it may be shown that the second
statement implies the third.\newline

Parts 3 and 4 are easy. In each of the parts 3$^{\prime }$ and 4$^{\prime }$
one direction is easy, the other may be shown by induction. \hspace{6.4cm}$\square $\newline

\noindent PROPOSITION 8.4.\textbf{\ }For $\Gamma \cup \{\sigma \}\subseteq
BN(C)$:\newline

1. If $\Gamma \stackrel{d}{\vdash }\sigma $ then $\Gamma \stackrel{d^{\prime
}}{\vdash }\sigma $.\newline

2. If $\Gamma \stackrel{d^{\prime }}{\vdash }\sigma $ then $\Gamma \stackrel{%
g}{\vdash }\sigma $.\newline

3. $\Gamma $ is $g$-consistent iff $\Gamma $ is $d^{\prime }$-consistent iff 
$\Gamma $ is $d$-consistent.

(So for $e\in \{d,d^{\prime },g\}$ the prefix ``$e$-'' may be deleted from ``%
$e$-consistent'', ``$e$-inconsistent'' and ``$e$-contradictory'').\newline

\noindent \textit{Proof.} Part 1 is obvious, and part 3 is an easy
consequence of parts 1 and 2 above and part 4 of proposition 1.11.

Part 2 is immediate if $\Gamma $ is $d$-inconsistent. To complete the proof
assume that $\Gamma $ is $d$-consistent and proceed by course of values
induction. Let $\Gamma \vdash \sigma $ and let $<\rho _{i}>_{i\in n}$ be a $%
d^{\prime }$-deduction of $\sigma $ from $\Gamma $. If the annotation of $%
\rho _{n-1}$ ($=\sigma $) is that it belongs to $\Gamma $ or that it is
the consequent of a $d$-rule whose premises are previous sentences, the result
easily follows.

It remains to assume that the annotation of $\rho _{n-1}$ is that it is the
consequent of a new rule. The completion of the proof depends on the
specific rule in use. Following is a proof in the case of Darii. The other
cases are similar or easier.

Let $\rho _{n-1}=Iac$ and let its annotation be that it follows from $%
Iab,Abc $ by Darii. By the induction hypothesis $\Gamma \stackrel{g}{\vdash }%
Iab,Abc$. By part 3 of proposition 1.11, $\Gamma \stackrel{d}{\vdash }Abc$,
and by the definition of $\stackrel{g}{\vdash }$, there is $\eta \in BN(C)$
such that $\Gamma ,Eab\stackrel{d}{\vdash }\eta ,\widehat{\eta }$. Since $%
\Gamma $ is $d$-consistent then, in view of parts 1 and 4 of lemma 8.3,
there are $c^{\prime },c^{\prime \prime }\in C$ such that $\{\eta ,\widehat{%
\eta }\}=\{Ic^{\prime }c^{\prime \prime },Ec^{\prime }c^{\prime \prime }\}$.
By lemma 8.3, $\Gamma \stackrel{d}{\vdash }Ic^{\prime }c^{\prime \prime }$
and there is $Ea^{\prime }b^{\prime }\in BN(C)$ such that $\{Ea^{\prime
}b^{\prime },Eb^{\prime }a^{\prime }\}\cap (\Gamma \cup \{Eab\})\neq \phi $
and $\Gamma \cup \{Eab\}\stackrel{d}{\vdash }Ac^{\prime }a^{\prime
},Ac^{\prime \prime }b^{\prime }$, hence $\Gamma \stackrel{d}{\vdash }%
Ac^{\prime }a^{\prime },Ac^{\prime \prime }b^{\prime }$. In view of the $d$%
-consistency of $\Gamma $, lemma 8.3 implies that $Eab\in \{Ea^{\prime
}b^{\prime },Eb^{\prime }a^{\prime }\}$. Let $Eab=Ea^{\prime }b^{\prime }$
(the other case is similar), then $\Gamma \stackrel{d}{\vdash }Ac^{\prime
}a,Ac^{\prime \prime }b$. But $\Gamma \stackrel{d}{\vdash }Abc$, then $%
\Gamma ,Eac\stackrel{d}{\vdash }Ec^{\prime }c^{\prime \prime }$, hence $%
\Gamma \stackrel{g}{\vdash }Iac$. \hspace{6.6cm} $\square $
\newline

In view of the $g$-deduction completeness with respect to the class of Venn
models, part 2 of the above proposition is an immediate consequence of:%
\newline

\noindent PROPOSITION 8.5.\textbf{\ }The $d^{\prime }$-deduction is sound
with respect to the class of Venn models (hence with respect to the class of
order models).\newline

\noindent \textit{Proof.} Routine. \hspace{10.1cm} $\square $
\newline

\noindent REMARK 8.6.\textbf{\ }The converse of part 2 of proposition 8.4
does not always hold. For if $\Gamma $ is inconsistent then $C_{\Gamma
^{g}}=C$, while it is easy to find an inconsistent $\Gamma $ such that $%
C_{\Gamma ^{d^{\prime }}}=C_{\Gamma }\neq C$. Also the weaker statement: $%
\Gamma ^{d^{\prime }}\cap BN(C_{\Gamma })=\Gamma ^{g}\cap BN(C_{\Gamma })$,
does not always hold. A counter example is $\Gamma =\{Aab,Oab\}$.

The consistency of $\Gamma $ solves the problem as the following theorem
shows (cf. Smith 1983).\newline

\noindent THEOREM 8.7.\textbf{\ }For consistent $\Gamma $, $\quad \Gamma
^{d^{\prime }}=\Gamma ^{g}$.\newline

\noindent \textit{Proof.} The inclusion of $\Gamma ^{d^{\prime }}$ in $%
\Gamma ^{g}$ is guaranteed by part 2 of proposition 8.4. For the other
direction assume that $\Gamma $ is consistent and $\Gamma \stackrel{g}{%
\vdash }\sigma $. If $\sigma $ is universal the result follows by part 3 of
proposition 1.11 and part 1 of proposition 8.4. So it remains to deal with
the particulars. The consistency of $\Gamma $ restricts what to be considered
to the following:
\newline

Case 1. $\sigma $ is $Iab$ for some $a,b\in C$. By the method made use of in
the proof of part 2 proposition 8.4, consideration may be restricted to the
following subcase only. There are $c,c^{\prime }\in C$ such that $\Gamma 
\stackrel{d}{\vdash }Icc^{\prime },Aca,Ac^{\prime }b$, which implies that $%
\Gamma \stackrel{d^{\prime }}{\vdash }Iab$.\newline

Case 2. $\sigma $ is $Oab$ for some $a,b\in C$. In this case $\Gamma ,Aab%
\stackrel{d}{\vdash }\rho ,\widehat{\rho }$ for some $\rho \in BN(C)$.
Distinguish between two subcases.

Subcase 2.1. For some $c,c^{\prime }\in C$, $\{\rho ,\widehat{\rho }%
\}=\{Acc^{\prime },Occ^{\prime }\}$. Then $\Gamma \stackrel{d}{\vdash }%
Aca,Abc^{\prime },Occ^{\prime }$ which implies that $\Gamma \stackrel{%
d^{\prime }}{\vdash }Oab$.\newline

Subcase 2.2. For some $c,c^{\prime }\in C$, \quad $\{\rho ,\widehat{\rho }%
\}=\{Ecc^{\prime },Icc^{\prime }\}$. This subcase may be divided into the
following three subsubcases.

Subsubcase 2.2.1. $\Gamma \stackrel{d}{\vdash }Ecc^{\prime },\Gamma 
\stackrel{d}{\nvdash }Icc^{\prime }$. Then $\Gamma \stackrel{d}{\vdash }%
Ac^{\prime }a,Abc,Ecc^{\prime }$ which implies that $\Gamma \stackrel{%
d^{\prime }}{\vdash }Oab$.

Subsubcase 2.2.2. $\Gamma \stackrel{d}{\nvdash }Ecc^{\prime }$ and $\Gamma 
\stackrel{d}{\vdash }Icc^{\prime }$. Then there is $Ea^{\prime }b^{\prime
}\in BN(C)$ such that $\Gamma \stackrel{d}{\vdash }Ea^{\prime }b^{\prime }$
and $\Gamma ,Aab\stackrel{d}{\vdash }Aca^{\prime },Ac^{\prime }b^{\prime }$;
while $\Gamma \stackrel{d}{\nvdash }Aca^{\prime }$ or $\Gamma \stackrel{d}{%
\nvdash }Ac^{\prime }b^{\prime }$, but -by the consistency of $\Gamma $- not
both.

This subsubcase may be further divided into two subsubsubcases.

Subsubsubcase 2.2.2.1. $\Gamma \stackrel{d}{\nvdash }Aca^{\prime }$ but $%
\Gamma \stackrel{d}{\vdash }Ac^{\prime }b^{\prime }$. Then $\Gamma \stackrel{%
d}{\vdash }Icc^{\prime },Ea^{\prime }b^{\prime },Ac^{\prime }b^{\prime },$ $%
Aca,Aba^{\prime }$ from which $\Gamma \stackrel{d^{\prime }}{\vdash }Oab$
follows.

Subsubsubcase 2.2.2.2. $\Gamma \stackrel{d}{\vdash }Aca^{\prime }$ but $%
\Gamma \stackrel{d}{\nvdash }Ac^{\prime }b^{\prime }$. Similar to
subsubsubcase 2.2.2.1.

Subsubcase 2.2.3. $\Gamma \stackrel{d}{\nvdash }Ecc^{\prime }$ and $\Gamma 
\stackrel{d}{\nvdash }Icc^{\prime }$. Then there is $Ea^{\prime }b^{\prime
}\in BN(C)$ such that $\Gamma \stackrel{d}{\vdash }Ac^{\prime
}a,Abc,Ea^{\prime }b^{\prime }$ and $\Gamma ,Aab\stackrel{d}{\vdash }%
Aca^{\prime },Ac^{\prime }b^{\prime }$, while $\Gamma \stackrel{d}{\nvdash }%
Aca^{\prime }$ or $\Gamma \stackrel{d}{\nvdash }Ac^{\prime }b^{\prime }$.
But the consistency of $\Gamma $ implies that $\Gamma \stackrel{d}{\vdash }%
Ac^{\prime }b^{\prime }$, then $\Gamma \stackrel{d}{\nvdash }Aca^{\prime }$,
then $\Gamma \stackrel{d}{\vdash }Aca,Aba^{\prime }$. In particular, $\Gamma 
\stackrel{d}{\vdash }Ac^{\prime }a,Ac^{\prime }b^{\prime },Aba^{\prime
},Ea^{\prime }b^{\prime }$, hence the result. \hspace{12.1cm} $\square $
\newline

\noindent REMARK 8.8.\textbf{\ }In a different context, Smith (1983):

1. Excluded subcase 2.1 under the claim that it is impossible that $\Gamma 
\stackrel{d}{\vdash }Occ^{\prime }$.

2. Subsubcase 2.2.3 was deemed to be impossible.\newline

\textbf{9. Models of }$NF(C)$\textbf{\ revisited. }An $NF(C)$-structure $%
\mathfrak{B}$ is said to be a $d^{\prime }$-model if for every $\Gamma \cup
\{\sigma \}\subseteq BN(C)$, $\Gamma \vDash _{\mathfrak{B}}\sigma $ whenever $%
\Gamma \stackrel{d^{\prime }}{\vdash }\sigma $.

An immediate consequence of this definition is:\newline

\noindent PROPOSITION 9.1.\textbf{\ }An $NF(C)$-structure \noindent $\mathfrak{B}%
=<B,A^{*},E^{*},I^{*},O^{*},\mu >$ is a $d^{\prime }$-model iff it is a $d$%
-model (hence satisfying conditions 1-4 of proposition 2.3) and:

5. $(I^{*\mu }|A^{*\mu })\subseteq I^{*\mu }\subseteq \stackrel{\smallsmile 
}{I^{*\mu }}$.

6. $(I^{*\mu }|E^{*\mu })\cup (O^{*\mu }|\stackrel{\smallsmile }{A^{*\mu }}%
)\cup (\stackrel{\smallsmile }{A^{*\mu }}|O^{*\mu })\subseteq O^{*\mu }$. \hspace{4.6cm} $%
\square $
\newline

Along the lines of the proofs of lemma 2.5, theorem 2.6 and theorem 2.10,
the following may be proved:
\newline

\noindent THEOREM 9.2.\textbf{\ }For every $\Gamma \cup \{\sigma \}\subseteq
BN(C)$:\newline

1. $\mathfrak{B}_{\Gamma ^{d^{\prime }}}$ is a $d^{\prime }$-model (of $\Gamma
^{d^{\prime }}$, hence of $\Gamma $).\newline

2. $d^{\prime }$-deduction is sound and complete with respect to the class
of $d^{\prime }$-models. That is $\Gamma \stackrel{d^{\prime }}{\vdash }%
\sigma $ iff $\Gamma \stackrel{d^{\prime }}{\vDash }\sigma $.\newline

3. $\Gamma \stackrel{d^{\prime }}{\vDash }\sigma $ iff $\Gamma _{1}\stackrel{%
d^{\prime }}{\vDash }\sigma $ for some finite $\Gamma _{1}\subseteq \Gamma $.

(This is called $d^{\prime }$-compactness). \hspace{7cm} $\square $
\newline

\noindent REMARK 9.3.\textbf{\ }All remarks given in remarks and definitions
2.11 hold with ``$d^{\prime }$'' replacing ``$d$''. All proofs of the
original versions essentially go through; the only exception is the first
remark, whose modified version may be proved by part 2 of proposition 8.4.%
\newline

\noindent THEOREM 9.4.\textbf{\ }An $NF(C)$-structure $\mathfrak{B}%
(=<B,A^{*},E^{*},I^{*},O^{*},\mu >)$ is a $g$-model iff it is a $d^{\prime }$%
-model and:\newline

1. $A^{*\mu }\cap O^{*\mu }=\phi =E^{*\mu }\cap I^{*\mu }$, or\newline

2. $A^{*\mu }=E^{*\mu }=I^{*\mu }=O^{*\mu }=\mu (C)\times \mu (C)$.\newline

\noindent \textit{Proof.} Only if: By part 2 of proposition 8.4 and an
obvious generalization of part 3 of remarks and definitions 2.11.

If: Every $NF(C)$-structure which satisfies condition 2 is a $g$-model. So,
assume that $\mathfrak{B}$ is a $d^{\prime }$-model which satisfies condition 1.
To see that it is a $g$-model, let $\Gamma \cup \{\sigma \}\subseteq BN(C)$, 
$\Gamma \stackrel{g}{\vdash }\sigma $ and $\mathfrak{B}\vDash \Gamma $. By
remark 9.3, $\Gamma $ is consistent, hence by theorem 8.7, $\Gamma \stackrel{d^{\prime}} { \vdash }
\sigma $, hence $\mathfrak{B}\vDash \sigma $. \hspace{4cm} $\square $
\newline

Theorem 9.4 fully characterizes the class of $g$-models, as was promised
after the proof of theorem 2.10.\newline

\noindent DEFINITIONS and remarks 9.5.\textbf{\ }\newline

1. For an $NF(C)$-structure $\mathfrak{B}$ and a relation symbol $W\in
\{A,E,I,O\}$, define $Bt^{W}\mathfrak{B}$ (the basic $W$-theory of $\mathfrak{B}$), $%
Bt^{+}\mathfrak{B}$ (the basic positive theory of $\mathfrak{B}$), $Bt^{-}\mathfrak{B}$
(the basic negative theory of $\mathfrak{B}$) and $Bt\mathfrak{B}$ (the basic theory
of $\mathfrak{B}$) as follows:\newline

$Bt^{W}\mathfrak{B}=\{Wab\in BN(C):\mathfrak{B}\vDash Wab\}$.\newline

$Bt^{+}\mathfrak{B}=Bt^{A}\mathfrak{B}\cup Bt^{I}\mathfrak{B}$.\newline

$Bt^{-}\mathfrak{B}=Bt^{E}\mathfrak{B}\cup Bt^{O}\mathfrak{B}$.\newline

$Bt\mathfrak{B}=Bt^{+}\mathfrak{B}\cup Bt^{-}\mathfrak{B}$.\newline

So two $NF(C)$-structures are $B$-equivalent iff they have the same basic
theory.\newline

For $i\in 2$ let $\mathfrak{B}_{i}(=<B_{i},A_{i},E_{i},I_{i},O_{i},\mu _{i}>)$
be an $NF(C)$-structure.\newline

2. $\mathfrak{B}_{o}$ is said to be a substructure of $\mathfrak{B}_{1}$ and $\mathfrak{B}_{1}$ is said to be a superstructure of $\mathfrak{B}_{o}$ if $B_{o}\subseteq
B_{1},\mu _{o}=\mu _{1}$ and for every $W\in \{A,E,I,O\}$, $W_{o}=W_{1}\cap
(B_{o}\times B_{o})$. If, morever, $B_{o}=$ Range $\mu _{1}$ $(=$ Range $\mu
_{o})$, $\mathfrak{B}_{o}$ is said to be a core substructure of $\mathfrak{B}_{1}$.
Obviously each $NF(C)$-structure has a unique core substructure, to be
called its core substructure. $\mathfrak{B}_{o}$ is a core substructure of some $%
NF(C)$-structure iff it is the core substructure of itself iff $B_{o}=$
Range $\mu _{o}$. In this case $\mathfrak{B}_{o}$ is said to be a core
structure. Obviously every canonical structure is a core structure.

$\mathfrak{B}_{o}$, $\mathfrak{B}_{1}$ have the same core substructure iff $\mu
_{o}=\mu _{1}$ and $Bt\mathfrak{B}_{o}=Bt\mathfrak{B}_{1}$.\newline

3. If $\mathfrak{B}_{o}$ is a substructure of \noindent $\mathfrak{B}_{1}$ then they
have the same core substructure and the three structures have the same basic
theory. Hence for $e\in \{d,d^{\prime },g\}$ if one of them is an $e$-model,
so also are the other two.

In this case $\mathfrak{B}_{o}$ is said to be an $e$-submodel of $\mathfrak{B}_{1}$,
and $\mathfrak{B}_{1}$ is said to be an $e$-supermodel of $\mathfrak{B}_{o}$; and
the core substructure is said also to be the core $e$-submodel. If a core
structure is an $e$-model, it is said to be a core $e$-model.\newline

4. $\mathfrak{B}_{o}$ is said to be a positive semisubstructure of $\mathfrak{B}_{1}$
and $\mathfrak{B}_{1}$ is said to be a positive semisuperstructure of $\mathfrak{B}%
_{o}$ if $B_{o}\subseteq B_{1}$, $\mu _{o}=\mu _{1}$ and:\newline

$W_{o}=W_{1}\cap (B_{o}\times B_{o})$ \qquad for every $W\in \{A,I\}$,

$W_{o}\subseteq W_{1}\cap (B_{o}\times B_{o})$ \qquad for every $W\in
\{E,O\} $.\newline

In this case $Bt^{+}\mathfrak{B}_{o}=Bt^{+}\mathfrak{B}_{1}$ and $Bt^{-}\mathfrak{B}%
_{o}\subseteq Bt^{-}\mathfrak{B}_{1}$. For each $e\in \{d,d^{\prime },g\}$ if,
in addition, $\mathfrak{B}_{o}$ and $\mathfrak{B}_{1}$ are both $e$-models, it is
said also that $\mathfrak{B}_{o}$ is a positive $e$-semisubmodel of $\mathfrak{B}%
_{1} $ and $\mathfrak{B}_{1}$ is a positive $e$-semisupermodel of $\mathfrak{B}_{o}$.%
\newline

\noindent THEOREM 9.6.\textbf{\ }For each $e\in \{d^{\prime },g\}$ if $f{%
B}_{o}$ (defined as above) is a consistent core $e$-model then there is an
order model $\mathfrak{B}_{1}(=<B_{1},A_{1},\mu _{1}>)$ such that:\newline

1. $\mathfrak{B}_{1}$ is a positive $e$-semisupermodel of $\mathfrak{B}_{o}$.

2. If $A_{o}$ is a partial ordering, then so also is $A_{1}$.

3. If $\mathfrak{B}_{o}$ is complete, then it is the core $e$-submodel of $\mathfrak{B}_{1}$.\newline

For $e=d$, the above holds after weakening part 1 to become:

1$^{\prime }$. $B_{o}\subseteq B_{1},\mu _{o}=\mu _{1},Bt\mathfrak{B}%
_{o}\subseteq Bt\mathfrak{B}_{1}$, and $A_{o}=A_{1}\cap (B_{o}\times B_{o})$;
hence $Bt^{A}\mathfrak{B}_{o}=Bt^{A}\mathfrak{B}_{1}$.\newline

\noindent \textit{Proof.} Let $e\in \{d,d^{\prime },g\}$ and let $\mathfrak{B}%
_{o}$ be a consistent core $e$-model. Put:

$B^{\prime }=\{\{a_{o},a_{1}\}\subseteq B_{o}:<a_{o},a_{1}>\in I_{o}$ or $%
<a_{1},a_{o}>\in I_{o}$, and $\{a_{o},a_{1}\}$ has no $A_{o}$-lower bound$\}$%
,

$B_{1}=B_{o}\cup B^{\prime }$ ($B_{o},B^{\prime }$ may be assumed 
disjoint),

$A_{1}=A_{o}\cup \mathfrak{l}_{B^{\prime }}\cup \{<\{a_{o},a_{1}\},a_{2}>\in
B^{\prime }\times B_{o}:<a_{o},a_{2}>\in A_{o}$ or $<a_{1},a_{2}>\in A_{o}\}$%
,

$\mu _{1}=\mu _{o}$\newline

$A_{1}$ is reflexive on $B_{1}$ since $A_{o}$ is reflexive on $B_{o}$. To
prove the transitivity of $A_{1}$, let $<b_{o},b_{1}>,<b_{1},b_{2}>\in A_{1}$%
. If $<b_{o},b_{1}>$ or $<b_{1},b_{2}>$ belongs to $\mathfrak{l}_{B^{\prime }}$
then $<b_{o},b_{2}>\in A_{1}$, else $<b_{1},b_{2}>\in A_{o}$. If $%
<b_{o},b_{1}>\in A_{o}$ then $<b_{o},b_{2}>\in A_{1}$. It remains to
consider the case where $b_{o}=\{a_{o},a_{1}\}$ for some $\{a_{o},a_{1}\}\in
B^{\prime }$ such that $<a_{o},b_{1}>\in A_{o}$ or $<a_{1},b_{1}>\in A_{o}$,
in both cases $<b_{o},b_{2}>\in A_{1}$. So $A_{1}$ is transitive. Hence $%
<B_{1},A_{1},\mu _{1}>$ is an order model, which is to be denoted by ``$%
\mathfrak{B}_{1}$''.\newline

To prove part 2 it suffices to notice that if $<b_{o},b_{1}>,<b_{1},b_{o}>%
\in A_{1}$ then they both belong to $A_{o}$ or both belong to $\mathfrak{l}
_{B^{\prime }}$.\newline

To prove parts 1, 1$^{\prime }$ notice that $B_{o}\subseteq B_{1}$ and, by
the disjointness of $B_{o},B^{\prime }$, $A_{o}=A_{1}\cap (B_{o}\times
B_{o}) $. Let $<a_{o},a_{1}>\in I_{o}$. If $\{a_{o},a_{1}\}$ has an $A_{o}$%
-lower bound then it is an $A_{1}$-lower bound, else the element $%
\{a_{o},a_{1}\}\in B_{1}$ is an $A_{1}$-lower bound of the subset $%
\{a_{o},a_{1}\}\subseteq B_{1}$. In both cases $<a_{o},a_{1}>\in I_{1}$,
hence $I_{o}\subseteq I_{1}\cap (B_{o}\times B_{o})$.\newline

At this point the proof forks into two branches:
\newline

(i) Assume $e\in \{d^{\prime },g\}$ and let $<a_{o},a_{1}>\in I_{1}\cap
(B_{o}\times B_{o})$. To show that $<a_{o},a_{1}>\in I_{o}$ several cases
have to be considered, following is one of them, the others are similar or
easier.

There is $<a_{2},a_{3}>\in I_{o}$ such that $<a_{2},a_{o}>,<a_{3},a_{1}>\in
A_{o}$. Since $B_{o}=$ Range $\mu _{o}$ then, by theorem 9.4 and part 5 of
proposition 9.1, $<a_{o},a_{1}>\in I_{o}$. So $I_{1}\cap (B_{o}\times
B_{o})\subseteq I_{o}$. Hence $I_{o}=I_{1}\cap (B_{o}\times B_{o})$.

That $E_{o}\subseteq E_{1}\cap (B_{o}\times B_{o})$ and $O_{o}\subseteq
O_{1}\cap (B_{o}\times B_{o})$ is guaranteed by the consistency of $\mathfrak{B}%
_{o}$. This completes the proof of 1.\newline

(ii) The other branch is $e=d$. To show that $E_{o}\subseteq E_{1}$ assume
that there is $<a_{o},a_{1}>\in (E_{o}-E_{1})$, then $<a_{o},a_{1}>\in I_{1}$%
, then $\{a_{o},a_{1}\}$ has an $A_{1}$-lower bound. To show that this is
absurd, several cases have to be considered; following is one of them, the
others are easier or similar.

There is $<a_{2},a_{3}>\in I_{o}$ such that $<a_{3},a_{o}>,<a_{2},a_{1}>\in
A_{o}$. Since $B_{o}=$ Range $\mu _{o}$ then, by parts 3, 4 of proposition
2.3, $<a_{2},a_{3}>\in E_{o}$ which contradicts the consistency of $\mathfrak{B}%
_{o}$.

That $O_{o}\subseteq O_{1}$ is guaranteed by the consistency of $\mathfrak{B}%
_{o} $, since $A_{o}=A_{1}\cap (B_{o}\times B_{o})$. This completes the
proof of 1$^{\prime }$ and ends the forkation.\newline

For $e\in \{d,d^{\prime },g\}$, if $\mathfrak{B}_{o}$ is complete then ``$%
\subseteq $'' may be replaced by ``$=$'' at the appropriate places, which
proves part 3. \hspace{5.1cm} $\square $
\newline

Taking the relationship between the $e$-models ($e\in \{d,d^{\prime },g\}$)
and their respective core $e$-submodels into consideration, a weaker result,
which holds for a wider class of $e$-models, immediately follows:
\newline

\noindent COROLLARY 9.7.\textbf{\ }For $e\in \{d,d^{\prime },g\}$, if $\mathfrak{B}$ is an $e$-model whose core $e$-submodel is consistent, then there is an
order model $\mathfrak{B}^{\prime }$ such that \newline
$Bt\mathfrak{B}\subseteq Bt\mathfrak{B}^{\prime }$. Moreover, 
\begin{eqnarray*}
Bt^{+}\mathfrak{B} &=&Bt^{+}\mathfrak{B}^{\prime }\text{ \qquad if }e\in \{d^{\prime
},g\}\text{,} \\
Bt^{A}\mathfrak{B} &=&Bt^{A}\mathfrak{B}^{\prime }\text{ \qquad if }e=d\text{.} \hspace{7.8cm} 
\square
\end{eqnarray*}
\newline

In view of the last part of subsection 2.5, the above corollary may be
immediately strengthened as follows: 
\newline

\noindent COROLLARY 9.8.\textbf{\ }In the above corollary ``an order model''
may be replaced by ``a partial order model which is a c.o.m and a Venn model
at the same time''. \hspace{12cm} $\square $
\newline

Part 4 of theorem 1.12 may be extended to the case $e=d^{\prime }$, to
get a result similar to that obtained there for the case $e=g$; the
result obtained (there) for the case $e=d$ is weaker. Call the collection of
these three results ``syntactical congruence''.

Syntactical congruence together with the definitions of core $e$-models ($%
e\in \{d,d^{\prime },g\}$) yield semantical congruence as formulated by
parts 1 and 2 of the next theorem. Part 3 of the same theorem (whose proof
is straightforward) strengthens the conclusion of part 2, under some
additional condition. Alternatively, semantical congruence may be directly
proved by the characterizations of $e$-models ($e\in \{d,d^{\prime },g\}$)
given in propositions 2.3 and 9.1 and theorem 9.4.\newline

\noindent THEOREM 9.9.\textbf{\ }Let $e\in \{d,d^{\prime },g\}$ and let $%
\mathfrak{B}=<B,A^{*},E^{*},I^{*},O^{*},\mu >$ be a core $e$-model (consistent
or not). Put $\sim $ $=A^{*}\cap \stackrel{\smallsmile }{A^{*}}$, then:

1. $\sim $ is a congruence relation on $<B,A^{*},E^{*},\mu >$ and $%
A^{*}/\sim $ is a partial ordering on $B/\sim $.

Moreover, for $e\in \{d^{\prime },g\}$:

2. $\sim $ is a congruence relation on $\mathfrak{B}$. The mapping $b\longmapsto
b/\sim $ is an epimorphism from $\mathfrak{B}$ onto $\mathfrak{B}/\sim $. Hence $%
\mathfrak{B}/\sim $ is a core $e$-model which is basically equivalent to $\mathfrak{B}$.

3. If $\mathfrak{B}$ is, in addition, an order model, then $\mathfrak{B}/\sim $ is
also a partial order model. \hspace{11.8cm} $\square$
\newline

For $e\in \{d^{\prime },g\}$, semantical congruence makes it possible to
replace ``$\mathfrak{B}_{o}$'' in theorem 9.6 by ``$\mathfrak{B}_{o}/\sim $''. This
provides, for $e\in \{d^{\prime },g\}$, an alternative proof of a weaker
form of corollary 9.8, where the partial order model may be neither concrete
nor Venn.

The corresponding weaker result for the case $e=d$ may likewise be obtained,
but the alternative proof is a bit more involved.\newline

\noindent REMARKS and definitions 9.10.\textbf{\ }\newline

1. Theorem 9.6 (or corollary 9.7) and corollary 9.8 (or its weaker forms)
provide, respectively, direct ways to order models and partial order models
for consistent $\Gamma (\subseteq BN(C))$. Simply in each of them let the
core $e$-model be the canonical structure $\mathfrak{B}_{\Gamma ^{e}}(e\in
\{d,d^{\prime },g\})$. In the case of corollary 9.8 the partial order model
may be required to be a concrete order model and a Venn model at the same
time.\newline

2. Let $e\in \{d,d^{\prime },g\}$ and let $\mathfrak{C}$ be a class of $NF(C)$%
-structures, then:

\quad 1. $e$ is said to be $\mathfrak{C}$-strongly semantically complete if for
every $\Gamma \subseteq BN(C)$ there is $\mathfrak{B}\in \mathfrak{C}$ such that $Bt%
\mathfrak{B}=\Gamma ^{e}$.

\quad 2. $e$ is said to be $\mathfrak{C}$-syntactically complete if for every $%
\Gamma \cup \{\sigma \}\subseteq BN(C)$, $\Gamma \stackrel{e}{\vdash }\sigma 
$ whenever $\Gamma \vDash _{\mathfrak{C}}\sigma $.

\quad 3. $e$ is said to be $\mathfrak{C}$-consistently syntactically complete if
for every $e$-consistent $\Gamma \subseteq BN(C)$ and every $\sigma \in
BN(C) $, $\Gamma \stackrel{e}{\vdash }\sigma $ whenever $\Gamma \vDash _{%
\mathfrak{C}}\sigma $.

\quad 4. $e$ is said to be $\mathfrak{C}$-consistently semantically complete if
every $e$-consistent $\Gamma \subseteq BN(C)$ has a model in $\mathfrak{C}$.

For $i\in \{1,2,3\}$, the condition given in clause $i$ implies the
condition given in clause $i+1$.\newline

3. Put:

$Or$ = the class of all order models,

$Po$ = the class of all partial order models,

$Le$ = the set of all Leibniz models,

$Co$ = the class of all concrete order models,

$Ve$ = the class of all Venn models.

And for $e\in \{d,d^{\prime },g\}$ put:

$Be$ = $\{\mathfrak{B}_{\Gamma ^{e}}:\Gamma \subseteq BN(C)\}$.

Also put:

$M=\{Or,Po,Le,Co,Ve\}\cup \{Be:e\in \{d,d^{\prime },g\}\}$.\newline

4. $Le\cup Co\subseteq Po\subseteq Or$, $\quad Bg\subseteq Bd^{\prime
}\subseteq Bd$.\newline

5. Every element of $\bigcup M$ is a $d$-model.

\quad Every element of $(\bigcup M-Bd)\cup Bd^{\prime }$ is a $d^{\prime }$%
-model.

\quad Every element of $(\bigcup M-Bd)\cup Bg$ is a $g$-model.\newline

6. For $e\in \{d,d^{\prime },g\}$, $e$ is $Be$-strongly semantically
complete.

\quad For $\mathfrak{C}\in M-\{Le\}$, $d$ (respectively $d^{\prime },g$) is $%
\mathfrak{C}$-consistently semantically (respectively consistently
syntactically, syntactically) complete. If $C$ is finite, the exclusion of $%
Le$ may be dropped.\newline

7. For $e\in \{d,d^{\prime },g\}$ and $\Gamma \subseteq BN(C)$, $\Gamma $ is
said to be $e$-syntactically complete if for every $\sigma \in BN(C)$, $%
\Gamma \stackrel{e}{\vdash }\sigma $ or $\Gamma \stackrel{e}{\vdash }%
\widehat{\sigma }$.\newline

8. For $e\in \{d,d^{\prime },g\}$ and $\mathfrak{C}\in M-\{Le\}$, if $\Gamma $
is consistent and $e$-syntactically complete then there is $\mathfrak{B}\in 
\mathfrak{C}$ such that $Bt\mathfrak{B}=\Gamma ^{e}$. If, moreover, $C_{\Gamma }$ is
finite, the exclusion of $Le$ may be dropped.\newline

\textbf{10. Decidability revisited.}\newline

\noindent THEOREM 10.1.\textbf{\ }For each $e\in \{d,d^{\prime },g\}$ there
is a polynomial (of degree at most 8) time algorithm to decide for any $%
<\Gamma ,\sigma >\in \wp (BN(C))\times BN(C)$ whether $\Gamma \stackrel{e}{%
\vdash }\sigma $, provided that $\Gamma $ is essentially finite and:

$\Gamma \cap (\{Ecc:c\in (C-C_{\Gamma })\}\cup \{Occ:c\in (C-C_{\Gamma
})\})=\phi $.\newline

\noindent \textit{Proof.} For $e=d$ a proof may be obtained by slightly
modifying the appropriate parts of the proof of theorem 3.2.

In view of lemma 8.3, a proof for the case $e=d^{\prime }$ may be obtained
along the same lines as above.

In view of remarks 3.1, the first part of this theorem may be made use of to
determine whether $\Gamma $ is inconsistent. If yes, $\Gamma \stackrel{g}{%
\vdash }\sigma $; else $\Gamma \stackrel{g}{\vdash }\sigma $ iff $\Gamma 
\stackrel{d^{\prime }}{\vdash }\sigma $, by theorem 8.7. \hspace{11.5cm} $\square $
\newline

\textbf{11. Sorites. }Soriteses are well known in Aristotelian syllogistic
(see Hurley, P. J. 1982, p. 201; Rosenthal, M. and Yudin, R (eds.) 1967, p.
423; also cf. Boger, G. 1998, pp. 197-8; Smiley, T.J. 1973, pp. 139-40).

The notion of a sorites may be explicated as follows.\newline

\noindent DEFINITION 11.1.\textbf{\ }Let $e\in \{d,d^{\prime }\}$ and let $%
\Gamma \subseteq BN(C)$. An annotation of an $e$-deduction $<\sigma
_{i}>_{i\in k}$ from $\Gamma $ is said to be an $e$-sorites annotation if
the following conditions are satisfied:

1. $\sigma _{i}\neq \sigma _{j}$ whenever $i\neq j$ $(i,j\in k)$.

2. For $i\in k-1$, $\sigma _{i}$ is involved in the annotation of another
sentence in the following and only in the following way.

2.1. If $1\leq i\leq k-3$ then exactly one of the following holds:

2.1.1. $\sigma _{i+1}$ is annotated as the consequent of $\sigma _{i}$ by
some $e$-rule with one premise.

2.1.2. $\sigma _{i+1}$ is annotated as the consequent of $\sigma _{i-1}$, $%
\sigma _{i}$ or $\sigma _{i}$, $\sigma _{i-1}$ by some $e$-rule with two
premises.

2.1.3. $\sigma _{i+2}$ is annotated as the consequent of $\sigma _{i}$, $%
\sigma _{i+1}$ or $\sigma _{i+1}$, $\sigma _{i}$ by some $e$-rule with two
premises.

2.2. If $1\leq i=k-2$ then exactly one of 2.1.1 and 2.1.2 holds.

2.3. If $i=0$ then exactly one of the following holds:

2.3.1. $k=2$ and 2.1.1 holds.

2.3.2. $k>2$ and exactly one of 2.1.1 and 2.1.3 holds.\newline

An $e$-sorites from $\Gamma $ is an $e$-deduction from $\Gamma $ which
admits a sorites annotation. An $e$-sorites of $\sigma (\in BN(C))$ from $%
\Gamma $ is an $e$-deduction of $\sigma $ from $\Gamma $ which is an $e$%
-sorites. In case there is such a sorites, we write ``$\Gamma \stackrel{es}{%
\vdash }\sigma $''.

Condition 2 of the above definition entails that, with the exception of the
last sentence, every sentence occurring in an $e$-sorites from $\Gamma $ is
made use of exactly once as a premise of some application of some $e$-rule,
and in this (hence in each) application the premise or the premises
immediately precede the consequent.

For $e\in \{d,d^{\prime }\}$ there is, obviously, a set $\Gamma \subseteq
BN(C)$ and an $e$-deduction from $\Gamma $ which is not an $e$-sorites from $%
\Gamma $. So the best we may hope for is to find an $e$-sorites of $\sigma $
from $\Gamma$, for every $\Gamma \cup \{\sigma \}\subseteq BN(C)$ such that $%
\Gamma \stackrel{e}{\vdash }\sigma $. Even this is not always attainable.

Let $\Gamma _{o}=\{Aca,Ebc\}$ and $\Gamma _{1}=\{Acx,Ebx,Ica\}$, then for $%
i\in 2$, $\Gamma _{i}$ is consistent and $\Gamma _{i}\stackrel{d^{\prime }}{%
\vdash }Oab$, but not $\Gamma _{i}\stackrel{d^{\prime }s}{\vdash }Oab$. For $%
i=0$, adding the rule $\frac{Eab}{Oab}$ ($E$-sub) as an additional rule of
inference will solve the problem. Same holds for $i=1$ if, instead, $\frac{%
Iba,Ebc}{Oac}$ (Ferison) is added.\newline

\textbf{11.1. Further extension of direct deduction.} Taking into
consideration the following two rules of inference.

10. $E$-sub \qquad \qquad 11. Ferison

\noindent in addition to the ten rules of inference of $d^{\prime }$, the $%
d^{\prime \prime }$-deduction relation ``$\stackrel{d^{\prime \prime }}{%
\vdash }$'' may be defined along the lines of the definition of ``$\stackrel{%
d}{\vdash }$''. Likewise all the other definitions involving ``$d$'' may be
modified in an obvious way to give corresponding definitions involving ``$%
d^{\prime \prime }$''.
\newline

\noindent PROPOSITION 11.2.\textbf{\ }For $\Gamma \cup \{\sigma \}\subseteq
BN(C)$, 
\[
\Gamma \stackrel{d^{\prime }}{\vdash }\sigma \text{ \qquad iff \qquad }%
\Gamma \stackrel{d^{\prime \prime }}{\vdash }\sigma . 
\]

\noindent \textit{Proof.} One direction is obvious, the other is easy. \hspace{4.5cm} $%
\square $
\newline

\noindent DEFINITION/remark 11.3.\textbf{\ }The $d^{\prime \prime }$-models
may be defined along the lines of the definition of the $d^{\prime }$-models.

By the above proposition they are the same.\newline

\noindent THEOREM 11.4.\textbf{\ }Let $\Gamma \cup \{\sigma \}\subseteq
BN(C) $ and $e\in \{d,d^{\prime },d^{\prime \prime }\}$ then: 
\[
\Gamma \stackrel{es}{\vdash }\sigma \text{ \qquad whenever \qquad }\Gamma 
\stackrel{e}{\vdash }\sigma 
\]
provided one of the following conditions holds:

1. $\sigma $ is affirmative,

2. $\sigma $ is universal negative and $\Gamma $ is consistent,

3. $\sigma $ is particular negative, $e=d$ or $\Gamma $ is consistent and $%
e=d^{\prime \prime }$.

The other direction unconditionally holds, so the two sides are equivalent
if $\Gamma $ is consistent and $e\in \{d,d^{\prime \prime }\}$.\newline

\noindent \textit{Proof.} Assume $\Gamma \stackrel{e}{\vdash }\sigma $.
Distinguish between the following cases.

1. $\sigma =Aab$, for some $a,b\in C$. By proposition 11.2 and lemma 8.3
there is an $a$-$b$ chain in $\Gamma $, $<c_{i}>_{i\in n}$ say. We may
assume that this chain is injective. If $n=1$, then there is an $e$-sorites
of $\sigma $ from $\Gamma $ of length 1. Else $n\geq 2$; define $<\rho
_{i}>_{i\in 2n-3}$ as follows: 
\begin{eqnarray*}
\rho _{2j} &=&Ac_{o}c_{j+1}=Aac_{j+1}\text{ \qquad ~~}j\in n-1, \\
\rho _{2j+1} &=&Ac_{j+1}c_{j+2}\text{ \qquad \qquad \qquad }j\in n-2.
\end{eqnarray*}

Then $<\rho _{i}>_{i\in 2n-3}$ is an $e$-sorites of $\sigma $ from $\Gamma $.%
\newline

2. $\sigma =Iab$, for some $a,b\in C$. If $e=d$, the result is an easy
consequence of part 1 of this proof and lemma 8.3. Else, by proposition 11.2
we may assume that $e=d^{\prime }$. By lemma 8.3 it suffices to deal with
the following three subcases (for some $a^{\prime },b^{\prime }\in C$):

2.1. $Ia^{\prime }b^{\prime }\in \Gamma $ and $\Gamma \stackrel{d}{\vdash }%
Aa^{\prime }a,Ab^{\prime }b$,

2.2. $Ib^{\prime }a^{\prime }\in \Gamma $ and $\Gamma \stackrel{d}{\vdash }%
Aa^{\prime }a,Ab^{\prime }b$,

2.3. $\Gamma \stackrel{d}{\vdash }Aa^{\prime }a,Aa^{\prime }b$.

Assume 2.1 (the other two subcases are not harder), then there are an $%
a^{\prime }$-$a$ chain and a $b^{\prime }$-$b$ chain in $\Gamma $, let them
be, respectively $<c_{i}>_{i\in k}$ and $<c_{j}^{\prime }>_{j\in l}$. We may
assume that the ranges of these two chains are disjoint, otherwise this
subcase will be reduced to subcase 2.3. Also we may assume that each of
these two chains is injective. The following is a $d^{\prime }$ (hence a $%
d^{\prime \prime }$)-sorites of $\sigma $ from $\Gamma $: $Ab^{\prime
}c_{1}^{\prime },Ac_{1}^{\prime }c_{2}^{\prime },Ab^{\prime }c_{2}^{\prime
},...,Ab^{\prime }c_{l-2}^{\prime },Ac_{l-2}^{\prime }b,Ab^{\prime
}b,Ia^{\prime }b^{\prime },Ia^{\prime }b,Iba^{\prime },$ $Aa^{\prime
}c_{1}^{{}},Ibc_{1},...,$ $Ibc_{k-2},Ac_{k-2}a,Iba,Iab$.\newline

3. $\sigma =Eab$, for some $a,b\in C$ and $\Gamma $ is consistent. By
proposition 11.2 and lemma 8.3 there are $a^{\prime },b^{\prime }\in C$ such
that $\{Ea^{\prime }b^{\prime },Eb^{\prime }a^{\prime }\}\cap \Gamma \neq
\phi $ and there are an $a$-$a^{\prime }$ chain and a $b$-$b^{\prime }$
chain in $\Gamma $, let them be, respectively, $<c_{i}>_{i\in k}$ and $%
<c_{j}^{\prime }>_{j\in l}$. By the consistency of $\Gamma $, the ranges of
the two chains are disjoint. Moreover, we may assume that each of them is
injective. Let $Ea^{\prime }b^{\prime }\in \Gamma $ (the other case is not
harder), then the following is an $e$-sorites of $\sigma $ from $\Gamma $: $%
Ea^{\prime }b^{\prime },Ac_{k-2}a^{\prime },Ec_{k-2}b^{\prime
},...,Ec_{1}b^{\prime },Aac_{1}^{{}},Eab^{\prime },$ $Eb^{\prime
}a,Ac_{l-2}^{\prime }b^{\prime },Ec_{l-2}^{\prime }a,...,Ec_{1}^{\prime
}a,Abc_{1}^{\prime },Eba,Eab$.\newline

4. $\sigma =Oab$, for some $a,b\in C$. If $e=d$, then there is a one line $e$%
-sorites of $\sigma $ from $\Gamma $. Else assume that $\Gamma $ is
consistent and $e=d^{\prime \prime }$, by proposition 11.2 and lemma 8.3 it
suffices to deal with the following two subcases.

4.1. There are $a^{\prime },b^{\prime }\in C$ such that $Oa^{\prime
}b^{\prime }\in \Gamma $ and $\Gamma \stackrel{d}{\vdash }Aa^{\prime
}a,Abb^{\prime }$. Making use of Bocardo and Baroco it may be shown, along
the lines of part 3 of this proof, that there is a $d^{\prime }$ (hence a $%
d^{\prime \prime }$)-sorites of $\sigma $ from $\Gamma $.

4.2. There is $c\in C$ such that $\Gamma \stackrel{d^{\prime }}{\vdash }%
Ica,Ecb$. As in part 3 of this proof, there are $c^{\prime },b^{\prime }\in
C $ such that: 
\begin{equation}
\{Ec^{\prime }b^{\prime },Eb^{\prime }c^{\prime }\}\cap \Gamma \neq \phi 
\text{ \qquad and \qquad }\Gamma \stackrel{d}{\vdash }Acc^{\prime
},Abb^{\prime }  \tag{$*$}
\end{equation}
By lemma 8.3 it suffices to deal with the following two subsubcases.

4.2.1. For some $c^{\prime \prime }\in C$, $\Gamma \stackrel{d}{\vdash }%
Ac^{\prime \prime }c,Ac^{\prime \prime }a$. By this and ($*$), $\Gamma 
\stackrel{d}{\vdash }Abb^{\prime },Ac^{\prime \prime }c^{\prime },$ $%
Ac^{\prime \prime }a$. So there are $b$-$b^{\prime },c^{\prime \prime }$-$%
c^{\prime }$ and $c^{\prime \prime }$-$a$ injective $\Gamma $-chains; let
them be $<x_{i}>_{i\in k}, <y_{i}>_{i\in l}$ and $<z_{i}>_{i\in m}$
respectively.

By the consistency of $\Gamma $, the range of $<x_{i}>_{i\in k}$ and the
union of the ranges of $<y_{i}>_{i\in l}$ and $<z_{i}>_{i\in m}$ are
disjoint. Assume that the ranges of $<y_{i}>_{i\in l}$ and $<z_{i}>_{i\in m}$
have $c^{\prime \prime }$ only in common (the other case is similar).

If $Ec^{\prime }b^{\prime }\in \Gamma $, then there is a $d^{\prime }$
(hence a $d^{\prime \prime }$)-sorites of $\sigma $ from $\Gamma $. Else $%
Eb^{\prime }c^{\prime }\in \Gamma $ and the following is a $d^{\prime \prime }
$-sorites of $\sigma $ from $\Gamma $. $%
Abx_{1},Ax_{1}x_{2},Abx_{2},...,Abx_{k-2},$ $Ax_{k-2}b^{\prime },Abb^{\prime
},$ $Eb^{\prime }c^{\prime },Ebc^{\prime },Ec^{\prime }b,Ay_{l-2}c^{\prime
},Ey_{l-2}b,...,Ey_{1}b,Ac^{\prime \prime }y_{1},Ec^{\prime \prime
}b,Oc^{\prime \prime }b$ \newline
(here $E$-sub is made use of), $Ac^{\prime \prime
}z_{1},Oz_{1}b,Az_{1}z_{2},Oz_{2}b,...,Oz_{m-2}b,Az_{m-2}a,$ $Oab$.\newline

4.2.2. $\{Ic^{\prime \prime }a^{\prime },Ia^{\prime }c^{\prime \prime
}\}\cap \Gamma \neq \phi $ and $\Gamma \stackrel{d}{\vdash } Ac^{\prime \prime
}c,Aa^{\prime }a$, for some $c^{\prime \prime },a^{\prime }\in C$. By this
and (*), $\Gamma \stackrel{d}{\vdash }Abb^{\prime },Ac^{\prime \prime
}c^{\prime },Aa^{\prime }a$. So there are $b$-$b^{\prime },c^{\prime \prime
} $-$c^{\prime }$ and $a^{\prime }$-$a$ injective $\Gamma $-chains; let them
be $<x_{i}>_{i\in k},<y_{i}>_{i\in l}$ and $<z_{i}>_{i\in m}$ respectively.

By the consistency of $\Gamma $, the range of $<x_{i}>_{i\in k}$ and the
union of the ranges of $<y_{i}>_{i\in l}$ and $<z_{i}>_{i\in m}$ are
disjoint. If the ranges of $<y_{i}>_{i\in l}$ and $<z_{i}>_{i\in m}$ are not
disjoint, this case will be reduced to the above case; so assume that they
are disjoint.

If $Ia^{\prime }c^{\prime \prime }\in \Gamma $ then there is a $d^{\prime }$
(hence a $d^{\prime \prime }$)-sorites of $\sigma $ from $\Gamma $. Else $%
Ic^{\prime \prime }a^{\prime }\in \Gamma $, assume $Eb^{\prime }c^{\prime
}\in \Gamma $ (the other case is similar), then the following is a $d^{\prime
\prime }$-sorites of $\sigma $ from $\Gamma $. $%
Abx_{1},Ax_{1}x_{2},Abx_{2},...,Abb^{\prime },$ $Eb^{\prime }c^{\prime
},Ebc^{\prime },Ec^{\prime }b,Ay_{l-2}c^{\prime },$ $Ey_{l-2}b,...,Ec^{%
\prime \prime }b,Ic^{\prime \prime }a^{\prime },Oa^{\prime }b$ (here Ferison
is made use of), $Aa^{\prime }z_{1},Oz_{1}b,...,Oab$.$\square $\newline

\noindent PROPOSITION 11.5.\textbf{\ \quad }If $\Gamma (\subseteq BN(C))$ is
inconsistent then it is ds-inconsistent, in the sense that there is $\sigma
\in BN(C)$ such that $\Gamma \stackrel{ds}{\vdash }\sigma ,\widehat{\sigma }$%
.\newline

\noindent \textit{Proof.} Let $\Gamma $ be inconsistent, then there is a
universal $\rho \in BN(C)$ such that $\Gamma \stackrel{d}{\vdash }\rho ,%
\widehat{\rho }$. Distinguish between two cases:

1. $\rho =Aab$, for some $a,b\in C$. In this case the result is a direct
consequence of theorem 11.4.

2. $\rho =Eab$, for some $a,b\in C$. As in part 3 of the proof of theorem
11.4, there are $a^{\prime },b^{\prime }\in C$ such that $\{Ea^{\prime
}b^{\prime },Eb^{\prime }a^{\prime }\}\cap \Gamma \neq \phi $ and there are
injective $a$-$a^{\prime },b$-$b^{\prime }$ chains in $\Gamma $; let them
be, respectively, $<c_{i}>_{i\in k}$ and $<c_{j}^{\prime }>_{j\in l}$.

If the ranges of these chains are disjoint, the result follows by theorem
11.4 and the methods made use of in its proof. Else there is $c^{\prime
\prime }\in \{c_{i}:i\in k\}\cap \{c_{j}^{\prime }:j\in l\}$. Then there are
injective $c^{\prime \prime }$-$a^{\prime },c^{\prime \prime }$-$b^{\prime } 
$ chains in $\Gamma $. Along the lines of the proof of theorem 11.4 it may
be shown that $\Gamma \stackrel{ds}{\vdash }Ea^{\prime }c^{\prime \prime
},Ia^{\prime }c^{\prime \prime }$ (and $\Gamma \stackrel{ds}{\vdash }%
Eb^{\prime }c^{\prime \prime },Ib^{\prime }c^{\prime \prime }$). \hspace{10.8   cm} $\square $%
\newline

To show that the consistency condition in each of the parts 2,3 of theorem
11.4 cannot be completely dispensed with, we prove:\newline

\noindent PROPOSITION 11.6.\textbf{\ }Let $\Gamma \subseteq BN(C),e\in
\{d,d^{\prime },d^{\prime \prime }\}$ and $a,a^{\prime },b,b^{\prime
},c,c^{\prime }\in C$; and assume that $c\neq c^{\prime }$.\newline

1. If every $a$-$a^{\prime }$ chain in $\Gamma $ is a $<c,c^{\prime }>$
chain, then $Acc^{\prime }$ occurs as an assumption in every $e$-deduction
of $Aaa^{\prime }$ from $\Gamma $; moreover it is made use of as a premise
in the deduction if it is different from $Aaa^{\prime }$.\newline

In parts 2 and 3 below, $Ebb^{\prime }$ is assumed to be the only universal
negative sentence in $\Gamma $.
\newline

2. If every $a$-$b$ chain and every $a$-$b^{\prime }$ chain in $\Gamma $ is
a $<c,c^{\prime }>$ chain, then $Acc^{\prime }$ occurs as an assumption and
is made use of as a premise in every $e$-deduction of $Eaa^{\prime }$ and
every $e$-deduction of $Ea^{\prime }a$ from $\Gamma $.\newline

3. If every $a$-$b$ chain, every $a$-$b^{\prime }$ chain, every $a^{\prime }$%
-$b$ chain and every $a^{\prime }$-$b^{\prime }$ chain in $\Gamma $ is a $%
<c,c^{\prime }>$ chain, then for every injective $e$-deduction $<\sigma
_{i}>_{i\in k}$ of $Eaa^{\prime }$ from $\Gamma $ there is $j\in k$ such
that $\sigma _{j}$ is made use of as a premise at least twice.$\newline
$

In parts 4 and 5 below, $Obb^{\prime }$ is assumed to be the only negative
sentence in $\Gamma $.\newline

4. If every $b$-$a$ chain or every $a^{\prime }$-$b^{\prime }$ chain in $%
\Gamma $ is a $<c,c^{\prime }>$ chain, then $Acc^{\prime }$ occurs as an
assumption and is made use of as a premise in every $e$-deduction of $%
Oaa^{\prime }$ from $\Gamma $.\newline

5. If every $b$-$a$ chain and every $a^{\prime }$-$b^{\prime }$ chain in $%
\Gamma $ is a $<c,c^{\prime }>$ chain, then for every injective $e$%
-deduction $<\sigma _{i}>_{i\in k}$ of $Oaa^{\prime }$ from $\Gamma $ there
is $j\in k$ such that $\sigma _{j}$ is made use of as a premise at least
twice.\newline

\noindent \textit{Proof.} Generalize the first part to become:

For every $u,u^{\prime }\in C$, if every $u$-$u^{\prime }$ chain in $\Gamma $
is a $<c,c^{\prime }>$ chain, then $Acc^{\prime }$ occurs as an assumption in
every $e$-deduction of $Auu^{\prime }$ from $\Gamma $; moreover, it is made
use of as a premise in the deduction if it is different from $Auu^{\prime }$.

The stronger statement may be easily proved by course of values induction on
the length of the $e$-deduction.

Parts 2 and 4 may be proved similarly.\newline

Again generalize part 3 to become:

For every $u,u^{\prime }\in C$ if every $u$-$b$ chain, every $u$-$b^{\prime
} $ chain, every $u^{\prime }$-$b$ chain and every $u^{\prime }$-$b^{\prime
} $ chain in $\Gamma $ is a $<c,c^{\prime }>$ chain, then for every
injective $e$-deduction $<\sigma _{i}>_{i\in k}$ of $Euu^{\prime }$ from $%
\Gamma $ there is $j\in k$ such that $\sigma _{j}$ is made use of as a
premise at least twice.

The stronger statement may be proved by course of values induction on $k$ as
follows. Assume the required for $r<k$ and let $<\sigma _{i}>_{i\in k}$ be
an $e$-deduction of $Euu^{\prime }$ from $\Gamma $. Since $\{b,b^{\prime
}\},\{u,u^{\prime }\}$ are disjoint and $\sigma _{k-1}=Euu^{\prime }$, then
there are only two cases to consider:

1. For some $l<k-1,\sigma _{l}=Eu^{\prime }u$, in this case the result is
immediate by the induction hypothesis.

2. For some $l,m<k-1$ and some $v\in C$ it is the case that $l<m,$ $\{\sigma
_{l},\sigma _{m}\}=\{Auv,Evu^{\prime }\}$ and $\sigma _{k-1}$ is obtained
from them as the conclusion of applying the rule $\frac{Auv,Evu^{\prime }}{%
Euu^{\prime }}$.

If $\sigma _{l}$ is made use of as a premise in a step whose conclusion is $%
\sigma _{j}$ for some $j<k-1$, the result is immediate. Also if there is
some $u$-$v$ chain in $\Gamma $ which is not a $<c,c^{\prime }>$ chain, the
result follows by the induction hypothesis.

So it remains to assume that every $u$-$v$ chain in $\Gamma $ is a $%
<c,c^{\prime }>$ chain and for every $j<k-1$, $\sigma _{l}$ is not made use
of as a premise in the step which gives rise to $\sigma _{j}$. Put: 
\begin{eqnarray*}
\Delta _{o} &=&\{\sigma _{l}\}, \\
\Delta _{j+1} &=&\Delta _{j}\cup \{\sigma _{i}:i\in l\text{ and }\sigma _{i}%
\text{ is made use of as a premise} \\
&&\text{in a step whose conclusion is in }\Delta _{j}\}.
\end{eqnarray*}

Then for some $n$, $\Delta _{n+1}=\Delta _{n}$. Hence $<\sigma _{i}>_{i\in
l+1,\sigma _{i}\in \Delta _{n}}$ is an $e$-deduction of $\sigma _{l}$ from $%
\Gamma $, so by parts 1,2 above $Acc^{\prime }\in \Delta _{n}$.

Assume, towards a contradiction, that for every $i\in k,\sigma _{i}$ is made
use of as a premise at most once. Then $<\sigma _{i}>_{i\in m+1,\sigma
_{i}\notin \Delta _{n}}$ is an $e$-deduction of $\sigma _{m}$ from $\Gamma
-\Delta _{n}$. Again by parts 1,2 above, $Acc^{\prime }\notin \Delta _{n}$.
Hence the result.

Part 5 may be proved similarly. \hspace{7cm}
	 $\square $
\newline

EXAMPLES 11.7. To see that the consistency condition in each of the parts
2,3 of theorem 11.4 cannot be completely dispensed with, put:

$\Gamma _{o}=\{Aac,Aa^{\prime }c,Acc^{\prime },Ac^{\prime }b,Ac^{\prime
}b^{\prime },Ebb^{\prime }\}$ \quad , \quad $\sigma _{o}=Eaa^{\prime }$

$\Gamma _{1}=\{Aa^{\prime }c,Abc,Acc^{\prime },Ac^{\prime }a,Ac^{\prime
}b^{\prime },Obb^{\prime }\}$ \quad , \quad $\sigma _{1}=Oaa^{\prime }$.\newline

For $i\in 2$, $\Gamma _{i}\stackrel{d^{\prime }}{\vdash }\sigma _{i}$
(in fact $\Gamma _{o}\stackrel{d}{\vdash }\sigma _{o}$), but by the above
proposition $\Gamma _{i}\stackrel{d^{\prime \prime }s}{\nvdash }\sigma _{i}$.

This example still works even if $d^{\prime \prime }$ is augmented by all of
the Aristotelian syllogisms.

To see that consistency is not always necessary, just notice that whether $%
\Gamma $ is consistent or not, $\Gamma \stackrel{ds}{\vdash }\sigma $
whenever $\sigma \in \Gamma $.\newline

Following are basic properties of sorites.\newline

\noindent DEFINITIONS and remarks 11.8.\textbf{\ }Let $e\in \{d,d^{\prime
},d^{\prime \prime }\}$, $\Gamma \subseteq BN(C)$ and $k\in \Bbb{N}^{+}$,
and let $<\sigma _{i}>_{i\in k}$ be an $e$-sorites of $\sigma _{k-1}$ from $%
\Gamma $ according to some annotation.

1. Two annotations of an $e$-deduction from $\Gamma $ are said to be
essentially the same if the only difference between them is interchanging
``assumption'' (i.e. the corresponding sentence belongs to $\Gamma $) and
``A-Id'' in some places. An $e$-deduction from $\Gamma $ is said to have
essentially one, or unique, annotation (of some sort) if all of its
annotations (of this sort) are essentially the same.

2. For $1\leq \ell \leq k$, $<\sigma _{i}>_{i\in \ell }$ is an $e$-sorites
from $\Gamma $ according to the restriction of the given annotation iff $%
\ell =1$ or the annotation of $\sigma _{\ell -1}$ is neither ``A-Id'' nor
``assumption''.

3. Let $1\leq \ell <k$, then for at least one $j\in \{\ell ,\ell +1\}$, $%
<\sigma _{i}>_{i\in j}$ is an $e$-sorites from $\Gamma $ according to the
restriction of the given annotation.

4. For $e\in \{d,d^{\prime }\}$, every $e$-sorites from $\Gamma $ has
essentially one sorites annotation. This does not apply to $d^{\prime \prime
}$, for the $d^{\prime \prime }$-deduction $<Ixx,Exy,Oxy>$ has two $%
d^{\prime \prime }$-sorites annotations which are not essentially the same.
Only one of them is a $d^{\prime }$-sorites annotation.

5. Ferison and $E$-sub are the only $d^{\prime \prime }$-rules which are not 
$d^{\prime }$-rules. In every $d^{\prime \prime }$-sorites annotation at
most one of them is made use of, at most once.

6. In each $d^{\prime \prime }$-sorites at most one triple of the form $%
<Ixx,Exy,Oxy>$ or $<Exy,Ixx,Oxy>$ occurs, at most once.

If no such triple occurs, the sorities will have an essentially unique $%
d^{\prime \prime }$-sorites annotation. Else all of its $d^{\prime \prime }$%
-sorites annotations are essentially the same, with the only exception that
an occurrence of ``$Oxy$'' may be annotated as the consequence of the
preceding two sentences by Ferio (which is a $d^{\prime }$-rule) in some of
them and by Ferison (which is not) in the others.\newline

\textbf{12. Independence.}\newline

\noindent DEFINITIONS 12.1.\textbf{\ }Let $e$ be a deduction system and let $%
r$ be a rule of $e$. The deduction system obtained from $e$ by excluding $r$
will be denoted by ``$e_{r}$''.

1. $r$ is said to be derivable in $e$ if $\Gamma \stackrel{e_{r}}{\vdash }%
\sigma $ whenever $\Gamma $ is a set of antecedents of an instance of $r$,
and $\sigma $ is the corresponding conclusion. Otherwise $r$ is said to be
independent in $e$.

2. $e$ is said to be independent if each of its rules is independent in it.

3. $r$ is said to be weakly independent in $e$ if for some set $\Delta \cup
\{\rho \}$ of sentences, $\Delta \stackrel{e}{\vdash }\rho $ while $\Delta 
\stackrel{e_{r}}{\nvdash }\rho $.

4. $e$ is said to be weakly independent if each of its rules is weakly
independent in it.\newline

\noindent REMARKS 12.2.\textbf{\newline
}

1. Independence implies weak independence.

2. For $e\in \{d,d^{\prime },d^{\prime \prime }\}$, each rule $r$ of $e$ is
independent in $e$ iff it is weakly independent in $e$, hence $e$ is
independent iff it is weakly independent.

3. Each of $d,d^{\prime }$ is independent (cf. Glashoff (2005) where similar
results are obtained via brute force computation).

4. The independence of each of $ds$ and $d^{\prime }s$ is an immediate
consequence of the independence of each of $d$ and $d^{\prime }$
respectively.\newline

\noindent THEOREM 12.3.\textbf{\ }\newline

1. $E$-sub, Ferio and Ferison are derivable in $d^{\prime \prime }s$, hence
in $d^{\prime \prime }$. Each of the other rules of $d^{\prime \prime }$ is
independent in $d^{\prime \prime }$, hence in $d^{\prime \prime }s$.

2. $d^{\prime \prime }$ is not independent, hence not weakly independent.

3. $d^{\prime \prime }s$ is weakly independent; however, it is not
independent.\newline

\noindent \textit{Proof.}

1. The sequence $Aaa,Iaa,Eab,Oab$ shows that $E$-sub is derivable in $%
d^{\prime \prime }s$. The corresponding proofs for Ferio and Ferison are not
harder.\newline

Put $r$ = Bocardo and let $a,b,c$ be three pairwise distinct elements of $C$%
. Put $\Gamma =\{Aba,Obc\}$ and $\sigma =Oac$. It is easy to see that if $%
\Gamma \stackrel{d^{\prime \prime }}{\vdash }\rho $ then $\rho \in
\{Aba,Obc,Iab,Iba\}\cup \{Axx:x\in C\}\cup \{Ixx:x\in C\}$. From this the
independence of $r$ in $d^{\prime \prime }$ follows. Similarly the other
required results may be obtained.\newline

2. By part 1 above and part 2 of remarks 12.2.\newline

3. Part 1 above shows that $d^{\prime \prime }s$ is not independent. It
shows also that to prove the weak independence of $d^{\prime \prime }s$ it
suffices to deal with $E$-sub, Ferio and Ferison only.

Put $r=E$-sub and let $\Delta =\{Eab\}$ and $\rho =Oba$, for some distinct $%
a,b\in C$. $\Delta \stackrel{d^{\prime \prime }s}{\vdash }Oba$. To see that $%
\Delta \stackrel{d^{\prime \prime }s_{r}}{\nvdash }Oba$ notice that the only 
$d^{\prime \prime }s_{r}$ rules which yield an $O$-sentence are Ferio,
Baroco, Bocardo, and Ferison. To obtain $Oba$ by applying Baroco or Bocardo
the $O$-sentence occurring as one of the antecedents -in the present case-
will be the same as the conclusion, which is forbidden in sorites. To apply
Ferio or Ferison, the antencedents -in the present case- must be $Ibb$ and $%
Eba$. But if $\eta _{i},\eta _{i+1},\eta _{i+2}$ is a subsequence of a $%
d^{\prime \prime }s_{r}$ sorites deduction from $\Delta $ and the annotation
of $\eta _{i+2}$ is that it is obtained from $\eta _{i},\eta _{i+1}$ or $%
\eta _{i+1},\eta _{i}$ by some rule, then $\eta _{i+1}\in \Delta \cup
\{Acc:c\in C\}$. So neither Ferio nor Ferison is applicable, hence $\Delta 
\stackrel{d^{\prime \prime }s_{r}}{\nvdash }Oba$.

Next, put $r=$ Ferio (Ferison) and let $\Delta =\{Eab,Icb\}$ ($\{Eab,Ibc\}$)
and $\rho =Oca$ for some pairwise distinct $a,b,c\in C$. By a slight
modification of the above technique it may be shown that $\Delta \stackrel{%
d^{\prime \prime }s}{\vdash }\rho $, but $\Delta \stackrel{d^{\prime \prime
}s_{r}}{\nvdash }\rho $. \hspace{2.3cm} $\square $
\newline

\textbf{12.1. Independence of }$g$\textbf{\ and variations thereof. }$g$ and 
$d$ have the same deduction rules, but the notion of $g$-deduction is weaker
than that of $d$-deduction. By definition 1.9, for $\Gamma \cup \{\sigma
\}\subseteq BN(C)$, $\Gamma \stackrel{g}{\vdash }\sigma $ iff $\Gamma \cup \{%
\widehat{\sigma }\}$ is $d$-inconsistent. Likewise for each rule $r$ of $g$
(equivalently of $d$) define the $g_{r}$-deduction relation ``$\stackrel{%
g_{r}}{\vdash }$'' by: $\Gamma \stackrel{g_{r}}{\vdash }\sigma $ iff $\Gamma
\cup \{\widehat{\sigma }\}$ is $d_{r}$-inconsistent. From this and remarks
12.2 it easily follows that $r$ is independent in $g$ iff it is weakly
independent in $g$, hence $g$ is independent iff it is weakly independent.%
\newline

\noindent THEOREM 12.4.\textbf{\ }$g$ is independent.\newline

\noindent \textit{Proof.} Let $a,b,c$ be pairwise distinct elements of $C$,
and put $r=$ Barbara and $\Gamma =\{Aab,Abc\}$. $\Gamma \stackrel{g_{r}}{%
\nvdash }Aac$ iff $\Gamma \cup \{Oac\}$ is $d_{r}$-consistent. But the set
of all $d_{r}$-consequences of $\Gamma \cup \{Oac\}$ is $%
\{Aab,Abc,Oac,Iba,Icb\}\cup \{Axx:x\in C\}\cup \{Ixx:x\in C\}$, hence $%
\Gamma \cup \{Oac\}$ is $d_{r}$-consistent. Consequently Barbara is
independent in $g$.

The proofs of the independence of the other rules are similar or easier. \hspace{5pt}
$\square$
\newline

To get closer to the usual deduction systems, we introduce two new
deduction systems $g^{\prime }$, $g^{\prime \prime }$ and show that each of
them is equivalent to $g$ and discuss its independence.\newline

\textit{12.1.1. First variation on }$g$\textit{.} The deduction system $%
g^{\prime }$ is obtained by augmenting the system $d^{\prime }$ by the rule: 
$\frac{\rho ,\widehat{\rho }}{\sigma }$ (contradiction, Co for short).

Let $\Gamma \cup \{\sigma \}\subseteq BN(C)$. It is easy to see that if $%
\Gamma \stackrel{g^{\prime }}{\vdash }\sigma $ then there is a $g^{\prime }$%
-deduction of $\sigma $ from $\Gamma $ in which Co is never made use of or
it is made use of only at the last step; moreover, this applies to $%
g_{r}^{\prime }$ for each rule $r$ of $d^{\prime }$.\newline

\noindent THEOREM 12.5.\textbf{\ }The following are equivalent:

1. $\Gamma \stackrel{g}{\vdash }\sigma $,

2. $\Gamma \stackrel{g^{\prime }}{\vdash }\sigma $,

3. $\Gamma \stackrel{d^{\prime }}{\vdash }\sigma $ or $\Gamma $ is
inconsistent,

4. $\Gamma \cup \{\widehat{\sigma }\}$ is inconsistent.\newline

\noindent \textit{Proof.} Easy if $\Gamma $ is inconsistent; and in all
cases parts 1 and 4 are equivalent by definition 1.9.

Assume $\Gamma $ is consistent. By theorem 8.7, parts 1 and 3 are equivalent,
and by the definition of $g^{\prime }$, part 3 implies part 2. Finally
assume part 2, then there is a $g^{\prime }$-deduction of $\sigma $ from $%
\Gamma $ in which Co is never made use of, this implies part 3. \hspace{11.8cm} $\square $%
\newline

The following theorem settles the indepenence of $g^{\prime }$.\newline

\noindent THEOREM 12.6.\textbf{\ }\newline

1. Every rule of $g^{\prime }$ is independent in $g^{\prime }$ iff it is
weakly independent in $g^{\prime }$, hence $g^{\prime }$ is independent iff
it is weakly independent.

2. $\Gamma \stackrel{g_{co}}{\vdash }\sigma $ iff $\Gamma \stackrel{%
d^{\prime }}{\vdash }\sigma $.

3. For every rule $r$ of $d^{\prime }$, $\Gamma \stackrel{g_{r}^{\prime }}{%
\vdash }\sigma $ iff $\Gamma \stackrel{d_{r}^{\prime }}{\vdash }\sigma $ or $%
\Gamma $ is $d_{r}^{\prime }$-inconsistent.

4. $g^{\prime }$ is independent.\newline

\noindent \textit{Proof.} The proof of the first three parts is easy.

To prove the last part let $a,b\in C$, then $\frac{Aab,Oab}{Eab}$ is an
instance of Co. But by lemma 8.3, $\{Aab,Oab\}\stackrel{d^{\prime }}{\nvdash 
}Eab$. So by part 2 above $\{Aab,Oab\}\stackrel{g_{co}^{\prime }}{\nvdash }%
Eab$. Therefore Co is independent in $g^{\prime }$. To complete the proof
let $r$ be some other rule of $g^{\prime }$, then $r$ is a rule of $%
d^{\prime }$. By part 3 above the independence of $r$ in $g^{\prime }$ may
be proved by choosing a consistent set $\Gamma $ of antecedents of $r$ such
that $\Gamma \stackrel{d_{r}^{\prime }}{\nvdash }\sigma $, where $\sigma $
is the corresponding conclusion; which is always possible.\hspace{10pt} $\square $
\newline

\textit{12.1.2. Second variation on }$g$\textit{.} Though $g^{\prime }$ is
closer than $g$ to the contemporary deduction systems, it is not as close to
the Aristotelian spirit as $g$. Inspired by Gentzen-type sequent systems
(cf. Kleene, S.C. 1967, p. 306) we introduce a second variation $g^{\prime
\prime }$ on $g$, which will hopefully be close enough to both modern and
Aristotelian traditions. The deduction rules of $g^{\prime \prime }$ are:%
\newline

0$^{\prime }$. $\frac {}{\Gamma \vdash Aaa}$ $\qquad \qquad \quad ~($A-Id$%
^{\prime })$\newline

1$^{\prime }$. $\frac{\Gamma \vdash Aab}{\Gamma \vdash Iba}$ $\qquad \qquad
\quad ~(Apc^{\prime })$\newline

2$^{\prime }$. $\frac{\Gamma \vdash Eab}{\Gamma \vdash Eba}$ $\qquad \qquad
\quad ~(Ec^{\prime })$\newline

3$^{\prime }$. $\frac{\Gamma \vdash Aab,\Delta \vdash Abc}{\Gamma \cup
\Delta \vdash Aac}$ $\quad \qquad (Barbara^{\prime })$\newline

4$^{\prime }$. $\frac{\Gamma \vdash Aab,\Delta \vdash Ebc}{\Gamma \cup
\Delta \vdash Eac}$ $\quad \qquad (Celarent^{\prime })$\newline

5$^{\prime }$. $\frac {}{\Gamma \vdash \eta }$ $\quad \qquad \qquad \quad
~(Ass)$\newline

6$^{\prime }$. $\frac{\Gamma \cup \{\widehat{\sigma }\}\vdash \rho ,\Delta
\cup \{\widehat{\sigma }\}\vdash \widehat{\rho }}{\Gamma \cup \Delta \vdash
\sigma }$ $\quad ~(Raa)$\newline

where $a,b,c\in C$, $\Gamma \cup \Delta \cup \{\rho ,\sigma \}\subseteq
BN(C) $ and $\eta \in \Gamma $. ``Ass'' and ``Raa'' are abbreviations for
``Assumption'' and ``Reductio ad absurdum'' respectively. ``$\vdash $'' is
just a symbol, instead we could have made use of ordered pairs and write, 
e.g. ``$<\Gamma ,\sigma >$'' in place of ``$\Gamma \vdash \sigma $''.\newline

\noindent DEFINITION and remarks 12.7.\textbf{\ }Let $S$ be a set of
sequents, i.e. $S\subseteq \{$\noindent $\Gamma \vdash \sigma :\Gamma \cup
\{\sigma \mathbf{\}\subseteq }BN(C)\mathbf{\}}$.\newline

1. A $g^{\prime \prime }$-deduction from $S$ is a sequence $<\Gamma
_{i}\vdash \sigma _{i}>_{i\in k}$ of sequents, where $k\in \Bbb{N}$ and for
each $i\in k$, $\Gamma _{i}\vdash \sigma _{i}\in S$ or may be obtained from
preceding terms of the sequence by some $g^{\prime \prime }$-deduction rule.

If $k\neq 0$, $<\Gamma _{i}\vdash \sigma _{i}>_{i\in k}$ is said to be a $%
g^{\prime \prime }$-deduction of $\Gamma _{k-1}\vdash \sigma _{k-1}$ from $S 
$. In this case we write $S\stackrel{g^{\prime \prime }}{\Vdash }\Gamma
_{k-1}\vdash \sigma _{k-1}$.

2. We write ``$\Gamma \stackrel{g^{\prime \prime }}{\vdash }\sigma $'' for ``%
$\phi \stackrel{g^{\prime \prime }}{\Vdash }\Gamma \vdash \sigma $'', ``$%
g^{\prime \prime }$-deduction'' for ``$g^{\prime \prime }$-deduction from $%
\phi $'' and ``$g^{\prime \prime }$-deduction of $\Gamma _{k-1}\vdash \sigma
_{k-1}$ (or of $\Delta \vdash \rho $)'' for ``$g^{\prime \prime }$-deduction
of $\Gamma _{k-1}\vdash \sigma _{k-1}$ (or of $\Delta \vdash \rho $) from $%
\phi $''.\newline

3. The above definition and remark may be generalized to subsystems of $%
g^{\prime \prime }$.

4. The notions of derivability, independence and weak independence may be
extended to $g^{\prime \prime }$ in the obvious way.\newline

5. A deduction rule of $g^{\prime \prime }$ is independent in $g^{\prime
\prime }$ iff it is weakly independent in $g^{\prime \prime }$. Hence $%
g^{\prime \prime }$ is independent iff it is weakly independent.\newline

\noindent THEOREM 12.8.\textbf{\ }For every $\Gamma \cup \{\sigma
\}\subseteq BN(C)$: 
\[
\Gamma \stackrel{g}{\vdash }\sigma \text{ \qquad iff \qquad }\Gamma 
\stackrel{g^{\prime \prime }}{\vdash }\sigma \text{.}
\]

\noindent \textit{Proof.} Let $\Gamma \stackrel{g}{\vdash }\sigma $ then, by
definition 1.9, $\Gamma \cup \{\widehat{\sigma }\}\stackrel{d}{\vdash }\rho ,%
\widehat{\rho }$  \, for some $\rho \in BN(C)$. Let $<\xi _{i}>_{i\in k}$ and $%
<\eta _{j}>_{j\in l}$ be, respectively, $d$-deductions of $\rho $ and $%
\widehat{\rho }$ from $\Gamma \cup \{\widehat{\sigma }\}$. Then $<\Gamma
\cup \{\widehat{\sigma }\}\vdash \xi _{i}>_{i\in k}$ and $<\Gamma \cup \{%
\widehat{\sigma }\}\vdash \eta _{j}>_{j\in l}$ are, respectively, $g^{\prime
\prime }$-deductions of $\Gamma \cup \{\widehat{\sigma }\}\vdash \rho $ and $%
\Gamma \cup \{\widehat{\sigma }\}\vdash \widehat{\rho }$. To their
concatenation (which is a $g^{\prime \prime }$-deduction) add one more line
to obtain $\Gamma \vdash \sigma $ by Raa from lines $k-1$ and $k+l-1$. This
proves the only if direction.

To prove the other direction let $\Gamma \stackrel{g^{\prime \prime }}{%
\vdash }\sigma $, let $<\Delta _{i}\vdash \rho _{i}>_{i\in k}$ be a $%
g^{\prime \prime }$-deduction of $\Gamma \vdash \sigma $, and assume that $%
\Delta _{i}\stackrel{g}{\vdash }\rho _{i}$ for each $i<k-1$. To show that $%
\Delta _{k-1}\stackrel{g}{\vdash }\rho _{k-1}$ we deal with as many cases as
there are $g^{\prime \prime }$-deduction rules. Following we consider
Celarent$^{\prime }$ (4$^{\prime }$) and Raa (6$^{\prime }$), the other
cases are similar or easier.

Celarent$^{\prime }$: There are $a,b,c\in C$ and $j,l\in \Bbb{N}$ such that $%
j<l<k-1,$ $\Delta _{k-1}=\Delta _{j}\cup \Delta _{l},\{\rho _{j},\rho
_{l}\}=\{Aab,Ebc\}$ and $\rho _{k-1}=Eac$. By the above assumption $\Delta
_{j}\stackrel{g}{\vdash }\rho _{j}$ and $\Delta _{l}\stackrel{g}{\vdash }%
\rho _{l}$. Hence $\Delta _{k-1}(=\Delta _{j}\cup \Delta _{l})\stackrel{g}{%
\vdash }\rho _{j},\rho _{l}$. So by part 7 of proposition 1.11, $\Delta
_{k-1}\stackrel{g}{\vdash }\rho _{k-1}$.\newline

Raa: There are $\Sigma ,\Sigma ^{\prime },\eta ,j,l$ such that $\Sigma \cup
\Sigma ^{\prime }\cup \{\eta \}\subseteq BN(C)$; $j,l\in \Bbb{N}$; $j<l<k-1$%
, $\Delta _{j}=\Sigma \cup \{\widehat{\rho }_{k-1}\}$, $\Delta _{l}=\Sigma
^{\prime }\cup \{\widehat{\rho }_{k-1}\}$, $\{\rho _{j},\rho _{l}\}=\{\eta ,%
\widehat{\eta }\}$ and $\Delta _{k-1}=\Sigma \cup \Sigma ^{\prime }$. By the
above assumption $\Delta _{j}\stackrel{g}{\vdash }\rho _{j}$ and $\Delta _{l}%
\stackrel{g}{\vdash }\rho _{l}$. Hence $\Delta _{k-1}\cup \{\widehat{\rho }%
_{k-1}\}(=\Delta _{j}\cup \Delta _{l})\stackrel{g}{\vdash }\rho _{j},\rho
_{l}$. So by part 4 of proposition 1.11 and the relevant definitions, $%
\Delta _{k-1}\stackrel{g}{\vdash }\rho _{k-1}$.\hspace{7.7cm} $\square $
\newline

Notice that in the $g^{\prime \prime }$-deductions $<\Gamma \cup \{\widehat{%
\sigma }\}\vdash \xi _{i}>_{i\in k}$ and $<\Gamma \cup \{\widehat{\sigma }%
\}\vdash \eta _{j}>_{j\in l}$ which occur in the proof of the only if
direction of the above theorem, only the rules $0^{\prime }-5^{\prime }$ are
made use of. Moreover, if $\Gamma \stackrel{d}{\vdash }\sigma $, then there
is a $g^{\prime \prime }$-deduction of $\Gamma \vdash \sigma $ in which Raa
is never made use of.

This is essentially sufficient to prove the following:\newline

\noindent COROLLARY 12.9. If $\Gamma \stackrel{g^{\prime \prime }}{\vdash }%
\sigma $ then there is a $g^{\prime \prime }$-deduction of $\Gamma \vdash
\sigma $ in which Raa is never made use of or is made use of only in the
last step. \hspace{2.0cm} $\square $
\newline

This section is concluded by proving the independence of $g^{\prime \prime }$%
.\newline

\noindent THEOREM 12.10.\textbf{\ }Let $\Gamma \cup \{\sigma \}\subseteq
BN(C)$, $r$ be a $d$-deduction rule, and $r^{\prime }$ be the corresponding $%
g^{\prime \prime }$-deduction rule.

1. $\Gamma \stackrel{g_{r^{\prime }}^{\prime \prime }}{\vdash }\sigma $
\quad ~~iff \quad $\Gamma \stackrel{g_{r}}{\vdash }\sigma $,

2. $\Gamma \stackrel{g_{Ass}^{\prime \prime }}{\vdash }\sigma $ \quad iff
\quad $\phi \stackrel{g^{\prime \prime }}{\vdash }\sigma $ \quad iff \quad $%
\phi \stackrel{g}{\vdash }\sigma $ \quad iff \quad $\phi \stackrel{d}{\vdash 
}\sigma $,

3. $\Gamma \stackrel{g_{Raa}^{\prime \prime }}{\vdash }\sigma $ \quad iff
\quad $\Gamma \stackrel{d}{\vdash }\sigma $,

4. $g^{\prime \prime }$ is independent.\newline

\noindent \textit{Proof.} The proofs of the first and the third parts are
along the lines of the proof of theorem 12.8 noting that proposition 1.11
still holds after replacing ``$d$''$,$``$g$'' by ``$d_{r}$''$,$``$g_{r}$''
respectively. For part 2 it is sufficient to notice that each of the four
statements holds iff $\sigma $ is of the form $Ycc$ for some $c\in C$ and
some $Y\in \{A,I\}$.\newline

To prove the last part we consider three cases:

Rules $0^{\prime }-4^{\prime }$: Let $r^{\prime }$ be one of these rules and
let $r$ be the corresponding $d$-rule. Since $r$ is independent in $g$, there
is a set $\Gamma $ of antecedents of an instance of $r$ such that $\Gamma 
\stackrel{g_{r}}{\nvdash }\sigma $, where $\sigma $ is the corresponding
conclusion. Put $S=\{\Gamma \vdash \rho :\rho \in \Gamma \}$ then $S$ is a
set of antecedents of $r^{\prime }$ and $\Gamma \vdash \sigma $ is the
corresponding conclusion. By parts 2,3 of definition and remarks 12.7 and
part 1 above:

$S\stackrel{g_{r^{\prime }}^{\prime \prime }}{\Vdash }\Gamma \vdash \sigma $
\quad iff \quad $\phi \stackrel{g_{r^{\prime }}^{\prime \prime }}{\Vdash }%
\Gamma \vdash \sigma $ \quad iff \quad $\Gamma \stackrel{g_{r^{\prime
}}^{\prime \prime }}{\vdash }\sigma $ \quad iff \quad $\Gamma \stackrel{g_{r}%
}{\vdash }\sigma $.

But $\Gamma \stackrel{g_{r}}{\nvdash }\sigma $, hence $r^{\prime }$ is
independent in $g^{\prime \prime }$.

Rule Ass: For $a,b\in C$, $\phi \stackrel{g_{{}}^{\prime \prime }}{\Vdash }%
\{Oab\}\vdash Oab$ while, by 2 above, $\phi \stackrel{g_{Ass}^{\prime \prime
}}{\nVdash }\{Oab\}\vdash Oab$. Hence Ass is independent in $g^{\prime
\prime }$.\newline

Rule Raa: For distinct elements $a,b$ of $C$ let $\sigma =Iab$, $\rho =Iba$, 
$\Gamma =\{\rho \}$ and $S=\{\Gamma \cup \{\widehat{\sigma }\}\vdash \rho
,\Gamma \cup \{\widehat{\sigma }\}\vdash \widehat{\rho }\}$ then S is a set
of antecedents of Raa and $\Gamma \vdash \sigma $ is the corresponding
conclusion. By a slight modification of the proof given above for the rules $%
0^{\prime }-4^{\prime }$ it may be shown that Raa is independent in $%
g^{\prime \prime }$. \hspace{12.5cm} $\square $
\newline

\textbf{13. Algebraic semantics of AAS, a prelude. }The most well known
attempt to algebraically interpret Aristotelian syllogistic is that of Boole
(1948, first published 1847); however, it is not the first. More than a
century and a half earlier, this area of research was pioneered by Leibniz
(Kneale and Kneale 1966, pp. 338-45; Lenzen 2004). Following is a discussion
of the subject in general; the works of Leibniz and Boole will be briefly
discussed in section 17 below.

Regarding the central role played by order models in the semantics of $NF(C)$%
, they will be our starting point for algebraization. Each underlying order
structure of an order model will induce an algebra which may be expanded to
make the interpretation of $NF(C)$ possible.

The simplicity of order models stems from the fact that all relations are
determined by only one of them, namely the interpretation of $A$, which is
compatible with the Aristotelian view that Barbara is the essential
syllogism. Likewise, algebras defined in this section will each have one
(partial) binary operation and no others.\newline

\noindent DEFINITIONS and remarks 13.1.\newline

1. Let $B$ be a non-empty set and let $\oplus $ be a function from a subset
of $B\times B$ to $B$, then $\oplus $ is said to be a partial binary
operation on $B$, and $<B,\oplus >$ is said to be a partial algebra.\newline

2. Let $\stackrel{\circ }{\leq }$ be a binary relation on a set $B$, the
partial binary operation $+_{\stackrel{\circ }{\leq }}$ induced by $%
\stackrel{\circ }{\leq }$ on $B$ is defined by: 
\[
+_{\stackrel{\circ }{\leq }}:\text{ }\stackrel{\circ }{\leq }\text{ }%
\rightarrow \text{ }B 
\]
\[
\text{ \quad \qquad }a+_{\stackrel{\circ }{\leq }}b=a 
\]
$+_{\stackrel{\circ }{\leq }}$ is commutative (see 3.3 below) if $\stackrel{%
\circ }{\leq }$ is antisymmetric. $<B,+_{\stackrel{\circ }{\leq }}>$ is
called the partial algebra induced by $<B,\stackrel{\circ }{\leq }>$.\newline

3. Let $<B,\stackrel{\circ }{\leq }>$ be an order structure, then $<B,+_{%
\stackrel{\circ }{\leq }}>$ satisfies:

\quad 1. Right associativity: 
\[
(a+_{\stackrel{\circ }{\leq }}b)+_{\stackrel{\circ }{\leq }}c=a+_{\stackrel{%
\circ }{\leq }}(b+_{\stackrel{\circ }{\leq }}c) 
\]
in the sense that for every $a,b,c\in B$ if the rhs exists, so does the lhs
and they are equal.

\quad 2. Idempotence: 
\[
a+_{\stackrel{\circ }{\leq }}a=a\text{ \qquad all }a\in B. 
\]
If, moreover, $\stackrel{\circ }{\leq }$ is antisymmetric, then $<B,+_{%
\stackrel{\circ }{\leq }}>$ satisfies:

\quad 3. Commutativity: 
\[
a+_{\stackrel{\circ }{\leq }}b=b+_{\stackrel{\circ }{\leq }}a 
\]
in the sense that for every $a,b\in B$, if both sides exist they are equal.

Honouring Leibniz, a partial algebra $<B,\oplus >$ satisfying conditions 1
and 2 will be called a Leibniz algebra (LA for short). If, moreover, it
satisfies condition 3 it will be called a commutative Leibniz algebra (CLA
for short).

So, $<B,+_{\stackrel{\circ }{\leq }}>$ is a LA if $<B,\stackrel{\circ }{\leq 
}>$ is an order structure; moreover, it is a CLA if $\stackrel{\circ }{\leq }
$ is antisymmetric.\newline

4. An idempotent partial algebra $<B,\oplus >$ will be called a weak Leibniz
algebra (WLA for short) if it satisfies:

\quad 1. Weak right associativity: 
\[
(a\oplus b)\oplus c=a\oplus (b\oplus c) 
\]
in the sense that for every $a,b,c\in B$ if $(a\oplus b)$ and the rhs both
exist, then the lhs exists and equals the rhs.

If, moreover, $<B,\oplus >$ is commutative (in the sense of condition 3.3
above), it will be called a commutative weak Leibniz algebra (CWLA for
short).

Obviously every LA (CLA) is a WLA (CWLA).\newline

5. With abuse of notation, ``LA'', ``CLA'', ``WLA'' and ``CWLA'' will denote
also the classes of all LAs, CLAs, WLAs and CWLAs respectively; what is
intended will be clear from the context.

Abuses of notations such as this may take place later on without further
notice.\newline

6. Let $<B,\oplus >$ be a partial algbra, the binary relation $\leq _{\oplus
}$ induced by $\oplus $ on $B$ is defined by: 
\[
\leq _{\oplus }=\{<a,b>\in B:<<a,b>,a>\in \oplus \}, 
\]
so $a\leq _{\oplus }b$ iff $a\oplus b=a$ (in the sense that the lhs exists
and equals the rhs). Obviously, $\leq _{\oplus }$ is antisymmetric if $%
\oplus $ is commutative.\newline

7. Let $\stackrel{\circ }{\leq }$ and $\oplus $ be, respectively, a binary
relation and a partial binary operation on a set $B$, and let $+_{\stackrel{%
\circ }{\leq }}$, $\leq _{\oplus }$, $\leq _{+_{\stackrel{\circ }{\leq }}}$
and $+_{\leq _{\oplus }}$ be as defined above. Then:

\quad 1. $\leq _{+_{\stackrel{\circ }{\leq }}}=$ $\stackrel{\circ }{\leq }$.

\quad 2. $+_{\leq _{\oplus }}\subseteq \oplus $; moreover, $+_{\leq _{\oplus
}}$ is commutative if $\oplus $ is.\newline

8. Let $<B,\oplus >$ be a WLA, then $<B,\leq _{\oplus }>$ is an order
structure, called the order structure induced by $<B,\oplus >$. Moreover, $%
\leq _{\oplus }$ is antisymmetric if $\oplus $ is commutative.\newline

\newpage

\textbf{14. Algebraic interpretation of }$NF(C)$\textbf{.}\newline

\noindent DEFINITION 14.1. Let $\mathfrak{B}=<B,\oplus >$ be a WLA and let $\mu
:C\rightarrow B$. The structure $\mathfrak{B}^{\mu }=<B,\oplus ,\mu >$ is said
to be a weak Leibniz structure (WLS for short). The reduct $\mathfrak{B}$ shall
be called the WLA base of $\mathfrak{B}^{\mu }$.

Leibniz structures (LS for short), commutative Leibniz structures (CLS for
short) and commutative weak Leibniz structures (CWLS for short) are defined
analogously.

The following definition shows how $NF(C)$ may be interpreted in these
structures. So they may, and will, be considered as $NF(C)$-structures and
will be treated like other $NF(C)$-structures when dealing with semantics. In
particular, all semantical notions (such as ``$\mathfrak{B}^{\mu }$ is a (n
algebraic) model of $\Gamma $'' or ``$\Gamma \underset{\mathfrak{C}}{\vDash }%
\sigma $'', for $\Gamma \cup \{\sigma \}\subseteq BN(C)$ and $\mathfrak{C}
\subseteq $WLS) will be assumed to be known.\newline

\noindent DEFINITION 14.2. Let $\mathfrak{B}^{\mu }=<B,\oplus ,\mu >$ be a WLS,
and let $a,b\in C$, then:\newline
\newline
\begin{tabular}{ll}
1. $\mathfrak{B}^{\mu }\vDash Aab$ & iff $\mu a\oplus \mu b$ exists and equals $%
\mu a$ \\ 
& (iff $<<\mu a,\mu b>,\mu a>\in \oplus ).$ \\ 
&  \\ 
2. $\mathfrak{B}^{\mu }\vDash Iab$ & iff the system of equations $x\oplus \mu
a=x $, \\ 
& $x\oplus \mu b=x$ has a solution \\ 
& (iff the equation $x\oplus \mu a=x\oplus \mu b$ has a solution, \\ 
& iff the equation $x\oplus \mu a=y\oplus \mu b$ has a solution). \\ 
&  \\ 
3. $\mathfrak{B}^{\mu }\vDash Eab$ & iff $\mathfrak{B}^{\mu }\nvDash Iab.$ \\ 
&  \\ 
4. $\mathfrak{B}^{\mu }\vDash Oab$ & iff $\mathfrak{B}^{\mu }\nvDash Aab.$%
\end{tabular}
\newline
\newline

\noindent REMARKS 14.3.

\begin{description}
\item  1. The order (partial order) model $<B,\stackrel{\circ }{\leq },\mu >$
is basically equivalent to the WLS (CLS), $<B,+_{\stackrel{\circ }{\leq }%
},\mu >$.

\item  2. The WLS (CLS), $<B,\oplus ,\mu >$, is basically equivalent to the
order (partial order) model $<B,\leq _{\oplus },\mu >$.

\item  3. Consequently, every WLS (hence every LS, every CWLS and every CLS)
is an $e$-model for $e\in \{d,d^{\prime },d^{\prime \prime },g\}$.

\item  4. In the light of remarks and definitions 9.10 it may be easily seen
that:

\begin{description}
\item  1. For $e\in \{d,d^{\prime },d^{\prime \prime },g\}$, $e$ is sound
wrt WLS, hence wrt every subclass of it.

\item  2. $g$ is CLS-syntactically complete.

\item  3. For $e\in \{d^{\prime },d^{\prime \prime },g\}$, $e$ is
CLS-consistently syntactically complete.

\item  4. For $e\in \{d,d^{\prime },d^{\prime \prime },g\}$, $e$ is
CLS-consistently semantically complete.
\end{description}
\end{description}

In clauses 2-4, CLS may be replaced by any class intermediate between it
and WLS.\newline

\textbf{15. Annihilators: Embedding the partial into a total. }An
annihilator of a (partial) binary operation $*$ on a set $B$ is an element $%
b\in B$ such that: 
\[
x*b=b=b*x\text{ \qquad all }x\in B 
\]
Obviously $*$ has at most one annihilator.

An annihilator algebra is an ordered triple $\mathfrak{B}=<B,*,b>$ such that the
reduct $^{r}\mathfrak{B}=<B,*>$ is a partial algebra, and $b$ is an annihilator
of $*$.

The subreduct $^{sr}\mathfrak{B}$ of $\mathfrak{B}$ is the ordered pair $<B^{\prime
},*^{\prime }>$, where: 
\[
B^{\prime }=B-\{b\}\text{ \qquad , \qquad }*^{\prime }=*\cap (B^{\prime
}\times B^{\prime })\times B^{\prime } 
\]
Here, and in the sequel, $B^{\prime }$ is assumed to be non-empty.\newline

\noindent DEFINITIONS 15.1.\newline

1. An annihilator Leibniz algebra (ALA for short) is an annihilator algebra
whose subreduct is a LA.

Annihilator commutative Leibniz algebras (ACLA for short), annihilator weak
Leibniz algebras (AWLA for short) and annihilator commutative weak Leibniz
algebras (ACWLA for short) are defined analogously.\newline

2. A Leibniz algebra with annihilator (LAA for short) is an annihilator
algebra whose reduct is a LA.

Commutative Leibniz algebras with annihilators (CLAA for short), weak
Leibniz algebras with annihilators (WLAA for short) and commutative weak
Leibniz algebras with annihilators (CWLAA for short) are defined analogously.

\vspace{10pt}
As usual, an algebra or a structure based on an algebra is said to be total
if each of its operations is total. ``TLA'' will stand for ``total Leibniz
algebra'', ``TLS'' will stand for ``total Leibniz structure'' and similarly
for the other cases.\newline

\noindent REMARKS 15.2.\newline

1. The subreduct of an annihilator algebra is a LA (respectively CLA, WLA or
CWLA) if the reduct is.

Hence LAA $\subseteq $ ALA, and similarly for the other cases.\newline

2. Let $\mathfrak{B}=<B,*,b>$ be a total annihilator algebra whose subreduct
also is total. Then $\mathfrak{B}$ is LAA iff it is ALA; ``LA'' may be replaced
by ``CLA'', ``WLA'' or ``CWLA''.\newline

\begin{tabular}[t]{ll}
3. TCLAA = TCWLAA & = ICSGA (idempotent commutative \\ 
& \qquad \qquad ~~~semigroups with annihilators) \\ 
& = OSLA (operational semilattices \\ 
& \qquad \quad \quad \ with annihilators)
\end{tabular}

The order structures induced by these algebras are lower semilattices with
smallest elements.

In the above equations ``C'', the last ``A'' or both, may be dropped every-
where (the corresponding parenthetic clause is to be modified accordingly).%
\newline

The following definition designates to each partial algebra a total
annihilator algebra in which it may be embedded.\newline

\noindent DEFINITION and remarks 15.3.\newline

1. For $i\in 2$, let $\mathfrak{B}_{i}(=<B_{i},*_{i}>)$ be a partial algebra. A
bijection $f$ from $B_{0}$ to $B_{1}$ is said to be an isomorphism from $%
\mathfrak{B}_{0}$ to $\mathfrak{B}_{1}$ if \noindent for every $x,y,z\in B_{0}$: 
\[
<<x,y>,z>\in *_{0}\text{ \quad iff \quad }<<fx,fy>,fz>\in *_{1}\text{.} 
\]

$\mathfrak{B}_{0}$ and $\mathfrak{B}_{1}$ are said to be isomorphic if there is an
isomorphism from one of them to the other.\newline

2. If two partial algebras are isomorphic and one of the partial binary
operations has an annihilator, then its image is an annihilator of the other.%
\newline

3. Two annihilator algebras are said to be isomorphic if their reducts are.%
\newline

4. Two total annihilator algebras are isomorphic iff their subreducts are.%
\newline

5. Every partial algebra $\mathfrak{B}(=<B,*>)$ is the subreduct of some total
annihilator algebra. For, let $0 \notin B$ and $_{0}B=B\cup \{0 \}$.

Put: 
\begin{eqnarray*}
_{0}*\text{ } &:&\text{ }_{0}B\times \text{ }_{0}B\rightarrow \text{ }_{0}B
\\
x\text{ }_{0}*\text{ }y &=&
\begin{cases}
x*y & \text{if} <x,y>\in  \text{Domain} \, *, \\ 
$0$ & \text{otherwise}
\end{cases}
\end{eqnarray*}

Then $_{0}\mathfrak{B}(=<$ $_{0}B,$ $_{0}*,0>)$ is a total annihilator algebra,
and $\mathfrak{B}$ is its subreduct.\newline
6. Every total annihilator algebra whose subreduct is isomorphic to $\mathfrak{B}
$, is isomorphic to $_{0}\mathfrak{B}$. This warrants calling $_{0}\mathfrak{B}$ the
total annihilator algebra induced by $\mathfrak{B}$.

The identity map on $B$ is an embedding of $\mathfrak{B}$ into $_{0}\mathfrak{B}$.%
\newline

7. $_{0}\mathfrak{B}$ is a TALA iff $\mathfrak{B}$ is a LA. ``LA'' may be replaced
by ``CLA'', ``WLA'' or ``CWLA''.\newline

\textbf{16. Back to algebraic interpretation. }Let $\mathfrak{B}(=<B,*,0>)$ be a
WLAA, then its reduct $^{r}\mathfrak{B}(=<B,*>)$ is a WLA. So $<B,*,\mu >$ is a
WLS, for every ${\mu }:C\rightarrow B$. Obviously, for all $a,b\in C$, $%
Iab$ is satisfied in this structure. Hence none of $d,d^{\prime },d^{\prime
\prime }$ nor $g$ is consistently semantically complete wrt any class of
such structures, though every one of them is sound wrt each of these
classes. Evidently expanding the structure to $<B,*,0,\mu >$ will not solve
the problem.

As a matter of fact, the annihilator is the source of the difficulty, and we
may get around it by not permitting the annihilator to be assigned as a
value corresponding to any element of $C$, nor accepting it as a solution of
any of the relevant equations below. An additional advantage of this
approach is to be able to consider the more general AWLA.\newline

\noindent DEFINITION 16.1. (non-annihilator interpretation of $NF(C)$ in
annihilator algebras)\newline

1. Let $\mathfrak{B}(=<B,*,0>)$ be an AWLA and let $\mu :C\rightarrow B^{\prime
}(=B-\{0\})$. The structure $\mathfrak{B}^{\mu }(=<B,*,0,\mu >)$ is called an
annihilator weak Leibniz structure (AWLS for short). The reduct $\mathfrak{B}$
of $\mathfrak{B}^{\mu }$ is called the AWLA base of $\mathfrak{B}^{\mu }$.

The structures based on the other algebras (total or not) are defined, and
their names are abbreviated, analogously.\newline

2. For each $a,b\in C$:

\begin{tabular}{ll}
1. $\mathfrak{B}^{\mu }\vDash Aab$ & iff $\mu a*\mu b$ exists and equals $\mu a$
\\ 
& (equivalently $<<\mu a,\mu b>,\mu a>\in *$). \\ 
2. $\mathfrak{B}^{\mu }\vDash Iab$ & iff the system of equations $x*\mu a=x$, \\ 
& $x*\mu b=x$ has a solution different from $0$ \\ 
& (iff the equation $x*\mu a=x*\mu b$ \\ 
& has a solution which makes $x*\mu b\neq 0$, \\ 
& iff the equation $x*\mu a=y*\mu b$ \\ 
& has a solution which makes $y*\mu b\neq 0$). \\ 
3. $\mathfrak{B}^{\mu }\vDash Eab$ & iff $\mathfrak{B}^{\mu }\nvDash Iab$. \\ 
4. $\mathfrak{B}^{\mu }\vDash Oab$ & iff $\mathfrak{B}^{\mu }\nvDash Aab$.
\end{tabular}

This shows how $NF(C)$ may be interpreted in the structures defined in part
1. So they may, and will, be considered as $NF(C)$-structures and will be
treated like other $NF(C)$-structures when dealing with semantics. In
particular, all semantical notions (such as ``$\mathfrak{B}^{\mu }$ is a (total
algebraic) model of $\Gamma $'' or ``$\Gamma \vDash _{\mathfrak{C}}\sigma $'',
for $\Gamma \cup \{\sigma \}\subseteq BN(C)$ and $\mathfrak{C}\subseteq $AWLS)
will be assumed to be known.\newline

\noindent REMARKS 16.2.\newline

1. $\mathfrak{B}^{\mu }$ and $^{sr}\mathfrak{B}^{\mu }$ are basically equivalent,
hence every AWLS is an $e$-model for $e\in \{d,d^{\prime },d^{\prime \prime
},g\}$.\newline

2. In part 4 of remarks 14.3, ``WLS'' and ``CLS'' may be, respectively,
replaced by ``TAWLS'' and ``TACLS''.\newline

3. If $\mathfrak{B}^{\mu }$ is a TCLSA (equivalently TCWLSA), the provisions
given in part 2 of definitions 16.1 may be simplified in the obvious way; in
particular, the second provision will be equivalent to $\mu a*\mu b\neq 0$.%
\newline

To investigate the relationship between TCLSA and the Venn models we make
use of a (n intermediate) subclass of TCLSA, namely the subclass of those
TCLSA based on OSLA which are reducts of Boolean algebras.

These reducts will be called Boolean-Leibniz algebras with annihilators,
BLAA for short. As usual BLSA is a Boolean-Leibniz structure with
annihilator, i.e. a LSA based on a BLAA.\newline

\noindent PROPOSITION 16.3.\textbf{\ }Every TCLAA may be embedded in a BLAA.%
\newline

\noindent \textit{Proof.} Let $\mathfrak{B}(=<B,*,0>)$ be a TCLAA. The mapping:

\begin{eqnarray*}
f &:&B\rightarrow \wp (B) \\
f(b) &=&\{x\in B:x*b=x\}-\{0\}
\end{eqnarray*}
is an embedding of $\mathfrak{B}$ in the BLAA: $<\wp (B),\cap ,\phi >$. \hspace{4cm} $\square $%
\newline

$<\wp (B),\cap ,\phi >$ will be called the BLAA corresponding to $\mathfrak{B}$
and will be denoted by ``$Bl(\mathfrak{B})$''. For $\mu :C\rightarrow B$, $<\wp
(B),\cap ,\phi ,f\mu >$ is a BLSA; it will be called the BLSA corresponding
to, the TCLSA, $\mathfrak{B}^{\mu }$ and will be denoted by ``$Bl(\mathfrak{B}^{\mu
})$''. The relevant definitions and part 3 of remarks 16.2 show that $\mathfrak{B}^{\mu }$ and $Bl(\mathfrak{B}^{\mu })$ are basically equivalent.\newline

\noindent THEOREM 16.4.\textbf{\ }Every TCLSA is basically equivalent to a
Venn model. And every Venn model is basically equivalent to a BLSA (hence to
a TCLSA); moreover, the BLSA may be assumed to be based on a concrete BLAA
whose universe is a power set.\newline

\noindent \textit{Proof.} Let $\mathfrak{B}^{\mu }(=<B,*,0,\mu >)$ be a TCLSA,
then $Bl(\mathfrak{B}^{\mu })$ is a BLSA and $\mathfrak{B}^{\prime }=<\wp (B)-\{\phi
\},f\mu >$ is a Venn model. They all are basically equivalent.

On the other hand, let $\mathfrak{B}(=<B,\mu >)$ be a Venn model, then $<\wp
(\bigcup B),\cap ,$ $\phi ,\mu >$ is a BLSA which is basically equivalent to
it. \hspace{3.8cm} $\square $
\newline
\newline
\noindent COROLLARY 16.5.\textbf{\ }In part 4 of remarks 14.3, ``WLS'' and
``CLS'' may be, respectively, replaced by ``TWLSA'' and ``BLSA'' (either the
superclass TCLSA, or the subclass consisting of those elements each of which
is based on a concrete BLAA whose universe is a power set, may replace BLSA). 
\hspace{1.5cm} $\square $
\newline

\textbf{17. Leibniz and Boole. }The calculus de continentibus et contentis,
or the calculus of identity and inclusion -which is an algebraic treatment
of concepts- was developed by Leibniz during 1679-90 (Kneale and Kneale
1966, p. 337). As may be gathered from a passage of the same reference (pp.
340-3), or from a translation of an original text of Leibniz (Lewis 1960,
pp. 297-305), this calculus is the theory of operational semilattices (OSL
for short) with applications to concepts; commutativity and idempotence are
explicitly stated, while associativity is implicitly taken for granted (the
aforementioned passage is abbreviated with some slight changes from the
aforementioned translation (Kneale and Kneale 1966, p. 343); notice that the
edition of Lewis' book referred to in Kneale and Kneale (1966) is earlier
than the one referred to above).

Kneale and Kneale (1966)'s assessment of this calculus is unfavorable. It
asserts (p. 337) that Leibniz ``intended, no doubt, to produce something
wider than traditional logic. [...]. But [...] he never succeeded in
producing a calculus which covered even the whole theory of syllogism.''. On
p. 345 this assertion is elaborated ``What he [Leibniz] produced was
certainly much less than he hoped to produce. For the last scheme [the
calculus de continentibus et contentis], lacking as it does any provision
for negation or for consideration of conjunction and disjunction together,
is still a fragment. So far from including all Aristotle's syllogistic
theory as a part, it contains no principle of syllogism except the first
[...]''.

Likewise, Lenzen (2004)'s assessment of the calculus de continentibus et
contentis is unfavorable. It asserts (p. 28) that this calculus ``remains a
very weak and uninteresting system [...]; thus it shall no longer be
considered here.''.

On the contrary, we have shown that neither negation (of terms) nor any
additional operations are needed to algebraically interpret AAS. It suffices
to require the OSL to possess an annihilator, i.e. to be OSLA. For the
structures based on the OSLA are the TCLSA and, by corollary 16.5, AAS is
both sound and complete with respect to them.

According to Kneale and Kneale (1966, p. 339) it may be seen that Leibniz
practically introduced annihilators when he interpreted $Eab$ as $ab$ ($\mu
a*\mu b$, in our terminology) is nothing.

Lenzen (2004) goes even further. It (pp. 2-3) asserts that Leibniz developed
stronger calculi, the most important of them (p.3) ``is L1, the full \textit{%
algebra of concepts} [...], L1 is deductively equivalent or isomorphic to
the ordinary algebra of sets. Since Leibniz happened to provide a complete
set of axioms for L1, he ``discovered'' the Boolean algebra 160 years before
Boole.''.

Moreover, Lenzen (2004) asserts that Leibniz succeeded in making use of his
logical theory to derive the basic laws of Aristotelian syllogism (p. 55).
In particular, the Aristotelian inferences may be derived as theorems of L1,
or the stronger calculus L2 (p. 56); a detailed discussion of the subject
may be found in Lenzen (2004, \S 8, pp. 55-73). Indeed, as we have shown,
AAS does not need all of this.

Boole did more than just algebraically interpreting AAS. In addition to
annihilators, which are sufficient for dealing with Aristotelian syllogisms
(which involve no term negation), he introduced complementation (which
corresponds to term negation) and a second binary operation. This is
possibly to:

1. be able to interpret all the Aristotelian categorical sentences into
equations (cf. Boole 1948, pp. 20-5),

2. deal with medieval categorical sentences which may involve term negation
(cf. Boole 1948, pp. 20, 27-47), or

3. deal with hypotheticals (cf. Boole 1948, pp. 48-59).

In addition to establishing the Aristotelian syllogistic rules, Boole (1948)
established some non-Aristotelian ones. For example (p. 37) $\frac{Ezy,Oyx}{%
Ox^{\prime }z}$, where ``$x^{\prime }$'' denotes ``not-$x$''.

Boole (1948) did not address the question of completeness, neither did he
consider consequences of more than two premises. However, it discussed (pp.
76-81) a general scheme to solve arbitrarily finite systems of simultaneous
equations in arbitrarily finitely many variables; applying, in particular,
Lagrange's method of indeterminate multipliers. This discussion took place
after making (p. 18) the confounding assertion ``[...] all the processes of
common algebra are applicable to the present [Boolean] system.''.

For one more confounding assertion see below.\newline

\textbf{18. Inadequacy: bounds of AAs. }Calling the symbols of its system
``elective symbols'' (p. 16), Boole (1948) makes (p. 59) another
confounding assertion: ``Every Proposition which language can express may be
represented by elective symbols, and the laws of combination of those
symbols are in all cases the same; but in one class of instances the symbols
have reference to collections of objects, in the other, to the truths of
constituent Proposition.''. This, probably, amounts -in modern language- to
asserting: Every proposition which language can express is equivalent to a
sentential combination of categorical sentences (SCCS for short).

SCCS should be taken seriously, since a stronger assertion has dominated
human thought over more than two millennia: Every argument can be put in a
syllogistic form. Even Bertrand Russell (1967, p. 198) asserts ``Of course
it would be possible to re-write mathematical arguments in syllogistic form,
but this would be very artificial and would not make them any more cogent.''.

Concerning these assertions, it is worthwile to bring to the fore what
Boche\'{n}ski (1968) calls attention to. On p. 63 it observes that
Artistotle ``says explicitly that not all logical entailment is
``Syllogistic''.''. Moreover it observes on the same page that
Aristotle declares that some logical entailments \textit{cannot} be reduced
to syllogisms. So it may be concluded that Artistotle himself contradicts
the aforementioned assertions of Boole and Russell, which makes making them
deeply confounding, and makes it more urgent for historians of thought to
investigate the matter.

Understanding SCCS depends on understanding the notion of categorical
sentences. If term negation is permitted, the sentences will be called
``Boolean categorical sentences'' and the corresponding assertion will be
denoted by ``SCBCS''. Otherwise, the sentences will be called ``Aristotelian
categorical sentences'' and the corresponding assertion will be denoted by
``SCACS''.

Hilbert and Ackermann (1950) formalizes the Boolean categorical sentences
(pp. 44-8) and informally refutes SCBCS (pp. 55-6).

To formally discuss SCACS (making use only of the methods developed above
and the well known results of sentential logic) augment the alphabet of the
language $C$ of the natural deduction formalization defined in section 1.4
above, by a ternary relation symbol $E^{\prime }$, and add $E^{\prime }abc$ (%
$a,b,c\in C$) to the set of sentences based on $C$. Denote the new set of
sentences by ``$BN^{\prime }(C)$''.

Intuitively, we like $E^{\prime }abc$ to mean that no a which is $b$, is $c$. This may be formalized as follows:

Interpret $BN^{\prime }(C)$ in a WLS $\mathfrak{B}^{\mu }=<B,*,\mu >$ by adding
the following provision to the provisions of definition 14.2.

\begin{tabular}{ll}
5. $\mathfrak{B}^{\mu }\vDash E^{\prime }abc$ & iff the system of equations: $%
x*\mu a=x$, \\ 
& $x*\mu b=x$ and $x*\mu c=x$ has no solution.
\end{tabular}
\qquad \newline

In an AWLS $\mathfrak{B}^{\mu }=<B,*,0,\mu >$, $BN^{\prime }(C)$ is interpreted
by adding the following provision to the provisions of part 2 of definitions
16.1.

\begin{tabular}{ll}
5. $\mathfrak{B}^{\mu }\vDash E^{\prime }abc$ & iff $0$ is the only solution of
the system \\ 
& of equations: $x*\mu a=x$, $x*\mu b=x$ \\ 
& and $x*\mu c=x$.
\end{tabular}
\newline

Recall that $0$ is not in the range of $\mu $; also notice that if $\mathfrak{B}%
^{\mu }$ is a TCLSA, then this provision is equivalent to

5$^{\prime }$. $\mu a*\mu b*\mu c=0$.

The other syntactical and semantical notions remain the same, or to be
appropriately modified in the obvious way.\newline

Let $\Gamma _{0},\Gamma _{1}\subseteq BN^{\prime }(C)$ and let $\mathfrak{D}%
\subseteq WLS\cup AWLS$. $\Gamma _{0}$ is said to $\mathfrak{D}$-imply $\Gamma
_{1}$ (symbolically $\Gamma _{0}\underset{\mathfrak{D}}{\vDash }\Gamma _{1}$)
if for every $D\in \mathfrak{D}$, $D\vDash \Gamma _{1}$ whenever $D\vDash \Gamma
_{0}$. $\Gamma _{0}$ is said to be $\mathfrak{D}$-equivalent to $\Gamma _{1}$,
or $\Gamma _{0},\Gamma _{1}$ are $\mathfrak{D}$-equivalent, if each of them $%
\mathfrak{D}$-implies the other. $\Gamma _{0}$ is said to be $\mathfrak{D}$-valid if 
$\phi \underset{\mathfrak{D}}{\vDash }\Gamma _{0}$, it is said to be $\mathfrak{D}$%
-consistent if $D\vDash \Gamma _{0}$ for some $D\in \mathfrak{D}$.

The above notions may be generalized, in the obvious way, to sets of
sentential combinations of elements of $BN^{\prime }(C)$. If $\Gamma _{0}$
or $\Gamma _{1}$ is a singleton, it may be replaced by its unique element,
e.g. ``$\rho \underset{\mathfrak{D}}{\vDash }\sigma $'' may replace ``$\{\rho
\}\underset{\mathfrak{D}}{\vDash }\{\sigma \}$''.

In what follows $c_{0},c_{1}$ and $c_{2}$ are assumed to be pairwise
distinct elements of $C$. For every $\mathfrak{D}\subseteq WLS\cup AWLS$, $%
Ec_{0}c_{1}$ $\mathfrak{D}$-implies $E^{\prime }c_{0}c_{1}c_{2}$. The converse
depends on $\mathfrak{D}$. In particular it does not hold for $\mathfrak{D}= BLSA$.
As a matter of fact we have the following:\newline

\noindent THEOREM 18.1.\textbf{\ }Let $\sigma $ be a sentential combination
of elements of $BN(C)$, then:

1. $E^{\prime }c_{0}c_{1}c_{2}$ is not $BLSA$-equivalent to $\sigma $, hence

2. $E^{\prime }c_{0}c_{1}c_{2}$ is not deductively equivalent to $\sigma $
(i.e. one of them does not deductively entail the other), for each deductive
system which is sound with respect to $BLSA$.\newline

To prove this, we first prove:\newline

\noindent LEMMA 18.2.\textbf{\ }Put:

$\Gamma _{0}=\{Ic_{0}c_{1},Ic_{1}c_{2},Ic_{2}c_{0}\}$ \quad and $\quad
\Gamma _{1}=\Gamma _{0}\cup \{E^{\prime }c_{0}c_{1}c_{2}\}$

then:

1. $\Gamma _{1}$ is $BLSA$-consistent.

2. $\Gamma _{1}$ is not $BLSA$-implied by any $BLSA$-consistent $\Gamma
\subseteq BN(C)$.\newline

\noindent \textit{Proof.} Part 1 is easy. To see part 2, assume that there is
a subset $\Gamma \subseteq BN(C)$ which is both $BLSA$-consistent and
$BLSA$-implies $\Gamma _{1}$. Then there is $\mathfrak{B}^{\mu }\in  BLSA$ which
is a model of $\Gamma \cup \Gamma _{1}$. By theorem 16.4 it may be assumed
that $\mathfrak{B}^{\mu }=<{\wp }(B),\cap ,\phi ,\mu >$ for some $B$.

Let $B^{\prime }=B\cup \{a\}$ for some $a\notin B$ and let $\mathfrak{B}^{\prime
\mu ^{\prime }}=<{\wp }(B^{\prime }),\cap ,\phi ,\mu ^{\prime }>$ where
for every $c\in C$,
\begin{eqnarray*}
\mu ^{\prime }(c)= 
\begin{cases}
\mu (c)\cup \{a\} & \text{if} \quad \mathfrak{B}^{\mu }\vDash Ac_{i}c \quad \text{for some} \quad i\in 3,
\\ 
\mu (c) & \text{otherwise}.
\end{cases}
\end{eqnarray*}

$\mathfrak{B}^{\prime \mu ^{\prime }}$ is a $BLSA$ which is basically equivalent
to $\mathfrak{B}^{\mu }$, hence it is a model of $\Gamma$; but it is not a
model of $\Gamma _{1}$. From this the required follows. \hspace{2.5cm} $\square $
\newline

\noindent \textit{Proof of theorem 18.1.} Assume that $E^{\prime
}c_{0}c_{1}c_{2}$ is $BLSA$-equivalent to a sentential combination of elements
of $BN(C)$, $\sigma $ say. Then $E^{\prime }c_{0}c_{1}c_{2}\wedge
Ic_{0}c_{1}\wedge Ic_{1}c_{2}\wedge Ic_{2}c_{0}$ ($\rho $ for short) is
$BLSA$-equivalent to $\sigma \wedge Ic_{0}c_{1}\wedge Ic_{1}c_{2}\wedge
Ic_{2}c_{0}$ ($\sigma _{1}$ for short) which also is a sentential
combination of elements of $BN(C)$.

By sentential logic, $\sigma _{1}$ may be assumed to be a disjunction of
conjunctions of elements of $BN(C)$ and their negations. Since $\rho $ is
$BLSA$-consistent and the negation of any element of $BN(C)$ is
$BLSA$-equivalent to some element of $BN(C)$, $\sigma _{1}$ may further be
assumed to be a non-empty disjunction of $BLSA$-consistent conjunctions of
elements of $BN(C)$. Consequently $\rho $ is $BLSA$-implied by each of these
conjunctions, which contradicts part 2 of lemma 18.2. From this the required
follows.$ \hspace{5.5cm} \square $
\newline

\textbf{Acknowledgements. }Several friends were kind enough to provide me
with references which proved to be very helpful. My deep gratitude is hereby
expressed to each of them: Wafik Lotfalla, Essawy Amasha, Sharon Amasha, and
Fawzy Hegab. I am most indebted to two more friends: Azza Khalifa for pointing out some misprints, and Ahmed Ghaleb for patiently and carefully proofreading the manuscript and transforming its scientific Workplace file into TEX.
\newline

\textbf{Appendix. }The following algorithm, to generate the first $n(>0)$
primes, may not be efficient, but it is simple, and its running time (see
below) makes it sufficient for our purposes.\newline

Input: $n$ (positive integer)

Output: $p$ (the strictly increasing list of the first $n$ primes)

Procedure:

Declare $i,j,k,m$ natural number parameters;

$p_{0}\leftarrow 2$;

If $n=1$ go to $***$

Else $p_{1}\leftarrow 3$, $i\leftarrow 1$, $m\leftarrow 2$

End If;

For$_{1}$ $i<n-1$ do

\qquad $k\leftarrow p_{i}+2,m\leftarrow mp_{i}$

For$_{2}$ $k\leq m+1$ do

\qquad $j\leftarrow 0$

For$_{3}$ $j\leq i$ do

\qquad If $p_{j}|k$ go to $*$

\qquad Else $j\leftarrow j+1$

\qquad End If;

Repeat

End For$_{3}$;

$*$ If $j>i$ go to $**$

\qquad Else $k\leftarrow k+2$

\qquad End If;

\qquad Repeat

\qquad End For$_{2}$;

$**$ $i\leftarrow i+1$

\qquad $p_{i}\leftarrow k$

\qquad Repeat

\qquad End For$_{1}$;

$***$ Print $p$;

\qquad End Algorithm.

The termination of this algorithm is guaranteed by the respective upper
bounds stipulated at the beginnings of the three For loops. The correctness
is guaranteed by the well known fact which goes back to Euclid's Elements: $%
p_{i+1}\leq 1+\stackrel{i}{\underset{j=0}{\Pi }}p_{j}$, together with the
simple fact that $p_{i+1}$ is the first (odd) integer greater than $p_{i}$,
which is not a multiple of any of $p_{0},...,p_{i}$.

To estimate the running time, notice that (Landau 1958, p. 91) for large $n$%
, $p_{n}<n^{2}$. For such $n$ the For$_{1}$ loop is iterated at most $n$
times, for each iteration the For$_{2}$ loop is iterated at most $n^{2}$
times, and for each of these iterations the For$_{3}$ loop is iterated at
most $n$ times. All the steps of the algorithm are simple assignment or
comparison steps, the only exception is the test $p_{i}|k$ which needs at
most $k$ ($\leq n^{2}$) simple steps. So the total running time is a
polynomial in $n$, of degree at most $1+2+1+2=6$.\newline

Department of Mathematics

Faculty of Science

Cairo University

Giza - Egypt

m0amer@hotmail.com \quad , \quad amer@sci.cu.edu.eg

URL: http://scholar.cu.edu.eg/?q=mohamedamer/


\begin{thebibliography}{99}
\bibitem{}  Adzic, M. and Dosen, K. (2016). G\"{o}del's Notre Dame course. 
\textit{The Bulletin of Symbolic Logic}, \textbf{22}(4), 469-481.

\bibitem{}  Bellucci, F. and Pietarinen, A. (2016). Existential graphs as an
Instrument of Logical analysis: part I. alpha. \textit{The Review of
Symbolic Logic}, \textbf{9}(2), 209-237.

\bibitem{}  Boche\'{n}ski, I.M. (1968). \textit{Ancient Formal Logic}.
Amsterdam: North Holland.

\bibitem{}  Boger, G. (1998). Completion, reduction and analysis: three
proof-theoretic processes in Aristotle's Prior Analytics. \textit{History
and Philosophy of Logic}, \textbf{19}, 187-226.

\bibitem{}  Boole, G. (1948). \textit{The Mathematical Analysis of Logic}.
Oxford: Basil Blackwell.

\bibitem{}  Corcoran, J. (1972). Completeness of an ancient logic. \textit{%
The Journal of Symbolic Logic}, \textbf{37}, 696-702.

\bibitem{}  Frege, G. (1967). Begriffsschrift, a formula language modeled
upon that of arithmetic for pure thought. In vanHeijenoort, J. ed., \textit{%
From Frege to G\"{o}del: A Source Book in Mathematical Logic, 1879-1931}.
Harvard University Press, Cambridge, MA, 1-82.

\bibitem{}  Glashoff, K. (2002). On Leibniz's characteristic numbers. 
\textit{Studia Leibnitiana}, \textbf{XXXIV/2}, 161-184.

\bibitem{}  --------------- (2005). Aristotelian syntax from a
computational-combinatorial point of view. \textit{Journal of Logic and
Computation}, \textbf{15}, 949-973.

\bibitem{}  --------------- (2007). On negation in Leibniz system of
characteristic numbers. na (to appear).

\bibitem{}  --------------- (2010). An intensional Leibniz semantics for
Aristotelian logic. \textit{The Review of Symbolic Logic}, \textbf{3}(2),
262-272.

\bibitem{}  Hilbert, D. and Ackermann, W. (1950). \textit{Principles of
Mathematical Logic}, 2$^{nd}$ ed. (Translated by Hammond, L.M., Leckie,
G.G., Steinhardt, F. and edited with notes by Luce, R.E.). Chelsea
Publishing Company, New York.

\bibitem{}  Hurley, P.J. (1982). \textit{A Concise Introduction to Logic}.
Belmont, CA: Wadsworth Publishing Company.

\bibitem{}  Kleene, S.C. (1967). \textit{Mathematical Logic}. New York: John
Wiley \& Sons, Inc.

\bibitem{}  Kneale, W. and Kneale, M. (1966). \textit{The Development of
Logic} (third impression). London: Oxford University Press.

\bibitem{}  Landau, E. (1958). \textit{Elementary Number Theory} (Goodman,
J.E., translator). New York: Chelsea Publishing Company.

\bibitem{}  Lenzen, W. (2004). Leibniz's logic. In \textit{Handbook of
the History of Logic}, vol. 3, \textit{The Rise of Modern Logic: From
Leibniz to Frege}; Gabbay, D.M. and Woods, J. editors. Elsevier, Noth
Holland, 1-83.

\bibitem{}  Lewis, C.I. (1960). \textit{A Survey of Symbolic Logic}. New
York: Dover Publications, Inc.

\bibitem{}  \L ukasiewicz, J. (1998). \textit{Aristotle's Syllogistic From
The Standpoint of Modern Formal Logic} (second edition, enlarged). Oxford
University Press Inc., New York (Special edition for Sandpiper Books Ltd.).

\bibitem{}  Martin, J.N. (1997). Aristotle's natural deduction reconsidered. 
\textit{History and Philosophy of Logic}, \textbf{18}, 1-15.

\bibitem{}  Pasquali, F. and Retor\'{e}, C. (2016). Aristotle's square of
opposition in the light of Hilbert's epsilon and tau quantifiers. \textit{%
Arxiv: 1606.08326 [pdf, ps, other], subjects: Logic (math. LO)} work
presented at: Aristotle 2400 years world congress, Thessaloniki, May 23-28,
2016.).

\bibitem{}  Rosenthal, M. and Yudin, P., editors (1967). \textit{A
Dictionary of Philosophy} (English translation, edited by Dixon, R.R., and
Saifulin, M.). Moscow: Progress Publishers.

\bibitem{}  Russell, B. (1967). \textit{A History of Western Philosophy} (a
Clarion book). New York: Simon and Schuster.

\bibitem{}  Russinoff, I.S. (1999). The syllogism's final solution. \textit{%
The Bulletin of Symbolic Logic}, \textbf{5}(4), 451-469.

\bibitem{}  Schumann, A. (2006). A lattice for the language of Aristotle's
syllogistic and a lattice for the language of Vasilev's syllogistic. \textit{%
Logic and Logical Philosophy}, \textbf{15}(1), 17-37.

\bibitem{}  Shepherdson, J.C. (1956). On the interpretation of Aristotelian
syllogistic. \textit{The Journal of Symbolic Logic}, \textbf{21}(2), 137-147.

\bibitem{}  Smiley, T. (1962). Syllogism and quantification. \textit{The
Journal of Symbolic Logic}, \textbf{27}(1), 58-72.

\bibitem{}  Smiley, T.J. (1973). What is a syllogism? \textit{Journal of
Phylosophical Logic}, \textbf{2}, 136-154.

\bibitem{}  Smith, R. (1983). Completeness of an ecthetic syllogistic. 
\textit{Notre Dame Journal of Formal Logic}, \textbf{24}, 224-232.

\bibitem{}  Sotirov, V. (1999). Arithmetizations of syllogistic $\grave{a}$ 
la Leibniz. \textit{Journal of Applied Non-classical Logics}, \textbf{9}(2),
387-405.

\bibitem{}  --------------- (2015). Leibniz arithmetized syllogistic: the
intensional semantics. \textit{Proceedings of the 10}$^{\mathit{th}}$\textit{%
\ Panhellenic Logic Symposium}, 54-58 (Samos, Greece, June 11-15, 2015).

\bibitem{}  Valencia, V. S. (2004). The Algebra of Logic. In 
\textit{Handbook of the History of Logic}, vol.3, \textit{The Rise of Modern
Logic: From Leibniz to Frege}; Gabbay, D.M. and Woods, J. editors. Elsevier,
North Holland, 389-544.

\bibitem{}  Venn, J. (1880). On the diagrammatic and mechanical
representations of propositions and reasoning. \textit{Philosophical Magazine%
}, series 5, 10:59, pp. 1-18, DOI: 10.1080/14786448008626877.
\end{thebibliography}
\end{document}